\documentclass[10pt,twoside]{amsart}

\usepackage[english]{babel}
\usepackage[all]{xy}
\usepackage{amssymb,amsmath,bbm,url,xr}
\usepackage[mathscr]{euscript}
\usepackage{mathrsfs}
\input cyracc.def
\DeclareFontFamily{U}{russian}{}
\DeclareFontShape{U}{russian}{m}{n}
        { <5><6> wncyr5
        <7><8><9> wncyr7
        <10><10.95><12><14.4><17.28><20.74><24.88> wncyr10 }{}
\DeclareSymbolFont{Russian}{U}{russian}{m}{n}
\DeclareSymbolFontAlphabet{\mathcyr}{Russian}
\makeatletter
\let\@math@cyr\mathcyr
\renewcommand{\mathcyr}[1]{\@math@cyr{\cyracc #1}}
\makeatother

\title{Orientation theory in arithmetic geometry}

\author{Fr\'ed\'eric D\'eglise}
\address{IMB (UMR5584), \\
 9 avenue Alain Savary, 2\\
1078 DIJON CEDEX FRANCE
France}
\email{frederic.deglise@u-bourgogne.fr}
\urladdr{http://deglise.perso.math.cnrs.fr/}
\thanks{Partially supported from the ANR (grant No. ANR-12-BS01-0002),
 and from the French "Investissements d'Avenir"
 program, project ISITE-BFC (contract ANR-lS-IDEX-OOOB).}

\date{November 2017.}

\newtheorem{thm}{Theorem}[subsection]
\newtheorem{prop}[thm]{Proposition}
\newtheorem{cor}[thm]{Corollary}
\newtheorem{lm}[thm]{Lemma}
\newtheorem{conj}{Conjecture}

\theoremstyle{definition}
\newtheorem{df}[thm]{Definition}
\newtheorem{ex}[thm]{Example}
\newtheorem{num}[thm]{}

\theoremstyle{remark}
\newtheorem{rem}[thm]{Remark}

\numberwithin{equation}{thm}

\newcommand{\NN} {\mathbb N}
\newcommand{\ZZ} {\mathbb Z}
\newcommand{\QQ} {\mathbb Q}
\renewcommand{\AA} {\mathbb A}
\newcommand{\GG} {\mathbb G_m }
\newcommand{\PP} {\mathbb P}

\newcommand{\T}{\mathscr T} 

\newcommand{\cM}{\mathcal M}
\newcommand{\E}{\mathbb E} 
\newcommand{\F}{\mathbb F} 

\newcommand{\un}{\mathbbm 1} 

\newcommand{\base}{\mathscr S} 

\renewcommand{\H}{\mathscr H_\bullet}
\newcommand{\SH}{S\mathscr{H}}

\newcommand{\sm}{\mathscr Sm}

\newcommand{\DM}{\mathrm{DM}}

\newcommand{\smod}[1]{#1\!-\!mod}
\newcommand{\wmod}[1]{#1\!-\!mod^w}

\newcommand{\ilim}[1] { \underset{#1}{\varinjlim} \ }

\DeclareMathOperator{\Th}{Th}
\DeclareMathOperator{\thom}{\mathfrak t}
\DeclareMathOperator{\rthom}{\bar{\mathfrak t}}
\DeclareMathOperator{\fdl}{\eta}
\DeclareMathOperator{\rfdl}{\bar{\eta}}
\DeclareMathOperator{\td}{Td}

\newcommand{\wdg}{\wedge}
\newcommand{\Sus}{\Sigma^\infty}
\newcommand{\MGL}{\mathbf{MGL}}
\newcommand{\KGL}{\mathbf{KGL}}
\newcommand{\HB}{\mathbf H_{\mathcyr{B}}}

\newcommand{\HBx}[1]{\mathbf H_{\mathcyr{B},#1}}
\newcommand{\HM}[1]{\mathbf H{#1}}
\newcommand{\HMx}[2]{\mathbf H{#1}_{#2}}
\newcommand{\HMS}[1]{\mathbf M{#1}}
\newcommand{\HMSx}[2]{\mathbf M{#1}_{#2}}

\DeclareMathOperator\Pic{Pic}
\newcommand{\spec}[1] {\operatorname{\mathrm{Spec}}(#1)}

\DeclareMathOperator{\Hom}{Hom}
\DeclareMathOperator{\uHom}{\underline{Hom}}

\newcommand{\pur}[1] { \mathfrak p_{(#1)} } 

\newcommand{\HH}{\mathbf{H}}

\newcommand{\et}{\text{\'et}}

\newcommand{\cont}{\text{cont}}
\DeclareMathOperator{\ch} {ch}
\DeclareMathOperator{\reg} {\mathscr Reg}
\DeclareMathOperator\coker{coKer}
\DeclareMathOperator{\res} {Res}

\newcommand{\h}  {\mathrm{h}}

\newcommand{\cupp}{{\scriptstyle \cup}}
\newcommand{\refp}{\cdot}

\newcommand{\eq}[1][r]
   {\ar@<-3pt>@{-}[#1]
    \ar@<-1pt>@{}[#1]|<{}="gauche"
    \ar@<+0pt>@{}[#1]|-{}="milieu"
    \ar@<+1pt>@{}[#1]|>{}="droite"
    \ar@/^2pt/@{-}"gauche";"milieu"
    \ar@/_2pt/@{-}"milieu";"droite"}

\begin{document}

\begin{abstract}
This work is devoted to study orientation theory
 in arithmetic geometry within the motivic homotopy theory
 of Morel and Voevodsky.
 The main tool is a formulation of the absolute purity
 property for an \emph{arithmetic cohomology theory},
 either represented by a cartesian section
  of the stable homotopy category
	or satisfying suitable axioms. 
	We give many examples, formulate conjectures and
  prove a  useful property of analytical invariance.
 Within this axiomatic, we thoroughly develop the theory
 of characteristic and fundamental classes, Gysin and residue
 morphisms. This is used to prove Riemann-Roch formulas,
 in Grothendieck style for arbitrary natural transformations
 of cohomologies, and a new one for residue morphisms.
 They are applied to rational motivic cohomology and \'etale
 rational $\ell$-adic cohomology, as expected by Grothendieck in
 \cite[XIV, 6.1]{SGA6}.
\end{abstract}

\maketitle

\tableofcontents

\section*{Introduction}

\noindent \textbf{History.} One of the most striking intuition
 of Riemann is that the
 natural domain of definition of abelian integrals are (branched) surfaces
 rather than the complex plane.
 Retrospectively, one is amazed that this single idea contained
 in seeds the modern development of both analytical and algebraic
 geometry, whose varieties
 are now studied through their sheaf of functions.
In this long and deep evolution, the Riemann-Roch formula
 played a catalytic role, or rather that of a lighthouse.

In his 1857 masterpiece, \cite{Riemann},
 Riemann studied his new complex functions,
 defined in modern terms on a compact Riemann surface $\Sigma$.
 Notably, he described the general form of these functions once
 we prescribed $m$ given (simple) poles in $\Sigma$.
 The striking new idea is the appearance of a geometrical
 invariant of $S$ that Riemann had discovered before, that we now
 know as the \emph{genus} $p$ of $\Sigma$.
 Riemann established that complex functions on $S$ with $m$
 given poles depend upon (at least) $m-p+1$ constants
 (see \emph{loc. cit.}, \textsection 5). A formula that we
 now read as:
\begin{equation} \tag{R}\label{eq:Rf}
l(D) \geq \deg(D)-p+1,
\end{equation}
where $D$ is the divisor on $S$ made by the formal sum
 of the $m$ given points, with degree $\deg(D)$ equal to $m$,
 and $l(D)$ is the dimension of the space of functions $f$ on
 $\Sigma$ whose associated divisor is $D$
 -- meaning it admits simple poles exactly at the $m$ points
 of the support of $D$.
A few years later, Roch in \cite{Roch}
 interpreted analytically the difference
 of the two members of \eqref{eq:Rf} as the ``number of linearly
 disjoint integrands which can vanish at the $m$ given poles''
 (see \cite[second par. p. 802]{Gray}). In modern terminology,
 this becomes the \emph{Riemann-Roch formula}, 
 which we write today as:
\begin{equation} \tag{RR}\label{eq:RRf}
l(D)-l(K-D)=\deg(D)-p+1
\end{equation}
where $K$ is the canonical divisor on $\Sigma$.

Looking through the glass of a century of research,
 one is amazed by the exceptional role that took up
 this simple formula.
 This is particularly visible in the algebraic reformulation
 of Riemann's ideas by Clebsch,
 and then Brill-Noether (Max), where one of the driving motivation
 was to define the genus of an algebraic curve 
 (complex plane projective)
 in order to prove the Riemann-Roch formula.
 The same problem is addressed slightly later in the development
 of algebraic surfaces by M.~Noether, and then by the Italian
 school (Castelnuovo, Enriques, Severi, ...), whose guide
 was to formulate the correct 
 extension of the Riemann inequality in dimension 2
 and in particular to find the good notion of genus.

But the most historically surprising application
 of the formula came almost eighty 
 years after its introduction by Riemann when F.K.~Schmidt extended
 it to the case of function fields $K$ over a finite field and
 use it to prove the rationality of the Zeta function associated
 with $K$ (1931). The formula followed the path
 opened up by the influential 1882 work of Dedekind-Weber,
 who developed the birational point of view initiated by Riemann
 by transporting his work to the purely arithmetical world of 
 function fields over the complex numbers. 
The impact of Schmidt's proof on the modern formulation
 of algebraic geometry is of primary importance,
 as it lead Weil to his work on abelian varieties and
 most of all to the formulation of his conjectures on Zeta functions.

As said by Dieudonn\'e in his history of algebraic geometry,
 the followers of Riemann split into several branches without
 much interactions (see \cite{Dieu}, beginning of chap. VI).
 So while the notion of cohomology in arithmetic geometry was
 slowly revealing itself, the algebraic geometers working
 on the Riemann-Roch problem for surfaces were discovering
 the theory of \emph{canonical classes}: M. Noether for surfaces (1886),
 Severi (1932), Segre and Todd in higher dimensions.
 Meanwhile, Poincar\'e introduced singular homology
 and topologists started to study \emph{characteristic classes}
 of vector bundles (Stiefel and Whitney 1935, Chern 1946)
 without any connections with the theory of canonical classes.
The unifying tool was to be the theory of sheaves invented by Leray
 during World War II. Only a few years after its introduction,
 this theory was fully developed first by Cartan and Serre for
 analytical varieties and secondly by Kodaira and Spencer for
 Kahlerian varieties. The problem of extending the Riemann-Roch
 formula in higher dimensions
 was crystallized in those years of boiling development around
 the notion of sheaves, in particular through the attempts
 of Kodaira and Serre. The first one had already solved the
 extension problem for Kahlerian varieties of dimension 2
 (1951) and dimension 3 (1952) and linked the problem with
 computations of the canonical classes of Todd, while Serre
 had remarked its link with duality, used Thom cobordism
 theory to treat a special case\footnote{This work is unpublished
 but see the account of \cite[VIII. 12.]{Dieu}.}
 and conjectured a general form for the extended Riemann-Roch formula.
 Shortly after these advances, it belonged to Hirzebruch 
 (1954, see also \cite{Hirz}) to prove (and make precise\footnote{In a letter
 to Kodaira and Spencer, Serre conjectured that the Euler characteristic
 should be expressed by some polynomial expression on Chern classes})
 the formula conjectured by Serre,
 formula that we now call the Hirzebruch-Riemann-Roch formula:
\begin{equation} \tag{HRR}\label{eq:HRRf}
\chi(X,E)=\deg\big(\ch(E).\td(E)\big)
\end{equation}
where $E$ is a vector bundle on an analytical variety $X$,
 $\chi(X,E)$ (resp. $\ch(E)$, $\td(E)$)
 its Euler characteristic (resp. Chern character, Todd class)
 with values in singular cohomology.
 The proof of Hirzebruch, rather technical, uses (and developed)
 the theory of characteristic classes (Chern, Todd classes,...)
 and makes use again of Thom cobordism theory.\footnote{To anticipate
 the content of this paper, one remarks that the universality of
 cobordism theory was already fully playing its role here.}

But the final revolution of Riemann's original problem
 was imagined by a single man whose ideas were to change
 completely our conception of it, Grothendieck.
 Shortly after the proof of Hirzebruch (see \cite{BS}),
 Grothendieck gave a new and meaningful interpretation
 of the formula, whose first practical interest was
 the simplicity of its proof and its validity for algebraic
 varieties over an arbitrary base field.
 The two main ideas introduced there,
 which had never been anticipated before,
 was first a relative formulation
 (\emph{i.e.} for a morphism rather than a single algebraic varieties)
 and secondly a purely cohomological interpretation 
 of the \eqref{eq:RRf}-formula by
 introducing a generalized cohomology theory that would soon
 become famous as $K$-theory. From a conceptual point of view,
 the Grothendieck-Riemann-Roch formula expresses
 the defect of functoriality of a natural transformation of cohomology
 theories with respect to the exceptional covariant functoriality;
 in the assumptions of Borel-Serre, given a proper morphism 
 $f:Y \rightarrow X$ of non singular quasi-projective varieties
 over any field $k$,
 $T_X$ (resp. $T_Y$) being the tangent bundle of $X$ (resp. $Y$),
 one has for any element $y$, in the K-group $K(Y)$ of 
 \emph{virtual vector bundles} over $Y$,
\begin{equation} \tag{GRR}\label{eq:GRRf}
\ch(f_*(y)).\td(T_X)=f_*\big(\ch(y).\td(T_Y))
\end{equation}
where $\ch$ denotes the Chern character from K-theory 
 to rational Chow groups, and $\td$ is the Todd class of
 a vector bundle.

During the period just described, history tells us that
 interactions between topology, geometry and algebraic geometry
 were very strong\footnote{It was after a seminar in Princeton
   which gathered most of the main characters discussed here
	 that Hirzebruch found his proof.}.
 Therefore, soon after the appearance of the \eqref{eq:GRRf}-formula,
 Atiyah and Hirzebruch introduced topological K-theory,
 that they immediately understood as a generalized cohomology
 theory\footnote{\emph{i.e.} a cohomology theory that satisfies
 all the axiom of Eilenberg-Steenrod except the dimension axiom;}
 and proved the topological formulation of the \eqref{eq:GRRf}-formula.
 It was soon realized that the covariant functoriality
 involved in the formula should be a consequence of Poincar\'e duality,
 on the model of the covariant functoriality discovered by Gysin
 in his study of sphere bundles (1942).
Consequently, a very general \eqref{eq:GRRf}-formula, in which one considers
 an arbitrary natural transformation of generalized cohomology theories,
 each equipped with a \emph{complex orientation} to get the usual theory
 of characteristic classes, was written by Dyer (see \cite{Dyer})
 -- and stated as a folklore theorem,
 only 4 years after the original formulation of Grothendieck !

History must stop at some point. We will end it by two cornerstones
 of which our work is a direct continuation. The first one is
 Quillen discover of the universality of Thom 
 complex cobordism theory in terms of formal group laws
 and oriented cohomology theories:  \cite{Qui}.
 The second one is Grothendieck's final extension
 of his Riemann-Roch formula to the arithmetic setting in
 \cite{SGA6}.

\bigskip

\noindent \textbf{Motivic stable homotopy and cohomology theories.} The
 purpose of this work is to extend the arithmetic formulation
 of the Grothendieck-Riemann-Roch formula of \cite{SGA6} 
 in the same way that by Dyer (again \cite{Dyer})
 extends the topological formulation of the Hirzebruch-Riemann-Roch formula
 following Atiyah-Hirzebruch.
 To that end, the natural framework is Morel-Voevodsky's
 stable motivic homotopy theory, as it is defined by a clear analogy
 with the ordinary stable homotopy category of topological spaces
 used by Dyer.\footnote{When we work over the field of complex numbers,
 the stable motivic homotopy category can be realized in the
 ordinary stable homotopy category so that any motivic construction
 or statement has a realization in the topological world.}

The objects of the stable homotopy category,
 both classical and motivic,
 are meant to represent cohomology theories.
 Called spectra, they form a triangulated category whose
  distinguished triangles correspond 
  to universal long exact sequences in cohomology.
 Similarly, all structures or properties
  of the stable homotopy category are reflected
  in the cohomologies representable by spectra.
 Probably the most important example of such a structure
  is the existence of
  a (symmetric) tensor product, called the \emph{smash product}:
  a (commutative) monoid\footnote{In this introduction and in the whole paper,
   all monoid structures will be assumed to be commutative.}
   on a spectra induces a product structure on its cohomology.
  These monoids are of primary importance
  in (motivic) stable homotopy; they are called \emph{(motivic) ring spectra}.

In the motivic setting, we work over a base scheme\footnote{In this introduction
 and in the whole paper, all schemes will be assumed to be Noetherian
 of finite dimension.} $S$; the motivic stable homotopy category
 of Morel-Voevodsky is denoted by $\SH(S)$.
 The starting point of this paper is that the representability of
  a cohomology theory has many interesting consequences.
 Let us first describe the obvious ones,
  for a given spectrum $\E$ 
  (see Prop. \ref{prop:ppty_product} for details):
  the cohomology represented by $\E$
  is a contravariant functor $\E^{**}$ from smooth $S$-schemes to bigraded
  abelian groups (the first index is the degree and the second one
  is called the \emph{twist}). It satisfies the homotopy invariance property
  with respect to the affine line $\AA^1$
  (the affine line $\AA^1$ is contractible),
  stability property with respect to the projective line $\PP^1$
  (seen as the analogue of the circle in topology).
 Moreover, it can be extended to a \emph{cohomology theory with support}:
  given a smooth $S$-scheme $X$ and a closed subscheme $Z \subset X$,
  one can define a bigraded abelian group $\E^{**}_Z(X)$ of cohomology
  classes with support in $Z$, in an appropriate functorial way
  and such that $\E^{**}_X(X)=\E^{**}(X)$ (see Def. \ref{df:coh_support} for details).
 This theory with support satisfies the (Nisnevich) excision property
  (analogue of the excision property of the Eilenberg-Steenrod axioms in topology,
   see Sec. \ref{sec:anaylitical_inv}, property (Nis) for details).

In this paper, we will use another important property
 of $\SH(S)$, its basic functoriality in $S$:
 it is a fibered category over the category of schemes.
 A cartesian section\footnote{Explicitly: the data of a ring spectrum $\E_X$
  for each scheme $X$,
  with a given transition isomorphism $f^*(\E_X) \simeq \E_Y$ for any morphism
  $f:Y \rightarrow X$.} with respect to this fibered structure,
  eventually restricted to a subcategory $\base$ of the category of schemes,
  will be called an \emph{absolute spectrum}
   (see Def. \ref{df:abs_ring_sp} and Ex. \ref{ex:abs_ring_spectra}
    for examples).
 Such an absolute spectrum represents a cohomology 
  which is defined over the whole category $\base$,
  allowing to avoid the restriction to smooth $S$-schemes.\footnote{Using
  a terminology introduced by Beilinson, this is an absolute cohomology;
  this justifies the terminology \emph{absolute spectrum}.}
 It still satisfies all the properties enumerated above.
  But moreover, under the presence of a ring structure
   on the absolute spectrum, we get an important product on cohomology
   with support which is not commonly used (but see \cite[IV]{SGA45}). We call it
   the \emph{refined product}: given closed subschemes $T \subset Z \subset X$,
   it has the form
\begin{equation} \tag{$*$}\label{eq:product1_intro}
\E^{**}_T(Z) \otimes \E^{**}_Z(X) \rightarrow \E^{**}_T(X)
\end{equation}
and will be an essential technical tool to our study of fundamental classes
 (see Par. \ref{num:products}).

Note that all the basic properties of the cohomologies representable
 by an absolute spectrum are gathered in Proposition \ref{prop:ppty_product}.

\bigskip

\noindent \textbf{Absolute purity.} The absolute purity conjecture
 of Grothendieck,
 formulated in \cite[I 3.1.4]{SGA5},
 has been a major problem in \'etale cohomology,
 because of its consequences on finiteness and duality
 as pointed out by Grothendieck (see \emph{loc. cit.}).
 It was solved by Thomason in \cite{Thomason_pur} under some assumptions
  (on the coefficients) and in full generality by Gabber in \cite{Fuj}.
 Roughly, the conjecture says that the cohomology of a regular
  scheme $X$ with support in a closed regular scheme $Z$ is isomorphic to
  the cohomology of $Z$.
 The case of smooth schemes over a field (or even over some base)
  can be treated easily (see \cite[XVI, 3.9]{SGA4}). The main problem
  in this conjecture is to treat the case of regular schemes of unequal
  characteristics, which are the objects of the so-called
  \emph{arithmetic geometry}.

This problem was confined in the \'etale setting until D.C.~Cisinski and
 the author discovered that a similar statement could be formulated
 and proved in the newly defined setting of rational mixed motives
 (see \cite[14.4.1]{CD3})\footnote{As in the \'etale setting, 
 this result has important consequences for rational mixed motives.
 Note the importance of absolute purity was anticipated by Ayoub
 in \cite{Ayoub1}.}. This naturally raises the question of extending
 the problem of absolute purity
 to any representable cohomology theories.

So we introduce here (Def. \ref{df:abs_pur}) 
 the property of \emph{absolute purity} for an absolute
 spectrum and any closed subscheme $Z \subset X$
 whose immersion is regular.
 The formulation of this property is
 the new technical ingredient introduced by this work.
 As in the case of \'etale cohomology, for smooth schemes 
 over a field (or over some base), the property is always fulfilled
 according to a fundamental result of Morel-Voevodsky.
 Thus the interest of this property lies in the case of arithmetic
  geometry.
 Fortunately, there are several cohomology theories which satisfies
  absolute purity in the arithmetic case. The matrix case is
  integral algebraic K-theory according to the localization theorem
  of Quillen. The cases of rational motivic cohomology and rational
  cobordism theory follows (see Ex. \ref{ex:abs_pure}).

We think that the problematic of absolute purity is
 an important question for homotopy theory. A new aspect of
 our definition of this property is that it is intrinsic and
 do not depend on the choice of a purity isomorphism.\footnote{This
 is due to the use of the deformation space.} Moreover, it is formulated
 for any cohomology theory without assuming the existence of an
 orientation (see below). In particular, we conjecture that this
 property holds integrally not only for the algebraic cobordism
 spectrum but also for the sphere spectrum
 (see conjectures \ref{conj:purityMGL} and \ref{conj:purityS0}
  p. \pageref{conj:purityMGL}).
 This result would have several interesting
 consequences: see Remark \ref{rem:conseq_abs_pur}.

An important ingredient in the proof of absolute purity
 by Gabber is the so-called analytical invariance of
 \'etale cohomology with support. While we do not attack
 the previous conjectures on absolute purity, we nevertheless
 prove the \emph{analytical invariance property} for any cohomology
 representable by an absolute spectrum (see Thm. \ref{thm:analytic}),
 extending a result already obtained by Wildeshaus \cite{Wil}.
 Note in particular that the result can be applied to the
 absolute spectrum representing \emph{rigid cohomology}
 over a field
 of characteristic $p>0$ -- according to Ex. \ref{ex:abs_pure}(1).
 Thus our proposition contains
 in particular Theorem 1.1 of \cite{Ouw}.

\bigskip

\noindent \textbf{Orientation theory: characteristic and fundamental classes}
At this point,
 we connect motivic homotopy theory with the two fundamental
 notions that were developed has a natural evolution of the Riemann-Roch
 problem: characteristic classes and fundamental classes.
Recall the first ones have first been studied in algebraic topology,
 and evolved naturally into the theory of orientation while the second one
 was the domain of algebraic geometry and evolved into intersection
 theory.

The idea of transporting orientation theory from algebraic topology
 to algebraic geometry is comparatively quite recent
 as it takes its origin in the first proof of the Milnor
 conjecture by Voevodsky (see \cite{Voe0}).
 The theory has grown out of an unpublished work of Morel
 which was developed by several authors.

Let us recall in this introduction the basics of this theory.
 An orientation $c$ of a ring spectrum $\E$ over a base scheme $S$
 is a cohomological class -- in degree $(2,1)$ -- of
 the infinite projective space with coefficients in $\E$
 (Def. \ref{df:orientation}).
 Giving such a class allows to derive \emph{canonically}
 a lot of interesting cohomological invariants and structures.
 First, we can define Chern classes associated with vector
 bundles -- the first Chern class follows directly from the orientation,
 and the other ones follow from the computation of the cohomology
 of any projective bundle according to the classical construction
 of Grothendieck: see Section \ref{sec:orientation&Chern}.

The connection with the work of Quillen appears at this point.
 As in topology, Chern classes of an arbitrary oriented spectrum
 need not be additive: the first Chern class of a tensor product
 of line bundles is not the sum of the Chern classes of each
 bundle. Instead, it is described according to a \emph{formal group
 law} that is canonically associated with the chosen orientation
 (see Par. \ref{num:FGL} for more details).
 This is connected with the theory of \emph{Thom classes}
 that one derives from Chern classes.\footnote{The Thom class
 of a vector bundle $E$, as a class in the projective completion
 of $E$, equals the top Chern class of the universal quotient
 bundle (see Ex. \ref{ex:thom&euler_quotient_bdl}).}
 Indeed, these classes uniquely define a canonical 
 structure of $\MGL$-algebra on $\E$ where $\MGL$ is the algebraic
 cobordism spectrum (analogue of the complex cobordism spectrum
 $\mathbf{MU}$ in topology). Orientations on $\E$ are
 in one to one correspondence with structures of $\MGL$-algebra
 and it is widely believed that the formal group law of $\MGL$
 is the universal formal group law.\footnote{According to works
 respectively of Levine and of Hoyois, given a field $k$
 of exponential characteristic $p$, 
 this is true for $\MGL[1/p]$ if one restricts to $k$-schemes.}

Secondly, the orientation determines fundamental classes.
 In fact, assuming the absolute purity property
 for a regular closed immersion $Z \subset X$,
 the Thom class of the normal bundle of $Z$ in $X$ gives us 
 directly the \emph{refined} fundamental class of $Z$ in $X$,
 as a cohomology class of $X$ with support $Z$.
 Using the refined product
 by this fundamental class 
 -- \eqref{eq:product1_intro} with $T=Z$ --
 gives us the usual purity isomorphism:
$$
\E^{**}(Z) \xrightarrow{\sim} \E^{**}_Z(X).
$$
With our formalism, all this follows easily. But note however
 that this is the first appearance of this form of the purity
 isomorphism in motivic homotopy theory
 -- and even in \'etale cohomology.\footnote{We were especially
 inspired by the formulation of Poincar\'e duality of Bloch and Ogus
 in \cite{BO}.}

These generalized fundamental classes
 are the trace of classical intersection theory,
 though they can be defined in very general theories such as algebraic
 K-theory or algebraic cobordism. We illustrate this concretely
 by proving some of the classical formulas known for Chow groups:
 the excess intersection formula (Cor. \ref{cor:pullback&fdl_class})
 and associativity (compatibility with composition,
 Th. \ref{thm:product_fdl}).
 Note that the excess intersection formula allows to get back
 the classical link between fundamental classes and characteristic
 classes: when $Z$ can be parametrized by a section of a vector bundle
 over $X$, the fundamental class of $Z$ in $X$ equals the top Chern
 class of the vector bundle (see Cor. \ref{cor:fdl&Chern_transversal_section}
 for details).
 Let us also recall that the compatibility with composition
  of fundamental classes is a key technical problem in orientation
	theory. The geometrical tool (the double deformation space)
	used here is not new but the use of the refined
	product \eqref{eq:product1_intro} allows both a finer expression
  of this result and an easier proof.

\bigskip
\noindent \textbf{Residues and Gysin morphisms.}
A formal consequence of the purity isomorphism 
 for a regular closed immersion $i:Z \rightarrow X$
 is the existence of a localization long exact sequence
 in cohomology, so that the theory described just above
 for an oriented absolute ring spectrum canonically
 leads to such a sequence. It is made of two interesting
 morphisms. The first one is a morphism that we have called the 
 \emph{residue map associated with $(X,Z)$}. Our interest to
 that kind of maps come from our comparison between cycle modules
 and homotopy sheaves with transfers (\cite{Deg9}) as our residues
 on sheaves corresponds to Rost residues on cycle modules.
 One illustrates this phenomenon by the computation of
 our residue maps when $X$ is a trait, $Z$ its closed point:
 for motivic cohomology in bidegree $(n,n)$, our residue map
 coincides with Milnor residue symbol. More generally,
 the residue map on symbols in classical cohomologies
 always agree with Milnor residue map
 (see Par. \ref{num:res&symbols} for details).
 More interestingly, we prove here that the residue
 we have defined by purely geometrical means agrees with
 \emph{Tate residue} on De Rham cohomology (Th. \ref{prop:compare_res_Tate})
 -- this is an application of the analytical invariance of
  our residue map. This implies in particular, that for divisors,
 the residue map we have defined here -- by deformation to the normal
 bundle and orientation theory -- agrees with that of Leray.
 These classical constructions gets also extended to rigid cohomology
 (see Ex. \ref{ex:residues_Weil}).

The second interesting map is the so-called
 \emph{Gysin morphism}\footnote{This is the usual terminology
 in algebraic geometry. It is named that way
  after the pioneering work of Gysin that
	we have described in the historical part.}
 associated with the immersion $i$;
  in other words, the covariant functoriality of cohomology.
 According to our formalism, it is simply equal
 to the multiplication by the refined fundamental class
 (which gives a class with support in $Z$)
 followed by the obvious map which forgets the support
 (see \eqref{eq:Gysin&rfdl} for details).
 Thus the good properties of (refined) fundamental classes
  give all the basic expected properties of these 
  particular Gysin maps.
 
Once all this ground work is in place,
 one can introduce the main construction of orientation
 theory, that of \emph{Gysin morphisms} for certain projective 
 morphisms. In our work, it comes mainly in two settings:
  the geometric one, for schemes smooth over some fixed base,
  and the arithmetic one, for any regular schemes.
 The way we phrased and use absolute purity allows us
  to treat these two cases in a single turn, thus building
	Gysin morphisms associated with any projective morphism between
  one of these two kind of schemes.
 In fact, the process to build these Gysin maps
  is very classical: you need to consider
  the case of closed immersions and of projections of a
  projective bundle, and then find the sufficient condition
  so that these two cases can be glued.
 However, because we deal with general oriented spectra,
  the case of the projection of a projective bundle is not
  trivial. As in our previous work on the subject (\cite{Deg8}),
  to treat that case, we use a duality argument that can
 be summarized as follows: because of the projective bundle formula,
 the cohomology of a projective bundle $P/X$ is a finite free module
 over the cohomology of the base: thus it is dualizable.
 Then, one shows that the Gysin morphism associated with
 the diagonal embedding of $P$ induces a duality pairing
 on $\E^{**}(P)$. Finally one defines the Gysin morphism
 of the projection $P \rightarrow X$ as the transpose of
 the pullback with respect to this duality
(see Def. \ref{df:Gysin_Pn}).
 The condition for gluing then follows,
  using a computation of characteristic classes,
  as well as the main properties of Gysin morphisms
  for closed immersions:
  compatibility with composition, projection formula,
  compatibility with transversal pullbacks,
  excess of intersection formula.

One of the basic examples of a representable
 cohomology theory is \emph{Beilinson motivic cohomology}
 which is representable by absolutely oriented ring spectrum
 because of the ground work \cite{CD3}. Our construction
  gives Gysin morphisms for this cohomology with respect to
  any projective morphism between regular schemes. This is
  an improvement of the constructions of \cite[Th. 9]{Sou}.
For more general cohomology theories, such as algebraic K-theory
 or algebraic cobordism,
 Gysin morphisms reveal a third kind
 of characteristic classes, the \emph{cobordism classes}
 (see Def. \ref{df:cobordism}). Note these classes are interesting
 only when the formal group law associated with the considered
 oriented ring spectrum is non additive -- see in particular
 formula \eqref{eq:Myschenko}. The non triviality of these
 classes in the case of non additive formal group law 
 explains why the definition of Gysin morphisms in our context
 is far more difficult than in the ordinary case.

As a prelude to the Riemann-Roch formula,
 note that, inspired by a result of Panin,
 we give a uniqueness statement for Gysin morphisms by simple axioms
 (see Th. \ref{thm:compare_Gysin}).
 This allows us to compare the Gysin morphisms defined here
  with more classical constructions
	(see the examples  in \ref{ex:comput_Gysin}).

\bigskip

\noindent \textbf{Riemann-Roch formulas.}
The beauty of the Grothendieck-Riemann-Roch formula lies
 in its generality and the simplicity of its proof.
 The same phenomena happens here and the reader will
  see that the main technical work was to define Gysin
  morphisms.

The topological interpretation of the Riemann-Roch
 formula can be summarized as an answer to the following question:
 what happens if we \emph{change the choice of orientation}
 on an absolutely pure oriented spectrum $\E$ ? \\
The answer is that all the data defined through the
 orientation theory, as described above, 
 change and in particular the Gysin morphism:
 the change of the later is exactly 
 measured by the Riemann-Roch formula
 (see Section \ref{sec:change_orientation} for that point of view).

The setting of our general Riemann-Roch formula is to consider
 two absolute oriented ring spectra $(\E,c)$, $(\F,d)$
 and an arbitrary morphism of ring spectra
 $\varphi:\E \rightarrow \F$.
 Then one realizes that $\varphi(c)$ is an orientation of $\F$,
 which does not necessarily coincide with the given orientation $d$.
 Thus, we come back to the question of changing the orientation
 on a given spectrum and understanding its effect on Gysin morphisms.
 In a word, this change of orientation is measured by the Todd
 class $\td_\varphi$  associated with the morphism $\varphi$
 (see Prop. \ref{prop:todd}).
 With this definition in hands,
 we prove the analogue of the classical Grothendieck-Riemann-Roch formula:
 for any local complete intersection, any projective morphism
 $f:Y \rightarrow X$ in $\base$ with virtual tangent bundle $\tau_f$,
 one has for a cohomology class $y \in \E^{**}(Y)$
 (see also Th. \ref{thm:projective_RR}):
\begin{equation} \tag{$\text{GRR}_\varphi$} \label{eq:RR_intro}
\varphi_X(f_*(y))=f_*\big(\td_\varphi(\tau_f) \cupp \varphi_Y(y)\big).
\end{equation}
At this point, the proof is quite easy and follows the initial proof
 of Borel and Serre: one reduces to the case of closed regular embeddings
 and projective smooth morphisms. In our situation, the last case follows
 by duality from the first one -- here the proof differs slightly from
 that of Borel and Serre. The case of closed embedding is by reduction
 to divisors in which case it is tautological.

Inspecting this last case,
 we fall onto the surprising result that
 the Riemann-Roch formula for closed regular embedding $i:Z \rightarrow X$
 has a companion formula which was never observed till now and involves 
 the residue map $\partial_{X,Z}$. We call it the residual Riemann-Roch
 formula.
 Using the Todd class of the normal
 bundle $N_ZX$ of $Z$ in $X$, it reads as follows
 for a cohomology class $u \in \E^{**}(X-Z)$
 (see also Th. \ref{thm:RR1}):
\begin{equation} \tag{$\partial$RR} \label{eq:dRR_intro}
\partial_{X,Z}\big(\varphi_U(u)\big)
=\td_\varphi(-N_ZX) \cupp \varphi_Z\big(\partial_{X,Z}(u)\big).
\end{equation}

\bigskip

\noindent \textbf{Applications and comparison with the work
 of Gillet and Soul\'e.}
Before going to the applications, 
 let us explain with more details the crucial situation of changing
 between two orientations $c$ and $d$ on an absolute ring spectrum $\E$.
 As explained above, the orientations $c$ and $d$ correspond
  over a base scheme $X$ to formal group laws $F_c$ and $F_d$
  with coefficients in the ring $\E^{**}(X)$.
  Moreover, it follows automatically that these formal group laws
   are isomorphic: say from $F_d$ to $F_c$,
   the isomorphism corresponds to a power series $\Psi(t)$ with coefficients
   in $\E^{**}(X)$ of the form $\Psi(t)=(t+\hdots)$.
  Then the Todd class which corresponds to changing the orientation
   from $c$ to $d$ is uniquely defined in terms of $\Psi(t)$
   (see section \ref{sec:change_orientation} for a detailed discussion).

This understanding of Todd classes allows to enlighten the classical
 case. When $\base$ is the category of regular schemes,
 Quillen algebraic K-theory (resp. Beilinson motivic cohomology)
 is representable by an oriented absolute ring spectrum $\KGL$
 (resp. $\HB$).
 Then the classical (higher) Chern character as defined by Gillet
 corresponds to an isomorphism of absolute ring spectra:
$$
\ch:\KGL \rightarrow \oplus_{i \in \ZZ} \HB(i)[2i]
$$
that was first given in these terms by Riou
 (in a slightly different form, see \cite[6.2.3.9]{Riou}).
The formal group law associated with $\KGL$ (resp. $\HB$)
 is the multiplicative (resp. additive) formal group law
 -- see Example \ref{ex:abs_ring_spectra_FGL}. It follows from the
 theory of formal group laws that there is only one isomorphism
 of the form $\Psi(t)$ from the additive to the multiplicative formal
 group law: this is the exponential power series and one recovers the
 classical definition of the Todd power series, with a conceptual
 explanation why it has this precise form.
The Riemann-Roch theorem that we get for $\ch$ is a cohomological
 version of the formulas obtained by Gillet, \cite[Th. 4.1]{Gillet}:
 in fact, one will recognize most of the principles of motivic homotopy
 theory in the work of Gillet. The main difference is that we have built
 the Gysin morphisms appearing in the formula while in \emph{loc. cit.}
 they are part of the axioms (of a ``duality theory with support'').
 Therefore, the cohomological formulation of our Riemann-Roch formula
  for $\ch$ is new, valid for any projective morphism between
  regular schemes. It extends the classical Grothendieck-Riemann-Roch
  formula for $K_0$ to higher degrees.
 Note this formula was proved by Riou in \cite{Riou} in
 the particular case of smooth projective morphisms
 and was obtained in full generality in \cite{HS}
 but the construction of Gysin morphisms in \emph{op. cit.}
 is merely a reference to \cite{Deg8}
 -- see Example 2.4 in \emph{op. cit.}
 Note finally we have given a special attention to the residual
  Riemann-Roch formula in the case of the Chern character
  resulting in some purely algebraic results
  (see Section \ref{sec:application_residual}).

Our other examples are more concerned with pure orientation theory.
 First the analysis of the change of orientations shows
 that on motivic cohomology, as well as on mixed Weil cohomologies
 (they are representable according to \cite{CD2}),
 there is \emph{only one} possible orientation
  (see \ref{cor:uniqueness_orientation} for the general statement). 
 Therefore,
 on these cohomologies, there is only one possible theory of Chern classes
 and Gysin morphisms (see also 
  Th. \ref{thm:compare_Gysin} and Ex \ref{ex:comput_Gysin}).
 This settles the point of unnecessary questions about signs, uniqueness,
  that frequently occur in this kind of situation.
We also give a new proof, as well as a conceptual explanation,
 of a formula of Quillen for computing the cobordism class
 of a non necessarily trivial projective bundle,
  but only with rational coefficients (see Th. \ref{thm:Quillen_formula}
  and Example \ref{ex:Quillen_cobord_formula}).\footnote{This is enough
   to get the integral formula as well as the case of schemes
   of characteristic $0$: see Remark \ref{rem:Quillen_formula}(2).}
 The idea of the proof is to look at the Riemann-Roch formula
  that one gets by changing the formal group law of algebraic cobordism
  (conjecturally the universal formal group law)
  to the additive one -- this is why we need rational coefficients.

\bigskip

\noindent \textbf{Comparison with the work of Panin
 and \'etale cohomology.}
The last setting to which our work can be compared is the fundamental
 work of Panin on oriented cohomology theories. In \cite{PaninRR}, Panin
 proves our general Riemann-Roch formula for cohomology theories
 satisfying suitable axioms but defined over the category
 of smooth $k$-schemes. So our work should be viewed as an extension
 of the axiomatic of oriented cohomology theory to the arithmetic case.
 This is what we prove in the last section of this paper:  
  we have extracted a list of all the properties of representable
  cohomology theories that we have used in this work,  
  dubbed here \emph{arithmetic cohomology} for short
  (see Def. \ref{df:coh_thy}).\footnote{A
  longer but more precise terminology is also introduced:
  \emph{absolutely pure oriented ringed cohomology with support}.}
 This axiomatic indeed generalizes the one of Panin
  and the proofs of this paper show that our results 
  still apply to it, yielding Gysin morphisms, residues,
  characteristic classes
  and their formulas such as the Riemann-Roch ones.
	
This point of view is in fact useful because it applies especially
 to the \'etale $l$-adic cohomology of $\ZZ[1/l]$-schemes.
 In fact, we do not know if this cohomology is representable
  by an \emph{absolute} ring spectrum.\footnote{According to the results
  of \cite{CD4}, it is representable over any scheme by a ring spectrum
  but we do not know these ring spectra form a cartesian section of $\SH$
   -- unless one restricts to schemes over a field.}
But however,
 it satisfies all the axioms of an arithmetic cohomology
 -- in particular because of Thomason result about absolute purity.
Thus, the constructions of this paper apply to that cohomology,
 and its rational version, and in particular give Gysin morphisms.
 We deduce the classical Riemann-Roch formula for rational $\ell$-adic
 cohomology of regular schemes, as expected by Grothendieck
 in \cite[XIV, 6.1]{SGA6} and even the higher Riemann-Roch formula
 as well as the residual Riemann-Roch formula \eqref{eq:dRR_intro}
 (see Cor. \ref{cor:etale_RR}).
Note finally that one can also apply
 the recent work of \cite{CD4} on $\h$-motives,
 to get that \'etale motivic cohomology with coefficients in any ring
 $R$ is an arithmetic cohomology.\footnote{Recall that in the case
 where $R$ is a torsion ring of characteristic exponent invertible
 on the schemes considered, this later cohomology agrees with the
 usual \'etale cohomology with coefficients in $R$.}
 The Gysin morphisms that one get extend to arbitrary coefficients,
 in the regular case, the recent construction of Gabber-Riou
 (\cite{GRiou}, see remark \ref{rem:compare_Riou} for the comparison).

\subsection*{Outline of the work}

The paper is organized as follows. In the first section,
 we recall the basics on the motivic stable homotopy category
 give the basic properties of representable cohomologies,
 and states the absolute purity property. We end-up
 the section with a discussion about analytical invariance.
 In section \ref{sec:orientation},
 we recall orientation theory in motivic homotopy theory,
 gives the construction of (refined) fundamental classes and their
 properties.
 In section \ref{sec:Gysin}, we define the residues
 and Gysin morphisms and study their properties.
 Section \ref{sec:RR} is centered around
 the general Riemann-Roch formula and its proof.
 In section \ref{sec:examples}, we treat examples and give
 applications of our formulas.
 Finally, Section \ref{sec:Panin} discusses the axiomatic of
 arithmetic cohomologies and the case of \'etale cohomology,
 motivic \'etale cohomology and continuous \'etale cohomology
 with various coefficients.

\section*{Notations and conventions} \label{sec:notations}

All schemes in this paper are assumed to be noetherian
 of finite dimension.
 We will say that an $S$-scheme $X$, or equivalently its structure morphism,
 is \emph{projective} if it admits an $S$-embedding into $\PP^n_S$ 
 for a suitable integer $n$.\footnote{For example, if one works with
 quasi-projective schemes over a noetherian affine scheme
 (or more generally a noetherian scheme which admits an ample line bundle), 
 then a morphism is proper if and only if it is projective with our
 convention -- use \cite[Cor. 5.3.3]{EGA2}.}
 
In the whole text, unless stated otherwise,
 $\base$ stands for a sub-category of the category of such schemes.
 We will assume that $\base$ is stable by blow-up and contains
 any open subscheme of (resp. projective bundle over) a scheme
 in $\base$. The category $\base$ can be the category of all schemes,
 especially in the examples and definitions which do not deal with
 absolute purity. On the contrary, when dealing with the absolute
 purity property, the relevant examples for applications area:
\begin{itemize}
\item the category $\reg$ of all regular schemes (\emph{i.e.}
 its local rings are regular.)
\item the category $\sm_S$ of smooth $S$-schemes,
 for an arbitrary base scheme $S$.
\end{itemize}

By convention, unless explicitly stated, when we speak of the rank
 of a vector bundle, the dimension of a morphism, the codimension
 of a closed subscheme (or a closed immersion), it will always be
 assumed to be \emph{constant}.\footnote{If one wants to
  avoid this convention, see remarks \ref{rem:non_constant_rk}
  and \ref{rem:non_constant_codim}.}

Given a projective bundle $P/X$ associated with a vector bundle $E/X$,
 we will call \emph{canonical line bundle} on $P/X$ the tautological
 line bundle $\lambda$ on $P$ characterized by the property
 $\lambda \subset P \times_X E$ (it corresponds with the notation
 $\lambda=\mathcal O(-1)$).

The letter $\NN$ denotes the set of non negative integers.

\section{Absolute cohomology and purity} \label{sec:abs_coh&pur}

\subsection{Functoriality in stable homotopy}

\begin{num} \label{num:functoriality_properties_SH}
Recall that the stable homotopy category of schemes defines a 2-functor
 from the category of schemes to the category of symmetric monoidal closed triangulated categories.
This means that for any morphism of schemes $f:T \rightarrow S$, we have a pullback functor
$$
f^*:\SH(S) \rightarrow \SH(T) 
$$
which is symmetric monoidal and such that for any composable morphisms of schemes $f,g$, we have the relation:
$f^*g^*=(gf)^*$.

We will use the following properties:
\begin{enumerate}
\item[(A1)] For any morphism (resp. smooth morphism) $f$,
 the functor $f^*$ admits a right (resp. left) adjoint denoted by $f_*$ (resp. $f_\sharp$).
\item[(A2)] For any cartesian square:
$$
\xymatrix@=14pt{
Y\ar^q[r]\ar_g[d] & X\ar^f[d] \\
T\ar^p[r] & S \\ 
}
$$
such that $f$ is a smooth morphism, the base change map
$$
p^* f_\sharp \rightarrow g_\sharp q^*
$$
is an isomorphism.
\item[(A3)] For any smooth morphism $f:Y \rightarrow X$ and any spectrum $\E$ over $Y$
 (resp. $\F$ over $X$), the canonical transformation:
$$
f_\sharp(\E \wdg f^*(\F)) \rightarrow f_\sharp(\E) \wdg \F
$$   
is an isomorphism.
\item[(A4)] For any closed immersion $i:Z \rightarrow X$ with complementary open immersion
 $j$,
 there exists a unique natural transformation
 $\partial_i:i_*i^* \rightarrow j_\sharp j^*[1]$
which fits in a distinguished triangle of the form:
$$
j_\sharp j^* \xrightarrow{ad'} Id \xrightarrow{ad} i_*i^*
 \xrightarrow{\partial_i} j_\sharp j^*[1]
$$
where $ad$ (resp. $ad'$) is the unit (resp. counit) map 
 of the adjunction $(i^*,i_*)$ (resp $(j_\sharp,j^*)$).
\end{enumerate}
Except for the last property,
 these are easy consequences of the construction of $\SH$
  -- see \cite{Ayoub2}.
Property (A4) is a consequence of \cite[\S 3, th. 2.21]{MV}
 -- see \cite[4.5.47]{Ayoub2} for details.

\begin{rem} \label{rem:nil-immersion}
One can deduce from (A4) that $i_*$ is fully faithful.
Moreover, when $i$ is a nil-immersion, $i^*$ is fully faithful,
 so that $(i^*,i_*)$ is an equivalence of categories. 
\end{rem}

Note also that we can derive from properties (A1) and (A4) the following ones:
\begin{enumerate}
\item[(A5)] For any closed immersion $i$,
 the functor $i_*$ admits a right adjoint denoted by $i^!$.
\item[(A6)] For any cartesian square:
\begin{equation} \label{eq:proto_cartesian}
\xymatrix@=14pt{
T\ar^k[r]\ar_g[d] & Y\ar^f[d] \\
Z\ar^i[r] & X
}
\end{equation}
such that $i$ is a closed immersion, the base change morphism 
$$
f^*i_* \rightarrow k_*g^*
$$
is an isomorphism.
\item[(A7)] For any closed immersion $i:Z \rightarrow X$ and any spectrum $\E$ over $Z$
 (resp. $\F$ over $X$), the canonical transformation:
$$
i_*(\E \wdg i^*(\F)) \rightarrow i_*(\E) \wdg \F
$$   
is an isomorphism.
\end{enumerate}
For these three last properties, we refer the reader to \cite{Ayoub1}
 or for a more compact reference to \cite{CD3}, 2.3.3, 2.3.8 and 2.3.15
 respectively.
\end{num}

\begin{num} \label{num:i_*_monoidal}
For any smooth $S$-scheme $X$, we denote by $\Sus_S X$ the infinite suspension spectrum
 associated with the sheaf represented by $X$ with a base point added.
 Recall that $\Sus_S X=f_\sharp(\un_X)$.

Consider again the notations of axiom (A4).
We simply denote by $X/U$ the cokernel of the map induced by $j$
 in the category of Nisnevich sheaves of sets over $\sm_S$.
 Note it is pointed by identifying $U$ with the base point.
As $j$ is a monomorphism, we deduce a homotopy cofiber sequence
$$
U \rightarrow X \rightarrow X/U
$$
in the $\AA^1$-local model category of simplicial sheaves over $S$
 (\cite [Sec. 3.2, p. 105]{MV}).
According to (A4) and Remark \ref{rem:nil-immersion},
 the canonical map 
$\un_Z \rightarrow i^* \Sus(X/U)$ is an isomorphism.

To simplify notations,
 we will still denote by $X/U$ the infinite suspension spectrum associated
 with the sheaf $X/U$.
Given any object $K$ of $\SH(X)$, we get a canonical isomorphism:
$i^*\big((X/U) \wdg K \big)=i^*(K)$ which, by adjunction, induces a map
\begin{equation} \label{eq:i_*_monoidal}
(X/U) \wdg K \rightarrow i_*i^*(K).
\end{equation}
In fact, the localization axiom (A4) for $i$ is equivalent to the fact
 that this map is an isomorphism and $i_*$ is fully faithful -- see \cite[2.3.15]{CD3}.  
\end{num}

\begin{rem}
In what follows, we will consider the isomorphisms listed in the above properties
 as identities unless it involves a non trivial commutativity statement.
\end{rem}

\subsection{Absolute cohomology} \label{sec:abs_coh}

As usual, a ring spectrum over a scheme $S$ will be a commutative monoid
 of the symmetric monoidal category $\SH(S)$. 
\begin{df} \label{df:abs_ring_sp}
An \emph{$\base$-absolute spectrum} (resp. \emph{ring spectrum})
 $\E$ is
 a collection of spectra (resp. ring spectra) $\E_X$ over $X$
 for a scheme $X$ in $\base$
 and the data for any morphism $f:Y \rightarrow X$ in $\base$
 of an isomorphism of spectra (resp. ring spectra)
 $\epsilon_f:f^*\E_X \rightarrow \E_Y$ satisfying the usual cocycle
 condition.\footnote{In other words,
  this is a cartesian section of the $\base$-fibered category of spectra
  (resp. ring spectra).}

A \emph{morphism} $\varphi:\E \rightarrow \F$ of $\base$-absolute spectra
 (resp. ring spectra)
 is a collection of morphisms $\varphi_X:\E_X \rightarrow \F_X$ of spectra 
 (resp. ring spectra)
 indexed by schemes $X$ in $\base$ such that for any morphism $f:Y \rightarrow X$ in $\base$,
 the following diagram commutes:
$$
\xymatrix@R=18pt@C=28pt{
{f^*\E_X}\ar^{f^*\varphi_X}[r]\ar_{\epsilon_f^\E}[d]
& {f^*\F_X}\ar^{\epsilon_f^\F}[d] \\
{\E_Y}\ar^{\varphi_Y}[r]
 & {\F_Y}.
}
$$
\end{df}
As the category $\base$ is fixed in the entire paper
 we will abusively say \emph{absolute} for $\base$-absolute.
However, when the category $\base$ is the category
 of all $S$-schemes, we will say $S$-absolute for
 $\base$-absolute.

\begin{rem}\label{rem:weak_spectrum}
In the whole paper,
 we will be interested only in $\base$-absolute spectra.
 However, it is also convenient to consider the
 sections $(\E_X)_{X \in \base}$ over $\base$ of the 
 stable homotopy category which are not necessarily cartesian
 (\emph{i.e.} the given transition morphisms
 $\epsilon_f:f^*(\E_X) \rightarrow \E_Y$ are not isomorphisms).
We will call these objects \emph{weak $\base$-absolute spectra}.
 As in the above definition, they can also be equiped 
 with a ring structure.
\end{rem}

\begin{ex} \label{ex:abs_ring_spectra}
\begin{enumerate}
\item Let $S$ be a fixed scheme.

Then, up to the choice of the isomorphisms $\epsilon_f$
 in the above definition,
 an $S$-absolute ring spectrum $\E$ is determined
 by its value on $S$. Reciprocally, given any
 ring spectrum $\E_S$ over $S$, we get a canonical
 $S$-absolute spectrum $\E$ by putting $\E_X=f^*(\E_S)$
 for any $f:X \rightarrow S$.
In the $S$-absolute case,
 we will frequently identify $\E_S$ and $\E$
 to simplify notations.

A basic example of this kind of situation
 is given by the concept of mixed Weil theory.
 Let us recall the setting. We let $k$ be a field (non necessarily perfect)
 and $K$ be a field of characteristic $0$.
 A $K$-linear mixed Weil theory $E$ over $k$
 is a presheaf of commutative differential graded $K$-algebras over
 the category of smooth affine $k$-algebras satisfying
 axioms: homotopy invariance, Nisnevich excision, dimension, stability
 and K\"unneth formula (see \cite{CD2}).
 To such a theory one canonically associates a ring spectrum $\E$
 over $k$. 
 
We then get an absolute ring spectrum over the category of
 $k$-schemes by the preceding procedure (cf. \cite [17.2.5]{CD3}).
 The original cohomology defined by $E$ for smooth affine $k$-schemes
  gets extended to any $k$-schemes. This extension is uniquely
  characterized by the $h$-descent property (cf. \cite [17.2.6]{CD3})
  and by the commutation with limit property (cf. Lemma \ref{lm:coh&lim} below).

Examples are given by the classical Weil theories: Betti,
 De Rham, geometric \'etale and rigid cohomologies.
\item The $0$-sphere spectrum $S^0$
 -- unit for the smash product --
 is obviously an absolute ring spectrum.
\item For any scheme $S$, one can consider one of the following
 ring spectra:
\begin{itemize}
\item The Beilinson motivic cohomology spectrum $\HBx S$ (see \cite{Riou,CD3}).
\item The homotopy invariant $K$-theory ring spectrum $\KGL_S$ (see \cite{Voe1,Riou}).\footnote{Recall
 it represents \emph{homotopy invariant K-theory} according to \cite{Cis}.}
\item The cobordism ring spectrum $\MGL_S$ (see \cite[2.1]{PPR}).
\end{itemize}
Then each of these examples defines an absolute ring spectrum
 denoted respectively by $\HB$, $\KGL$, $\MGL$
 (see the respective reference given above for this assertion).
\item In \cite[\textsection 6]{Voe1}, Voevodsky defined
 the \emph{motivic Eilenberg-Mac Lane spectrum} over any smooth $k$-scheme $S$
 where $k$ is a perfect field,
 representing motivic cohomology with integral coefficients.
 This construction was generalized to any base scheme $S$ in the work
 of \cite[Ex. 3.4]{DRO}.

In \cite[11.2.17]{CD3}, D.C.~Cisinski and the author gave
 a general construction of this spectrum relying on other ideas
 of Voevodsky.
 Let $\Lambda$ be a localization of the ring $\ZZ$
 and $S$ be any scheme.
 Then one defines a ring spectrum $\HMx \Lambda S$
 such that the abelian group
$$
\Hom_{\SH(S)}(S^0,\HMx \Lambda S(m)[n])
 =H_{\mathcal M}^{n,m}(S,\Lambda)
$$
is Voevodsky's motivic cohomology of $S$ with coefficients
 in $\Lambda$.\footnote{Note that, in the current state of the theory,
 the consideration of the ring $\Lambda$ for coefficients
 is crucial as, given a localization $\Lambda'$ of $\Lambda$, 
 we do not always have:
$$
\HMx {(\Lambda')} S
 \simeq \HMx \Lambda S \otimes_\Lambda \Lambda'.
$$
See however Proposition 11.2.19 of \emph{loc. cit.} for cases where
 this identification holds.
}
According to \cite[11.2.21]{CD3},
 for any morphism $f:S' \rightarrow S$ of schemes, there
 exists a  canonical morphism of ring spectra:
$$
\tau_f:f^*(\HMx \Lambda S)
 \rightarrow \HMx \Lambda {S'}.
$$
In other word, we have a weak absolute ring spectrum
 (over the category of all schemes) that we will denote
 $\HM\Lambda^w$.
 A fundamental conjecture of Voevodsky (cf \cite[Conj. 17]{Voe2}
 for the case of integral coefficients)
 can be reformulated saying that $\HM\Lambda^w$ is in fact
 a (strong) $\base$-absolute ring spectrum over the category
 of all schemes.

Unfortunately, this conjecture is still not known.
 In \cite[16.1.7]{CD3}, we prove it when $\Lambda=\QQ$
 if one restricts to the case of geometrically unibranch
 schemes.\footnote{In fact, we proved that the ring spectrum
 $\HMx \QQ S$ coincides with $\HBx S$ when $S$ is geometrically
 unibranch.}

On the other hand, by construction, the above map $\tau_f$
 is obviously an isomorphism when $f$ is smooth.\footnote{Using
 the notations of \cite[Par. 11.2.21]{CD3}, when $f$ is smooth,
 the exchange morphism $f^*\varphi_* \rightarrow \varphi_*f^*$
 used to construct $\tau_f$ is an isomorphism.}
Thus, if $\base$ is the category of smooth schemes
 over an arbitrary base scheme $S$,
 then the collection $\HMx \Lambda X$
 defines an $S$-absolute ring spectrum simply denoted by
 $\HM\Lambda_{/S}$.
 Moreover, from the perspective of this paper,
 we can safely extend $\HM\Lambda_{/S}$ as an $S$-absolute
 ring spectrum taking the various pullbacks of $\HMx \Lambda S$
 over any $S$-scheme as in the first Example above. 
This is particularly relevant if $S$ is the spectrum of a prime field
 (we refer the reader to \cite{CD4} for more details).
\item Recall finally a very interesting construction of Markus Spitzweck.
 Let $\Lambda$ be a localization of $\ZZ$.
 In \cite{Spit}, Spitzweck defines a ring spectrum $\HMSx \Lambda \ZZ$
 over $\ZZ$
 (\emph{loc. cit.}, Def. 4.27)
 whose pullbacks to any field $\spec k \rightarrow \spec \ZZ$
 is Voevodsky's Eilenberg-Mac Lane spectrum
 (\emph{loc. cit.} 6.7, 9.16 and 9.17). Thus, the absolute ring spectrum
 $\HMS \Lambda$
 obtained by considering for any scheme $f:X \rightarrow \spec \ZZ$
 the ring spectrum $f^*(\HMSx \Lambda \ZZ)$ is an interesting candidate
 for motivic cohomology in general. Note in particular that 
 according to \emph{loc. cit.}, 7.19, for any smooth $\ZZ$-scheme $X$,
$$
\HMS \ZZ^{n,m}(X)=CH_{2m-n}(X,m)
$$
where the right hand side is Bloch's higher Chow group as defined by Levine.
\end{enumerate}
\end{ex}

\begin{num}
In the remainder of this section,
 we consider an absolute ring spectrum $\E$
 (which can be weak in the sense or Remark \ref{rem:weak_spectrum}
  until Paragraph \ref{num:products}) and define structures
 on its associated cohomology theory.

As usual, we call closed (resp. open) pair any couple $(X,Z)$ (resp. $(X/U)$)
 such that $X$ is a scheme and $Z$ (resp. $U$) is a closed (resp. open) subscheme
 of $X$.
\end{num}
\begin{df} \label{df:coh_support}
Given a closed pair $(X,Z)$, corresponding to a closed immersion $i$,
 and a couple $(n,m) \in \ZZ^2$,
 we define the relative cohomology of $(X,Z)$ in bi-degree $(n,m)$ with coefficients
 in $\E$ by one of the following equivalent formulas: 
\begin{align*}
\E^{n,m}_Z(X):=&\Hom_X(X/X-Z,\E_X(m)[n]) \\
&=\Hom_X(i_*(\un_Z),\E_X(m)[n])
=\Hom_Z(\un_Z,i^!\E_X(m)[n]).
\end{align*}
When it clarifies the notation, we will also put:
 $\E^{n,m}(X,Z):=\E^{n,m}_Z(X)$.

When $X=Z$ we simply put as usual: $\E^{n,m}(X):=\E^{n,m}_X(X)$.
\end{df}
A morphism $\varphi:\E \rightarrow \F$
 of absolute ring spectra obviously induces 
 for any closed pair $(X,Z)$ a morphism
\begin{equation} \label{eq:nat_transfo_coh}
\varphi_*:\E_Z^{n,m}(X) \rightarrow \F_Z^{n,m}(X). 
\end{equation}

\begin{num} \label{num:pullback_support}
\textit{Contravariant functoriality}: Consider a closed pair $(X,Z)$
 and a morphism $f:Y \rightarrow X$ in $\base$.
Let $T=Y \times_X Z$ be the pullback in the category of schemes,
 considered as a closed subscheme of $Y$. 
Then we get a morphism of abelian groups:
$$
f^*:\E^{n,m}_Z(X) \rightarrow \E^{n,m}_T(Y)
$$
by one of the following equivalent definitions:
\begin{itemize}
\item According to Remark \ref{num:i_*_monoidal},
  $f^*(X/X-Z)=Y/Y-T$.
Thus, for a cohomology class $\rho:(X/X-Z) \rightarrow \E_X(m)[n]$,
 the pullback map $f^*(\rho)$ gives the desired map:
$$(Y/Y-T)=f^*(X/X-Z) \rightarrow f^*\E_X(m)[n]
 \xrightarrow{\epsilon_f} \E_Y(m)[n].$$
\item Consider the pullback square \eqref{eq:proto_cartesian}
 with $f$ as above
 and $i$ the immersion of $Z$ in $X$. Then property (A6) gives a canonical identification:
$f^*i_*(\un_Z)=k_*g^*(\un_Z)=k_*(\un_T).$
Thus, taking a cohomology class $\rho:i_*(\un_Z) \rightarrow \E_X(m)[n]$,
 the pullback map $f^*(\rho)$ gives the desired map:
$$k_*(\un_T)=f^*i_*(\un_Z) \rightarrow f^*\E_X(m)[n]
 \xrightarrow{\epsilon_f} \E_Y(m)[n].$$
\end{itemize}
\end{num}

\begin{num}\label{num:pushforward_support}
\textit{Covariant functoriality}: Consider closed immersions
 $T \xrightarrow \nu Z \xrightarrow i X$ and put $k=i \circ \nu$. \\  
 We define a pushforward in cohomology with support
$$
\nu_!:\E^{n,m}_T(X) \rightarrow \E^{n,m}_Z(X)
$$
by one of the following equivalent definitions:
\begin{itemize}
\item By functoriality of homotopy colimits,
 the immersion $(X-Z) \rightarrow (X-T)$ induces
 a canonical map $\bar \nu:(X/X-Z) \rightarrow (X/X-T)$ in $\SH(X)$.
 Then we associate to a cohomology class $\rho:(X/X-T) \rightarrow \E_X(m)[n]$
 the map $\nu_!(\rho):=\rho \circ \bar \nu$.
\item 
The unit map of the adjunction $(\nu^*,\nu_*)$ gives a morphism
 $$ad_\nu:i_*(\un_Z) \rightarrow i_*\nu_*\nu^*(\un_Z)=k_*(\un_T)$$ in $\SH(X)$.
 Then for any cohomology class $\rho:k_*(\un_T) \rightarrow \E_X(m)[n]$,
 we put: $\nu_!(\rho):=\rho \circ ad_\nu$.
\end{itemize}
\end{num}

\begin{num} \label{num:products}
\textit{Products}: Consider the assumption of the preceding paragraph,
 except that now we ask $\E$ is an absolute ring spectrum.
 We define a \emph{refined product}
 in cohomology with supports
\begin{equation} \label{eq:product1}
\E_T^{n,m}(Z) \otimes \E_Z^{s,t}(X) \rightarrow \E_T^{n+s,m+t}(X),
 (\lambda,\rho) \mapsto \lambda \refp \rho
\end{equation}
as follows.
Consider cohomology classes:
$$\lambda:\nu_*(\un_T) \rightarrow \E_Z(t)[s],
 \ \rho:i_*(\un_Z) \rightarrow \E_X(m)[n].$$
Applying $i_*$ to $\lambda$, we get a map:
$$
k_*(\un_T)=i_*\nu_*(\un_T)
 \rightarrow i_*(\E_Z(m)[n])
 \xrightarrow{\epsilon_i^{-1}} i_*i^*(\E_X(m)[n])
 \simeq \E_X \wdg i_*(\un_Z)(m)[n] 
$$
where the last identification uses (A7).
Let us simply denote by $i_*(\lambda)$ this composite map.
We define the product $\lambda \refp \rho$ as the following 
 composite morphism:
$$
k_*(\un_T) \xrightarrow{i_*(\lambda)} \E_X \wdg i_*(\un_Z)(m)[n]
 \xrightarrow{Id \wdg \rho} \E_X \wdg \E_X(m+s)[n+t]
 \xrightarrow{\mu_X} \E_X(m+s)[n+t]. 
$$
One can deduce from this product the usual product of cohomology
 with support as follows.
Assume we are given a cartesian square of closed immersions:
$$
\xymatrix@=14pt{
T\ar^{i'}[r]\ar_\nu[d] & Z'\ar^\iota[d] \\
Z\ar^i[r] & X
}
$$
Then one defines the \emph{cup-products} by the following formula:
\begin{equation} \label{eq:product2}
\E_{Z}^{n,m}(X) \otimes \E_{Z'}^{s,t}(X) \rightarrow \E_T^{n+s,m+t}(X),
 (\alpha,\beta) \mapsto \alpha \cupp \beta=\iota^*(\alpha) \refp \beta.
\end{equation}
One can also describe this product using the following identification in $\SH(X)$:
$$
(X/X-Z) \wdg (X/X-Z')=(X/X-T).
$$
This can be obtained by a direct computation of homotopy colimits
 or by applying formula \eqref{eq:i_*_monoidal} and (A6) as follows:
$$
(X/X-Z') \wdg (X/X-Z)=\iota_*(\un_{Z'}) \wdg i_*(\un_Z)
=\iota_*(\iota^*i_*(\un_Z))=\iota_*i'_*\nu^*(\un_Z)=(X/X-T).
$$ 
Then given $\alpha:(X/X-Z) \rightarrow \E(m)[n]$
 and $\beta:(X/X-Z') \rightarrow \E_X(t)[s]$, one checks that $\iota^*(\alpha) \refp \beta$
 is equal to the following composite map:
\begin{align*}
(X/X-T)&=(X/X-Z) \wdg (X/X-Z')
 \xrightarrow{\alpha \wdg \beta} \E_X \wdg \E_X(m+s)[n+s] \\
& \xrightarrow{\mu_X} \E_X(m+s)[n+r]. 
\end{align*}
Note that when $Z=Z'=X$
 the product \eqref{eq:product2} describes the usual cup-product in $\E$-cohomology.
\end{num}

\begin{rem}\label{rem:pushforward_support}
The need of the refined product just defined is
 the only reason for us to use absolute ring spectra.
 Indeed, we have not only used the
 existence of the structural map $\epsilon_f$ but
 also the fact it is an isomorphism.
\end{rem}

We have gathered the basic properties of these operations in the following proposition:
\begin{prop}\label{prop:ppty_product}
 Given an absolute ring spectrum $\E$,
  the following properties hold:
\begin{enumerate}
\item[(E1)] $f^*g^*=(gf)^*$, $\nu'_!\nu_!=(\nu'\nu)_!$ whenever defined.
\item[(E2)] When $\nu$ is a closed nil-immersion,
 $\nu_!$ is an isomorphism.
\item[(E3)] Consider the following cartesian squares:
$$
\xymatrix@=14pt{
T'\ar^{\nu'}[r]\ar[d] & Z'\ar[r]\ar^g[d] & X'\ar^f[d] \\
T\ar^\nu[r] & Z\ar[r] & X
}
$$
where horizontal maps are closed immersions.
Then for any cohomology class $\rho \in \E^{**}_T(X)$, $f^*\nu_!(\rho)=\nu'_!f^*(\rho)$.
\item[(E4)] Consider closed immersions:
 $W \xrightarrow \iota T \xrightarrow \nu Z \xrightarrow i X$. \\
Then for any triple $(\lambda,\alpha,\beta) \in \E^{**}_W(T) \times \E^{**}_T(Z) \times \E^{**}_Z(X)$,
one has: \\
  $\lambda \refp (\alpha \refp \beta)=(\lambda \refp \alpha) \refp \beta$. 
\item[(E5)] Under the assumption of (E3),
 for any couple $(\lambda,\rho) \in \E^{**}_T(Z) \times \E^{**}_Z(X)$,
 one has: \\ $f^*(\lambda \refp \rho)=g^*(\lambda) \refp f^*(\rho)$.
\item[(E6)] Under the assumption of (E4), 
 for any couple $(\lambda,\rho) \in \E^{**}_T(Z) \times \E^{**}_Z(X)$,
 one has: \\ $\nu_!(\lambda \refp \rho)=\nu_!(\lambda) \refp \rho$.
\item[(E7)] Consider the following diagram:
$$
\xymatrix@=14pt{
T'\ar^{\nu'}[r]\ar_h[d] & Z'\ar^g[d] &  \\
T\ar^\nu[r] & Z\ar^i[r] & X
}
$$
made of closed immersions and such that the square is cartesian.
Then for any couple $(\lambda,\rho) \in \E^{**}_T(Z) \times \E^{**}_{Z'}(X)$, one has:
 $h_!(g^*(\lambda) \refp \rho)=\lambda \refp g_!(\rho)$.
\end{enumerate}
\end{prop}
\begin{proof}
(E1) is clear and (E2) follows from Remark \ref{rem:nil-immersion}.

Consider (E3). Recall that $\nu_!$ (resp. $\nu'_!$) is induced
 by (pre)composition with the canonical map
 $\epsilon_X:(X/X-Z) \rightarrow (X/X-T)$ 
(resp. $\epsilon_Y:(Y/Y-Z') \rightarrow (Y/Y-T')$). Then (E3) simply follows
from the fact $f^*(\epsilon_X)=\epsilon_Y$.

Point (E4) follows easily using the associativity of the product
 on the ring spectrum $\E_X$.
Point (E5) follows easily using (A6) and point (E6) is clear.

Point (E7) is the most difficult one.
 Let us denote by $i'$ (resp. $k$, $k'$) the obvious embedding of
 $Z'$ (resp. $T$, $T'$) in $X$.
The cohomology classes of the statement to be proved
 can be written as follows: $\lambda:\nu_*(\un_T) \rightarrow \E_Z$,
 $\rho:i'_*(\un_{Z'}) \rightarrow \E_X$. \\
 We let $ad_g:i_*(\un_Z) \rightarrow i'_*(\un_{Z'})$
 (resp. $ad_h:k_*(\un_T) \rightarrow k'_*(\un_{T'}$)
 be the map induced by the unit of the adjunction $(g_*,g^*)$ (resp. $(h^*,h_*)$).
Then the left hand side of the relation to be proved is equal to:
\begin{align*}
k_*(\un_T) & \xrightarrow{ad_h} k'_*(\un_{T'})=i'_*g^*\nu_*(\un_Z)
 \xrightarrow{i'_*g^*(\lambda)} i'_*g^*(\E_Z)=i'_*(\E_{Z'})
 =\E_X \wdg i'_*(\un_{Z'}) \\
 & \xrightarrow{Id \wdg \rho} \E_X \wdg \E_X \xrightarrow \mu \E_X,
\end{align*}
while the right hand side is:
\begin{align*}
k_*(\un_T)=i_*\nu_*(\un_Z)
 & \xrightarrow{i_*(\lambda)} i_*(\E_Z)=\E_X \otimes i_*(\un_Z)
 \xrightarrow{Id \wdg ad_g} \E_X \wdg i'_*(\un_{Z'})
 \\
 & \xrightarrow{Id \wdg \rho} \E_X \wdg \E_X \xrightarrow \mu \E_X. 
\end{align*}
To check the identity,
 we prove the commutativity of the following diagram:
$$
\xymatrix@C=12pt@R=20pt{
k_*h_*h^*(\un_T)\ar@{=}[r]\ar@{}|{(1)}[rrd]
 & i_*g_*\nu'_*h^*(\un_T)\ar^{b^{-1}}_\sim[r]
 & i_*g_*g^*\nu_*(\un_T)\ar^{\lambda}[r]
 & i_*g_*g^*(\E_Z)\ar@{=}[r]
 & \E_X \wdg i_*g_*g^*(\un_Z) \\
k_*(\un_T)\ar@{=}[rr]\ar^{ad_h}[u]
 && i_*\nu_*(\un_T)\ar^{\lambda}[r]\ar^{ad_g}[u]
 & i_*(\E_Z)\ar@{=}[r]\ar_{ad_g}[u]
 & \E_X \wdg i_*(\un_Z)\ar_{Id \wdg ad_g}[u]
}
$$
where the maps with label $ad_g$ (resp. $ad_h$) are induced
 by the unit of the adjunction $(g^*,g_*)$ (resp. $(h^*,h_*)$)
 and the map $b$, whose inverse appears in the above diagram,
 stands for the base change isomorphism obtained from (A6).
Only the commutativity of part (1) is non trivial.
It follows from the description
of the base change map $b$ as the composite:
$$
\xymatrix@=24pt{
g^*\nu_*\ar^-{ad_h}[r]
 & g^*\nu_*h_*h^*\ar@{=}[r]
 & g^*g_*\nu'_*h^*\ar^-{ad'_g}[r] 
 & \nu'_*h^*
}
$$
and the relation between the unit and the counit
 $ad'_g$ of the adjunction $(g^*,g_*)$.
\end{proof}

\begin{num} \label{num:closed_pairs}
It is usually convenient to introduce the following notions.
We define a \emph{morphism of closed pairs} $\Delta:(Y,T) \rightarrow (X,Z)$
 as being a commutative diagram
$$
\xymatrix@=10pt{
T\ar[r]\ar_g[d] & Y\ar^f[d] \\
Z\ar[r] & X
}
$$
such that the induced map $T \xrightarrow \nu Y \times_X Z$ is
 a nil-immersion.\footnote{In other words,
 the diagram $\Delta$ is topologically cartesian.}
 We say that the morphism $\Delta$ is \emph{cartesian}
 when the above square is cartesian in the category of schemes
 -- \emph{i.e.} $\nu$ is an isomorphism. We also use the notation
 $(f,g)$ for $\Delta$ when we want to refer to the morphisms in the
 above square. \\
The composition of morphisms of closed pairs is given in categorical terms
 by the vertical composition of squares. 

Using (E2),
 we associate to $\Delta$ the following composite morphism of abelian groups:
$$
\Delta^*:\E^{n,m}_Z(X) \xrightarrow{f^*} \E^{n,m}_{Z \times_X Y}(Y)
 \xrightarrow{\nu_!^{-1}} \E^{n,m}_T(Y).
$$
The relation $\Delta^*\Theta^*=(\Theta\Delta)^*$
 is clear from (E1) and (E3).
\end{num}

\begin{rem}
Given a transformation $\varphi:\E \rightarrow \F$,
 it is clear that the associated morphism \eqref{eq:nat_transfo_coh}
 is natural with respect to contravariant and covariant functorialities.
 Moreover,
  it is compatible with all the products of Paragraph \ref{num:products}.   
\end{rem}

The following property is not essential to our purpose
 so that we do not list it among the fundamental axioms
 of an absolute cohomology.\footnote{It will only be applied
 in Example \ref{ex:abs_pure} in the case
 of $k$-absolute ring spectra.}
\begin{lm}\label{lm:coh&lim}
Let $(X_\alpha,Z_\alpha)_{\alpha \in A}$ be an essentially affine
 projective system of closed pairs of $\base$ whose projective
 limit\footnote{in the category of pairs of schemes,}
 $(X,Z)$ is still a closed pair in $\base$.

Then, given any $\base$-absolute ring spectrum $\E$,
 the canonical map:
$$
\ilim{\alpha \in A^{op}} \big(\E^{n,i}_{Z_\alpha}(X_\alpha)\big)
 \rightarrow \E^{n,i}_{Z}(X)
$$
is an isomorphism.
\end{lm}
\begin{proof}
According to \cite[4.3.6]{CD3},
 the category $\SH$ is a continuous motivic category in the
 sense of \emph{loc. cit.}, 4.3.2.
 Thus the lemma follows by applying \emph{lo. cit.}, 4.3.4
 to the projective systems
$$
(X_\alpha/X_\alpha-Z_\alpha)_{\alpha \in A}, 
 (\E_{X_\alpha})_{\alpha \in A}
$$
given that for any index $\alpha \in A$, the object
 $(X_\alpha/X_\alpha-Z_\alpha)$ of $\SH(X_\alpha)$ is
 constructible.
\end{proof}

\subsection{Absolute purity} \label{sec:abs_purity}

\begin{num} \label{num:def_space}
Let $(X,Z)$ be a closed pair in $\base$.
 We say that $(X,Z)$ is \emph{regular} if the inclusion
 $i:Z \rightarrow X$ is a regular embedding.

Assume $(X,Z)$ is regular.
We let $N_ZX$ (resp. $B_ZX$) be the normal cone (resp. blow-up)
 of $Z$ in $X$. Recall the definition of the \emph{deformation space} of $(X,Z)$
 as: $$D_ZX=B_{0 \times Z}(\AA^1_X)-B_ZX$$
 (see \cite{Ros}, or \cite[\S 4.1]{Deg8} for this presentation).\footnote{Note
 that according to our convention on $\base$, this is a scheme in $\base$.}
 This is a flat scheme over $\AA^1$ whose fiber over $1$
 (resp. $0$) is $X$ (resp. $N_ZX$).
 Note also that $D_ZZ=\AA^1_Z$ is a closed subscheme of $D_ZX$
 so that we finally get a \emph{deformation diagram}
 of closed pairs:
\begin{equation} \label{eq:deformation}
(X,Z) \xrightarrow{\sigma_1} (D_ZX,\AA^1_Z)
 \xleftarrow{\sigma_0} (N_ZX,Z),
\end{equation}
made of cartesian morphisms.
Note this diagram is natural with respect to cartesian
 morphisms of closed pairs.
\end{num}
\begin{df} \label{df:abs_pur}
Let $\E$ be a weak $\base$-absolute (ring) spectrum
 (Remark \ref{rem:weak_spectrum}).

For any closed pair $(X,Z)$,
 we say that $(X,Z)$ is $\E$-pure if $(X,Z)$ is regular
 and the morphisms 
$$
\E^{**}(X,Z) \xleftarrow{\sigma_1^*} \E^{**}(D_ZX,\AA^1_Z)
 \xrightarrow{\sigma_0^*} \E^{**}(N_ZX,Z)
$$
induced by the above deformation diagram
 (see notation in Definition \ref{df:coh_support})
 are isomorphisms of bigraded abelian groups.

We say that $\E$ is \emph{absolutely pure}
 if any regular closed pair in $\base$ is $\E$-pure.
\end{df}

Note the following trivial stability properties of spectra
 satisfying the absolute purity property with respect
 to a given closed pair:
\begin{prop}
Let $(X,Z)$ be a regular closed pair in $\base$.

Then the category of absolute spectra (resp. weak
 absolute spectra) $\E$ such that $(X,Z)$ is
 $\E$-pure is stable by suspension and twists,
 direct factors, infinite direct sums
 and distinguished triangles.\footnote{A distinguished triangle
 of (weak) absolute spectra is the datum of distinguished triangles
 over each schemes which are compatible
 with the structural base change maps.}
\end{prop}

\begin{ex} \label{ex:abs_pure}
We refer the reader to Example \ref{ex:abs_ring_spectra}
 for the absolute ring spectra appearing below:
\begin{enumerate}
\item Let $\E$ be any absolute ring spectrum
 and $S$ be a scheme.
 We will say that a closed $S$-pair $(X,Z)$
 is smooth if $X$ and $Z$ are smooth over $S$.

According to \cite[\S 3.2, th. 2.23, p. 115]{MV},
 any smooth closed $S$-pair $(X,Z)$ is $\E$-pure.
Thus,
 if $\base$ is the category of smooth $S$-schemes,
 any $\base$-absolute ring spectrum is absolutely pure.
 Note in particular this is the case for the $\sm_S$-absolute ring
 spectrum $\HH_\Lambda$, representing motivic cohomology
 with $\Lambda$-coefficients.
\item Let $k$ be a perfect field.

According to Popescu theorem, any regular closed pair $(X,Z)$
 over $k$ can be written as a projective limit of smooth closed pairs
 over $k$ provided $X$ is regular.
 Thus, according to the computation of the cohomology with
 supports of a projective limit (cf. Lemma \ref{lm:coh&lim}),
 we deduce from the previous example that any $k$-absolute ring
 spectrum is absolutely pure.
This is in particular the case for the $k$-absolute
 motivic Eilenberg-MacLane spectrum $\HH_\Lambda$.
\item According to \cite{CD3}, respectively Theorems 14.4.1 and 13.6.3,
  the absolute ring spectra $\HB$ and $\KGL$ are absolutely pure.
\item Recall from \cite[10.5]{NSO} that there is an isomorphism
 of absolute ring spectra:
$$
\MGL \otimes_\ZZ \QQ=\HB[b_1,b_2,\hdots]
$$
where $b_i$ is a generator of degree $(2i,i)$.
Then the preceding example implies $\MGL \otimes \QQ$ is absolutely pure.
We deduce from that example that any Landweber spectrum
(cf. \cite[th. 7.3]{NSO})
with rational coefficients is absolutely pure.
\end{enumerate}
\end{ex}

The only integral example of an absolutely pure ring spectrum
 is given by homotopy invariant K-theory. However, in view
 of the previous examples, we think it is reasonable to conjecture
 the following:
\begin{conj} \label{conj:purityMGL}
The absolute ring spectrum $\MGL$ is absolutely pure.
\end{conj}

\begin{conj} \label{conj:purityS0}
The absolute ring spectrum $S^0$ is absolutely pure.
\end{conj}

\begin{rem} \label{rem:conseq_abs_pur}
\begin{enumerate}
\item Recall the notion of cellular spectra, first introduced in \cite{DI}.
 Over a scheme $S$, one defines the category of cellular spectra
 as the smallest thick triangulated
 subcategory of $\SH(S)$ which contains the spheres that we denote here
 $S^0(m)[n]$
 for any integers $(n,m) \in \ZZ$.
 According to the previous proposition,
  Conjecture \ref{conj:purityS0} implies that
  any absolute spectrum $\E$ which is cellular over $\spec{\ZZ}$
  is absolutely pure.

In particular, Conjecture \ref{conj:purityS0} implies
 Conjecture \ref{conj:purityMGL}
 because $\MGL_S$ is cellular (see \cite[6.4]{DI})
 -- and the former conjecture would reprove the absolute purity
 for $\KGL$ (see \cite[6.4]{DI}).
\item Let $S$ be any scheme 
 and $\Lambda$ be the localization of $\ZZ$ at the primes which
 are not invertible on $S$.
 Let $\HMS \Lambda$ be the absolute ring spectrum
  of Spitzweck (see Ex. \ref{ex:abs_ring_spectra}(5).
 According to \cite[Cor. 11.4]{Spit}, $\HMSx \Lambda S$ is cellular.
 Thus, according point (1) above, 
 Conjecture \ref{conj:purityS0} implies that
 $\HMS \Lambda$ is absolutely pure over the category of $S$-schemes.
\item Supporting these conjectures: 
\begin{itemize}
\item $\MGL$ and $S^0$
 are absolutely pure for closed pairs which are smooth over some base.
\item The conjecture for $\MGL \otimes \QQ$ is true according to point (3)
 of the preceding example.
\item The conjecture for $S^0 \otimes \QQ$ is true if one restricts
 to base schemes over which $-1$ is a sum of squares.
 In fact, under this assumption, according to Morel theorem
 one has
 $S^0 \otimes \QQ \simeq \HB$ (see \cite[Cor. 16.2.14]{CD3})
 and we are reduced to point (3) of the above example.
\end{itemize}
\item Let us consider the Eilenberg-MacLane motivic
 ring spectrum, as a weak absolute ring spectrum,
 $\HM\ZZ^w$ -- Example \ref{ex:abs_ring_spectra}(4).
 Independently of the conjecture of Voevodsky
 which ask whether $\HM\ZZ^w$ is an absolute ring spectrum,
 it is interesting to ask if $\HM\ZZ^w$ is absolutely pure.
 Note that if this was true,
 then using the coniveau spectral sequence
 we will get for any regular scheme $S$ an isomorphism
$$
H^{2n,n}_\mathcal M(S,\ZZ) \simeq CH^n(S)
$$
where the right hand side is the group of $n$-codimensional cycles
 in $S$ modulo rational equivalence (see \cite{Ful, Gil}).\footnote{The
 absolute purity
 property allows to compute the $E_1$-term as the $\HM\ZZ^w$-cohomology
 of the residue fields  for $\HM\ZZ^w$, as $S$ is regular. As the motivic
 cohomology of fields is known, one gets the mentioned computation after
 identifying one of the differential with the classical divisor class
 map.}
 The existence of this isomorphism is particularly interesting
 as there is a well defined product of Voevodsky's motivic
 cohomology while it is still an open question to
 define a product on the classical Chow
  (see \cite[XIV, \textsection 8]{SGA6}).
\item We have separated the case of $\MGL$ with that of $S^0$
 because $\MGL$ is oriented (see next section) and $S^0$ is not.
Note however that in our formulation of absolute purity, we do not
 need any orientation. In the particular,
 one can see that the Conjecture \ref{conj:purityS0} is equivalent
 to ask the following:
\begin{enumerate}
\item[] For any closed immersion $i:Z \rightarrow X$ between regular schemes,
 with normal bundle $N_ZX$, there exists a canonical isomorphism in $\SH(Z)$:
$$
i^!(S^0) \simeq \Th(-N_ZX)
$$
-- see Par. \ref{num:Thom_space} for recall on the definition of 
 the right hand side.
\end{enumerate}
The resulting map in $\SH(S)$
$$
i_*(\Th(-N_ZX)) \rightarrow S^0
$$
would be called the (unoriented) fundamental class of $i$
 (relative to $S^0$).
This is the unoriented version of Definition \ref{df:fdl_class}.
 Moreover,
 one can see that this class would be universal
 among the fundamental classes constructed in this paper
  (using that any ring spectrum is an algebra over $S^0$).
\end{enumerate}
\end{rem}

\subsection{Analytical invariance} \label{sec:anaylitical_inv}

By construction of the stable homotopy category,
 an absolute cohomology theory satisfies cohomological descent
 for the Nisnevich topology. A convenient way to express this
 property uses the so called \emph{excision property}.
 Let us start by an elementary geometric fact which will link excision
  with analytical invariance.
\begin{prop}
Let $f:Y \rightarrow X$ be a morphism locally of finite type,
 $Z$ a closed subscheme of $X$, $T=f^{-1}(Z)$. Then the
 following assertions are equivalent:
\begin{enumerate}
\item[(i)] $f$ is \'etale at all points of the scheme $T$
 and the induced morphism $f|_T:T_{red} \rightarrow Z_{red}$ is an isomorphism.
\item[(ii)] the induced morphism
 $\hat f:\hat Y_T \rightarrow \hat X_Z$ between the respective formal completions
 of $Y$ at $T$ and $X$ at $Z$ is an isomorphism.
\end{enumerate}
\end{prop}
\begin{proof}
Note first that given any point $x \in Z$,
 the completion of the local ring of $\hat X_Z$ at $x$
 coincides with the completion of the local ring of
 $X$ at $x$ and the corresponding isomorphism is natural
 in $(X,Z,x)$.

Thus the equivalence of the assertions follow from
 \cite[17.6.3]{EGA4} which asserts that $f$ is \'etale
 at a point $y$ of $T$, $x=f(y)$, if and only if the induced
 morphism $\hat{\mathcal O}_{X,x} \rightarrow \hat{\mathcal O}_{Y,y}$
 between the respective completed local rings is an isomorphism.
\end{proof}

\begin{df} \label{df:excisive_morph}
Let $\Delta=(f,g):(Y,T) \rightarrow (X,Z)$ be a morphism
 of closed pairs. One says that $\Delta$ is \emph{excisive} if
 the morphism $f$ and the closed scheme $Z \subset X$
 satisfy the equivalent assertions of the preceding proposition.
\end{df}
Thus, generalizing slightly \cite[Lem. 1.6, p. 98]{MV},
 we get the following property of the absolute cohomology
 represented by $\E$:
\begin{enumerate}
\item[(Nis)] For any excisive morphism
 $\Delta:(Y,T) \rightarrow (X,Z)$,
 the associated pullback $\Delta^*:\E^{**}_Z(X) \rightarrow \E^{**}_T(Y)$
 is an isomorphism.
\end{enumerate}

\begin{rem}
In fact,
 it is well known since the work of Morel-Voevodsky
 that this property characterizes Nisnevich descent,
 though one only needs excisive morphisms which are globally
 \'etale (see for example \cite[Prop. 1.1.10]{CD2} for a precise
 statement).
\end{rem}

An interesting corollary of the absolute purity property 
 is the following stronger statement
 (see the end of this section for a stronger result):
\begin{prop}\label{prop:inv_completion}
Let $\E$ be an absolutely pure $\reg$-spectrum.
Let $X$ be a regular local scheme with closed point $x$,
 $\hat X$ be its completion at the point $x$
 and consider the canonical map $f:(X,x) \rightarrow (\hat X,x)$.

Then the induced morphism of cohomology with support:
$$
f^*:\E^{**}_x(\hat X) \rightarrow \E^{**}_x(X)
$$
is an isomorphism.
\end{prop}

\begin{rem}
This result will be used to compute residues in 
 Proposition \ref{prop:compare_res_Tate}.
\end{rem}

Using the Artin approximation property,
 one can give a much stronger result than 
 the preceding proposition according to
 an initial idea of Wildeshaus (cf. \cite[\textsection 5]{Wil},
 in the case of motives over a field).
\begin{thm}\label{thm:analytic}
Let $S$ be an excellent scheme,
 $(X,Z)$ and $(Y,T)$ be closed pairs
 made of $S$-schemes essentially of finite type.

Let $\E$ be an $S$-absolute spectrum.
 
Assume there exists an isomorphism
 $\mathfrak f:\hat Y_T \rightarrow \hat X_Z$
 between the respective formal completions.
Then there exists an isomorphism:
$$
\mathfrak f^*:\E^{**}_Z(X) \rightarrow \E^{**}_T(Y)
$$
which depends only on $\mathfrak f$.
\end{thm}
\begin{proof}
The following proof follows that of \cite[5.5]{Wil}
 using the more advanced theory we now have at our
 will.
 
We can assume $Z=T$, seen as a reduced scheme.
Let us fix a point $z \in Z$,
 and let $X_{(z)}$, $Y_{(z)}$ be the respective local
  schemes of $X$, $Y$ at $z$.
  According to \cite[2.4]{Swan}, the henselisation
  of $\mathcal O_{X,z}$ satisfies
  the Artin approximation property.
 Thus according to \cite[2.6]{Art}, there exists
 a common Nisnevich neighborhood $W_{(z)}$ of $(X_{(z)},z)$
 and $(Y_{(z)},z)$. Moreover, one can lift the situation
 in a neighborhood of $z$, both in $X$ and $Y$:
 there exists an $S$-scheme $W$ essentially of finite type,
  which lifts $W_{(z)}$ and fits into the following commutative
 diagram:
$$
\xymatrix@=12pt{
& W\ar^g[rd]\ar_f[ld] & \\
X'\ar@{^(->}[d] & Z'\ar@{^(->}[d]\ar[r]\ar[l]\ar[u] & Y'\ar@{^(->}[d] \\
X & Z\ar[r]\ar[l] & Y
}
$$
where the squares are cartesian, $X'$, $Y'$, $Z'$ are open neighborhood
of $z$ in $X$, $Y$, $Z$ respectively, and $f$ (resp. $g$) is a Nisnevich neighborhood
 of $Z'$ in $X$ (resp. $Y$).

Using this construction, one can further find Zariski hypercoverings
 $\mathcal X$, $\mathcal Y$, $\mathcal Z$ of $X$, $Y$, $Z$
 and a simplicial scheme $\mathcal W$
 which fits into the following commutative diagram:
\begin{equation}
\begin{split}
\label{eq:analytic_inv}
\xymatrix@=12pt{
& \mathcal W\ar^g[rd]\ar_f[ld] & \\
\mathcal X\ar_p[d] & \mathcal Z\ar[d]\ar[r]\ar[l]\ar[u]
 & \mathcal Y\ar^q[d] \\
X & Z\ar[r]\ar[l] & Y
}
\end{split}
\end{equation}
such that the squares are cartesian and for each integer $n \geq 0$,
 $\mathcal W_n$ is a Nisnevich neighborhood of $\mathcal Z_n$
 in $\mathcal X_n$ (resp. $\mathcal Y_n$).

To finish the proof, one has to use the fact that the stable
 homotopy category $\SH(S)$ can be extended to simplicial schemes 
 according to \cite[Chap. 4]{Ayoub2} or \cite[3.1]{CD3}.
 This implies that the cohomology $\E^{**}$ can be extended
 to simplicial $S$-schemes and simplicial $Z$-pairs.
 Then the following sequence gives us almost the desired isomorphism
$$
\mathfrak f^*_\mathcal W:\E^{**}_Z(X) \xrightarrow{p^*} \E^{**}_{\mathcal Z}(\mathcal X)
 \xrightarrow{f^*} \E^{**}_{\mathcal Z}(\mathcal W)
 \xrightarrow{(g^*)^{-1}}\E^{**}_{\mathcal Z}(\mathcal Y)
 \xrightarrow{(q^*)^{-1}} \E^{**}_Z(Y).
$$
Indeed, according to the Zariski descent property of $\SH$
 as formulated in \cite[3.2.7, 3.3.5]{CD3}, the maps $p^*$ and $q^*$ are isomorphisms.
 Moreover, we can derive from property (Nis) the fact the 
 maps $f^*$ and $g^*$ are isomorphisms (either we apply the Zariski
 descent spectral sequence or we argue directly in $\SH$).

To finally get $\mathfrak f^*$, one has to take
 the limit of the isomorphisms $\mathfrak f^*_\mathcal W$ over the filtering
 category of diagrams of the form \eqref{eq:analytic_inv}.
\end{proof}

\begin{rem}
\begin{enumerate}
\item The fact $\mathfrak f^*$ depends only on $\mathfrak f$
 can also be supplemented by the following cocycle condition:
 given composable isomorphisms $\mathfrak f$ and $\mathfrak g$ of certain formal completions,
 one gets: $(\mathfrak f \circ \mathfrak g)^*=\mathfrak g^* \circ \mathfrak f^*$.
\item Under the presence of a ring structure on $\E$,
 the isomorphism $\mathfrak f^*$ is compatible with
 products as defined in Paragraph \ref{num:products}.
\end{enumerate}
\end{rem}

As a corollary, we get the following reinforcement of 
 Prop. \ref{prop:inv_completion}:
\begin{cor}\label{cor:inv_completion}
Let $\E$ be an absolute spectrum.
Let $X$ be a local scheme with closed point $x$,
 $\hat X$ be its completion at the point $x$
 and consider the canonical map $f:(X,x) \rightarrow (\hat X,x)$.
Then the induced morphism of cohomology with support:
$$
f^*:\E^{**}_x(\hat X) \rightarrow \E^{**}_x(X)
$$
is an isomorphism.
\end{cor}

\begin{ex}
The preceding theorem and its corollary can be applied
 in particular to any of the absolute ring spectrum
 of Example \ref{ex:abs_ring_spectra}.
 In particular we get another proof of analytical invariance
 for De Rham cohomology and a proof in the case of
 rigid cohomology. 
 The later case was also proved independently
 by Ouwehand in \cite{Ouw}.

Note that it also holds for $K$-theory
 and algebraic cobordism. Even in the first case,
 this seems to be new.
\end{ex}

\begin{rem}
The preceding corollary is especially useful in dealing
 with absolute purity. In order to prove it for a given
 absolute ring spectrum, one easily reduces to the case of
 a closed pair $(X,x)$ where $X$ is a local regular scheme
 with closed point $x$. According to the preceding corollary,
 one derives that we can further assume that $X$ is the
 spectrum of a complete local scheme.
\end{rem}

\section{Orientation and characteristic classes}
\label{sec:orientation}

\subsection{Orientation theory and Chern classes}
\label{sec:orientation&Chern}

The considerations of this section are well known in motivic homotopy theory
 (see for example \cite{Bor, Vez} -- in the case of a base field).
 They were also studied, with a slightly different formalism,
 in our paper \cite{Deg8}.

\begin{num} \label{num:Pinfty}
Let $S$ be a scheme.
We will assume that the scheme $\PP^n_S$ is pointed by the infinite point
 (of homogeneous coordinates $[0:...:0:1]$ to fix ideas).
 Then we get a tower of pointed $S$-schemes
$$
\PP^1_S \rightarrow \PP^2_S \rightarrow \hdots \PP^n_S \xrightarrow{\iota_n} \PP^{n+1}_S
 \rightarrow \hdots
$$
where $\iota_n$ denotes the embedding of the last $n$-th coordinates.
The colimit of this tower in the category of pointed sheaves defines an object 
 $\PP^\infty_S$
 of the pointed homotopy category $\H(S)$ -- see \cite{MV}.
 We still denote by $\iota_1:\PP^1_S \rightarrow \PP^\infty_S$ 
 the induced map in the homotopy category.

The following definition is now basic in motivic homotopy theory:
\end{num}
\begin{df} \label{df:orientation}
Let $\E$ be an absolute ring spectrum
 with unit $\eta_S:S^0 \rightarrow \E_S$ over a scheme $S$ in $\base$.
 We can see $\eta_S$ as a class in the reduced cohomology
 $\tilde \E^{2,1}(\PP^1_S)$.\footnote{By definition of the Tate twist,
  $\E^{2,1}(\PP^1_S)=\E^{2,1}(S) \oplus \E^{0,0}(S)$.}

An \emph{orientation} of $\E$ over $S$ is a class $c_{S}$
 in the reduced cohomology
 $\tilde \E^{2,1}(\PP^\infty_S)$ such that $\iota_1^*(c_S)=\eta_S$.
 
An (absolute) \emph{orientation} of $\E$ is a family of
 classes $c=(c_{S})$ for any scheme $S$ in $\base$
 such that for any morphism $f:T \rightarrow S$,
 $f^*(c_{S})=c_{T}$.
In this situation, we also say that $(\E,c)$ is an 
 absolute oriented ring spectrum.
\end{df}

\begin{rem}
To give an orientation of an $S$-absolute ring spectrum $\E$
 it is enough and sufficient to give an orientation of $\E_S$.
\end{rem}

\begin{ex}\label{ex:abs_ring_spectra_orientation}
Let us review the absolute ring spectra
 of Example \ref{ex:abs_ring_spectra}:
\begin{enumerate}
\item According to \cite[2.2.8]{CD2}, the $k$-absolute ring spectrum
 associated with a Mixed Weil theory is canonically oriented.
\item Voevodsky's Einlenberg-Mac Lane motivic ring spectrum
 $\HH_\Lambda$,
 considered as an $S$-absolute ring spectrum,
 is oriented according to \cite[Sec. 11.3]{CD3}.
\item The cobordism ring spectrum $\MGL$ is canonically oriented.
 This follows from theoretical reasons
 (see recall in Prop. \ref{prop:orientation&MGL-modules})
  or can be directly seen from the construction
	 (see Par. \ref{num:orientation_MGL}).
\item The K-theory ring spectrum $\KGL$ is oriented.
 Let $S$ be a regular scheme ($S=\spec{\ZZ}$ would be enough).
By construction, we have a canonical isomorphism:
$$
K_0(S) \simeq \KGL^{0,0}(S).
$$
In particular, any line bundle $L/S$ defines an element $[L] \in \KGL^{0,0}(S)$.
Moreover, Bott periodicity theorem implies the existence of the following isomorphism:
$$
\KGL^{0,0}(S) \simeq K_0(S) \simeq \tilde K_0(\PP^1_S,\infty)
\simeq \widetilde{\KGL}^{0,0}(\PP^1_S,\infty).
$$
The image of $1$ by this isomorphism is an element $\beta$ in $\KGL^{-2,-1}(S)$
 called the \emph{Bott element}.
 Let $\lambda$ be the canonical line bundle on $\PP^\infty_S$.
 Then we define the orientation of $\KGL$ as follows:\footnote{This choice of orientation coincides
 with the one of \cite[Ex. 1.1.5]{ML}.
 In the literature however, one can find different choices of orientations of the
 ring spectrum $\KGL$. The present choice is justified by Example
 \ref{ex:gysin_H&KGL} as well as the correct form of the Todd class appearing in the
 Riemann-Roch theorem \ref{prop:classical_GRR}.}  
$$
c^\KGL:=\beta^{-1}.(1-[\lambda^\vee]) \in \KGL^{2,1}(\PP^\infty_S).
$$
\item Beilinson motivic cohomology spectrum $\HB$ is canonically oriented
 (see \cite[14.1.5]{CD3}, \cite{NSO}).
\end{enumerate}
It is worth to point out that we will show in Cor. \ref{cor:uniqueness_orientation}
 and Ex. \ref{ex:uniqueness_orientation} that the orientations involved
 in points (1), (2) and (5) are unique.
\end{ex}

\begin{rem}
The presence of the Bott element in formulas involving the spectrum
 $\KGL$ can be explained as follows.
 By construction (see \cite{Riou}),
 for any integer $n$ and any regular scheme $S$,
 one has a canonical contravariantly functorial isomorphism:
$$
\varphi^{n}_S:\KGL^{n,0}(S) \rightarrow K_{-n}(S)
$$
where the right hand side is the $n$-th Quillen $K$-theory of $S$.
According to the definition above, multiplication by $\beta$ on $\KGL^{**}(S)$
 induces an isomorphism.
 Thus, for any couple of integers $(r,n)$, we get a canonical isomorphism:
\begin{equation} \label{eq:iso_KGL&Ktheory}
\varphi^{n,i}_S:\KGL^{n,i}(S)
 \rightarrow \KGL^{n-2i,0}(S)
 \rightarrow K_{2i-n}(S),
 x \mapsto \varphi^{n-2i}_S(\beta^i.x).
\end{equation}
\end{rem}

\begin{num}
Recall from \cite[\S 4, 3.7]{MV} there exists a canonical isomorphism
 in $\H(S)$:
$$
\mathrm B \GG \simeq \PP^\infty_S
$$
where $\mathrm B\GG$ is the (Nisnevich) classifying space of $\GG$. 
This immediately gives an application
\begin{equation} \label{eq:Pic->H}
\Pic(S)\stackrel{(*)}=\Hom_{\H^s(S)}(S_+,\mathrm B\GG)
 \rightarrow \Hom_{\H(S)}(S_+,\mathrm B\GG)
 \simeq \Hom_{\H(S)}(S_+,\PP^\infty_S)
\end{equation}
where $\H^s(S)$ denotes the simplicial homotopy category
 and the first map is induced by the projection functor
 $\H^s(S) \rightarrow \H(S)$ -- the target category being
 the $\AA^1$-localization of the source category. We have used
 \cite[\S 4, 1.15]{MV} for the identification $(*)$.\footnote{
 Note also Morel and Voevodsky proved the map \eqref{eq:Pic->H} is
 an isomorphism
 whenever $S$ is regular ; \emph{op. cit.} Prop. 3.8.}
\end{num}

\begin{rem}
For any integer $n$, we let $\lambda_n$ be the canonical
 line bundle on $\PP^n_S$ (see Notations and conventions,
 p.~\pageref{sec:notations}).
Then the family $(\lambda_n)_{n \in \NN}$
 defines an element $\lambda$ of $\Pic(\PP^\infty_S)$
 -- which generates this group as a ring of formal power series
 over $\ZZ$ as soon as $S$ is local and regular.
The map \eqref{eq:Pic->H} is characterized by the fact
 it sends $\lambda$ to the canonical projection
  $\PP^\infty_{S+} \rightarrow \PP^\infty_S$.
\end{rem}

\begin{df} \label{df:1st_Chern_class}
Let $(\E,c)$ be an absolute oriented ring spectrum.

For any scheme $S$, we associate to the class $c$ a
 canonical morphism of sets:
\begin{align*}
c_1:\Pic(S) \xrightarrow{\eqref{eq:Pic->H}}  \Hom_{\H(S)}(S_+,\PP^\infty_S)
 & \xrightarrow{\Sus} \Hom_{\SH(S)}(\Sus S_+,\Sus \PP^\infty_S) \\
&  \xrightarrow{(c_{S})_*} \Hom_{\SH(S)}(\Sus S_+,\E_S(1)[2])=\E^{2,1}(S),
\end{align*}
called the \emph{first Chern class}.
\end{df}

\begin{ex}
In the case of the orientation of $\KGL$ defined
 in Ex. \ref{ex:abs_ring_spectra_orientation},
 one gets for any line bundle over a regular scheme $S$:
$$
c^\KGL_1(L)=\beta^{-1}.(1-[L^\vee]).
$$
\end{ex}

\begin{num} \label{num:basic_chern}
According to this definition and the preceding remark,
 we get the following properties:
\begin{enumerate}
\item[(a)] For any morphism of schemes $f:T \rightarrow S$
 and any line bundle $\lambda$ on $S$,
 $f^*c_1(\lambda)=c_1(f^{-1}\lambda)$.
\item[(b)] Let $n \geq 0$ be an integer, $\lambda_n$ be the line bundle over $\PP^n_S$ considered in the above remark 
 and $\nu_n:\PP^n_S \rightarrow \PP^{\infty}_S$ be the obvious morphism.
 
Then $c_1(\lambda_n)=\nu_n^*(c_{S})$
 as classes in $\E^{2,1}(\PP^n_S)$, according to the above remark.
\end{enumerate}
\end{num}

\begin{rem}
One must be careful that the relation
 $c_1(\lambda \otimes \lambda')=c_1(\lambda)+c_1(\lambda')$ does not necessarily hold.
 This is due to the fact that the second of the three maps considered
 in the definition of $c_1$
 is not a morphism of abelian groups.
 This remark will be made more precise latter (see \ref{prop:FGL}). 
\end{rem}

\begin{num}
Next we recall the projective bundle theorem for
 an absolute oriented ring spectrum $(\E,c)$.
 
Let $p:P \rightarrow S$ be a projective bundle of rank $n$
 with canonical line bundle $\lambda$.
 We define the following morphism:
$$
\epsilon_P:\oplus_{i=0}^n \E^{**}(X)
 \rightarrow \E^{**}(P),
 (x_0,\hdots,x_n)\mapsto \sum_i p^*(x_i).c_1(\lambda)^i.
$$
\end{num}
\begin{thm} \label{thm:PBF}
With the above assumptions and notations,
 the morphism $\epsilon_P$ is an isomorphism.
\end{thm}
In other words, $\E^{**}(P)$ is a free graded $\E^{**}(X)$-module
 with basis $(1,c_1(\lambda),\hdots,c_1(\lambda)^n)$ (as usual).
\begin{proof}
The proof, essentially due to Morel,
 is the same as the proof of Th. 3.1 in \cite{Deg8}.
 We recall the main steps for the convenience of the reader.

Using the Mayer-Vietoris long exact sequence associated to
 an open cover by two open subsets, we reduce to the case
 where $P$ is trivializable and then to the case $P=\PP^n_S$
 -- this uses only property \ref{num:basic_chern}(a).

Then the proof goes on by induction on $n$ ; the case $n=0$ is
 trivial and the case $n=1$ is an immediate consequence of the definition
 of the orientation $c_{S}$.

The principle of the induction is to use the following facts:
\begin{itemize}
\item Let $\PP^n_S/\PP^{n-1}_S$ be the cokernel of the embedding $\iota_{n-1}$
 in the category of pointed sheaves.
 Then the sequence
$$
\PP^{n-1}_S \xrightarrow{\iota_n} \PP^n_S \xrightarrow{\pi_n} \PP^n_S/\PP^{n-1}_S
$$
is homotopy exact ; in particular, it induces a long exact sequence:
$$
\hdots \rightarrow \E^{**}(\PP^n_S/\PP^{n-1}_S) \xrightarrow{\pi_n^*} \E^{**}(\PP^n_S)
 \xrightarrow{\iota_n^*} \E^{**}(\PP^{n-1}_S)
 \rightarrow \hdots
$$
\item There exists a canonical isomorphisms in $\H(S)$:
$$
\tau_n:\PP^n_S/\PP^{n-1}_S \rightarrow (\PP^1_S)^{\wedge,n}.
$$
\end{itemize}
Then we are reduce to prove the following relations:
\begin{itemize}
\item $c_1(\lambda_{n-1})^n=0$.
\item Using the isomorphism $\tau_n$,
 we get an isomorphism $\tau_n^*:\E^{2n,n}(\PP^n_S/\PP^{n-1}_S) \simeq \E^{0,0}(S)$.
 Then $\tau_n^*(\overline{c_1(\lambda_n)^n})=\eta_S$, the unit of the ring
 spectrum $\E_S$ --- where we have put:
 $\pi_n^*(\overline{c_1(\lambda_n)^n})=c_1(\lambda_n)^n$.
\end{itemize}
These relations can easily be deduced from the following lemma:
\begin{lm}[Morel] \label{lm:Morel}
Let $\delta_n:\PP^n_S \rightarrow (\PP^n_S)^{\wdg,n}$
 be the $n$-th diagonal
 of the pointed scheme $\PP^n_S/S$. Then the following square commutes in $\H(S)$:
$$
\xymatrix{
\PP^n_S\ar_{\pi_n}[d]\ar^-{\delta_n}[r] & (\PP^n_S)^{\wedge,n} \\
\PP^n_S/\PP^{n-1}_S\ar^{\tau_n}[r] & (\PP^1_S)^{\wedge,n}.\ar_{(\iota_1)^{\wedge,n}}[u]
}
$$
\end{lm}
For this lemma, we refer the reader to the proof of \cite[lem. 3.3]{Deg8}.
\end{proof}

As first remarked by Morel,
 the projective bundle theorem admits the following corollary
  -- whose proof can be easily adapted from \cite[Cor. 3.6]{Deg8}.
\begin{cor}
Let $\E$ be an orientable absolute ring spectrum.
Then for any scheme $X$ in $\base$ and any closed subschemes $Z$, $Z'$ of $X$,
 one has the following property:
$$
\forall (x,y) \in \E^{n,p}_Z(X) \times \E^{m,q}_{Z'}(X), \
 x \cupp y=(-1)^{nm}y \cupp x.
$$
\end{cor}

Following the method of Grothendieck,
 we can now introduce the following definition.
\begin{df} \label{df:Chern_classes}
Let $(\E,c)$ be an absolute oriented ring spectrum.

Let $E/S$ be a vector bundle of rank $n$.
We let $P=\PP(E)$ be the associated projective bundle,
 with projection $p$ and canonical line bundle $\lambda$.

Using the previous theorem,
 we define the \emph{Chern classes of $E/S$} with coefficients
 in $(\E,c)$ as the elements $c_i(E)$ of $\E^{2i,i}(S)$ 
 for an integer $i \in [0,n]$ such that
\begin{equation} \label{eq:Chern}
\sum_{i=0}^n p^*(c_i(E)).\big(-c_1(\lambda)\big)^{n-i}=0
\end{equation}
and $c_0(E)=1$.
We put $c_i(E)=0$ for $i \notin [0,n]$.
\end{df}
Indeed, the above theorem guarantees the existence and uniqueness
 of the Chern classes $c_i(E)$.
 
\begin{num} \label{num:prop_Chern}
We deduce from this definition the following (usual) properties of Chern classes:
\begin{enumerate}
\item[(a)] For any vector bundle $E/S$, and any morphism $f:T \rightarrow S$,
 $f^*c_i(E)=c_i(f^{-1}E)$.
\item[(b)] For any scheme $S$ and any isomorphism of vector bundles $E \simeq E'$
 over $S$, $c_i(E)=c_i(E')$.
\item[(c)] For any scheme $S$, any exact sequence of vector bundles:
$$
0 \rightarrow E' \rightarrow E \rightarrow E'' \rightarrow 0
$$
and integer $k \geq 0$, $c_k(E)=\sum_{i+j=k} c_i(E').c_j(E'')$.
\end{enumerate}
For details on the proof of formula (c) -- \emph{Whitney sum formula} --
 we refer the reader to \cite[3.13]{Deg8}.
\end{num}

\begin{rem}
Let $(\E,c)$ be an absolute oriented ring spectrum.

Let $Gr_S$ be the infinite Grassmannian, seen
 as a Nisnevich simplicial sheaf of sets over $\sm_S$.
  According to \cite[4.3.7]{MV}, it is isomorphic
  in $\H(S)$ to the classifying space
  $BGL_S$ of the infinite general linear group over $S$.
Moreover, when $S$ is regular and for any integer $n \geq 0$,
one gets a canonical isomorphism:
$$
\Hom_{\H(S)}(S^n,\ZZ \times Gr) \simeq K_n(S),
$$
according to \cite[4.3.13]{MV}, where $\ZZ \times Gr$
 is the product of $\ZZ$-copies of $Gr_S$.
 This map is compatible with pullbacks of regular schemes.

Using Chern classes,
 it is well known how to compute $\E^{**}(Gr_S)$ (see \cite[6.2]{NOS}):
$$
\E^{**}(Gr_S) \simeq \E^{**}(S)[[\mathfrak c_1,\hdots,\mathfrak c_n,\hdots]]
$$
where $\mathfrak c_i$ is a cohomology class of bidegree $(2i,i)$.
According to the above isomorphism, it corresponds to a map
in $\H(S)$:
$$
\mathfrak c_i:\ZZ \times Gr_S \rightarrow Gr_S \rightarrow \Omega^{\infty}(\E_S(i)[2i]).
$$

Moreover, using the techniques of Riou
 (cf. \cite{Riou}, Th. 1.1.6 as in the proof of Th. 6.2.1.2),
 one gets an isomorphism:
$$
\Hom_{\H(S)}(\ZZ \times Gr,\Omega^\infty(\E(i)[2i])
 \rightarrow \Hom(K_0(-),\E^{2i,i}(-))
$$
where the right hand side stands for the
 morphisms of presheaves of sets on $\sm_S$.
 Under this isomorphism, the map $\mathfrak c_i$
 corresponds to the natural transformation $c_i$
 that we have just defined.

As in \cite{Gillet},
 this allows to automatically extends Chern classes
 to higher Chern classes with support in a closed
 subscheme $Z \subset S$ as follows:
\begin{align*}
c_i^Z:K_n^{Z}(S)
 \simeq [S^n \wdg (S/S-Z),\ZZ \times Gr]
 &\xrightarrow{(\mathfrak c_i)_*}
  [S^n \wdg (S/S-Z),\Omega^\infty(\E(i)[2i])] \\
 &\simeq \E^{2i-n,i}_Z(S).
\end{align*}
\end{rem}

\begin{num} \label{num:FGL}
\textit{The associated formal group law.--}
The usual Segre embeddings
$$
\sigma_{nm}:\PP^n_X \times_X \PP^m_X \rightarrow \PP^{n+m+nm}_X,
$$
indexed by a pair of positive integers $(n,m)$,
induce a multiplication map
\begin{equation} \label{eq:segre}
\sigma:\PP^\infty_X \times_X \PP^\infty_X \rightarrow \PP^\infty_X.
\end{equation}
This gives a structure of an H-group to the object $\PP^\infty_X$
 of $\H(X)$ which in turn induces a group structure on the target
 of the map \eqref{eq:Pic->H}.
 We deduce from the previous paragraph and the equality
\begin{equation} \label{eq:segre_pullback}
\sigma_{nm}^{-1}(\lambda_{n+m+nm})=\lambda_n \times_X \lambda_m
\end{equation}
that \eqref{eq:Pic->H} is a morphism of abelian groups.

Suppose $(\E,c)$ is an absolute oriented ring spectrum.
According to Theorem \ref{thm:PBF},
 the pullback along $\sigma$ corresponds to a (comultiplication) morphism:
$$
\E^{**}(X)[[c]] \xrightarrow{\sigma^*} \E^{**}(X)[[x,y]]
$$
where $x$ and $y$ stands for the Chern classes of the two
canonical line bundles $p_1^{-1}(\lambda)$ and $p_2^{-1}(\lambda)$.
It follows that $\sigma^*$ is determined by $\sigma^*(c)$
 which is a power series of the form:
\begin{equation} \label{eq:FGL}
F_X(x,y)=\sum_{i,j \geq 0} a^X_{ij}.x^iy^j.
\end{equation}
As $\sigma^*$ is a comultiplication,
 one deduces that $F_X$ is a commutative formal group law
 with coefficients in $\E^{**}(X)$ -- see \cite[\textsection 3.7]{Deg8}.
 In fact, the classes $a_{ij}^S$ enjoy the following properties:
\begin{itemize}
\item $a^S_{ij}$ has bidegree $(2-2i-2j,1-i-j)$ in $\E^{**}(S)$,
\item $a^S_{0,1}=1$, $a^S_{0,i}=0$ if $i \neq 1$,
\item for every couple $(i,j)$, $a^S_{ij}=a^S_{ji}$,
\item for any morphism of schemes $f:T \rightarrow S$,
 $f^*(a^S_{ij})=a^T_{ij}$.
\end{itemize}
\end{num}

\begin{df} \label{df:fgl}
Given the notations above,
 we will say that $F_S$ is the \emph{formal group law} associated
 with the oriented ring spectrum $(\E,c)$ above $S$.

We will say that $(\E,c)$ (or just $c$) is \emph{additive}
 (resp. \emph{multiplicative} with parameter $u$)
 if for any scheme $S$, $F_S(x,y)=x+y$
 (resp. $F_S(x,y)=x+y+u.x.y$).
\end{df}

\begin{ex} \label{ex:abs_ring_spectra_FGL}
Consider the absolute oriented ring spectra
 of Example \ref{ex:abs_ring_spectra_orientation}:
\begin{enumerate}
 \item The $k$-absolute oriented ring spectrum associated
 with a Mixed Weil theory is additive (cf. \cite[2.2.10]{CD2}).

The absolute oriented ring spectrum $\HB$
 is also additive.
 This last fact follows from the definition. To be more precise,
 given any regular scheme $S$, the canonical bijection:
$$
\HB^{2,1}(S)\simeq Gr_\gamma^1 K_0(S)_\QQ \simeq Pic(S)_\QQ
$$
is in fact an isomorphism of abelian groups.
Note also that, restricting to the category of smooth $S$-schemes,
 for an arbitrary base $S$,
 the absolute oriented ring spectrum $\HH_\Lambda$ is additive
 according to \cite[11.3.5]{CD3}.
\item The absolute oriented ring spectrum $\KGL$ is multiplicative
 with parameter $(-\beta)$
 as follows from the easy computation,
 with $l=[L^\vee]$ and $l'=[L^{\prime\vee}]$:
\begin{align*}
c_1^\KGL(L \otimes L')&=\beta^{-1}.(1-[(L \otimes L')^\vee])
=\beta^{-1}.(1-ll') \\
&=\beta ^{-1}((1-l)+(1-l')-(1-l).(1-l')) \\
&=c_1^\KGL(L)+c_1^\KGL(L')-\beta.c_1^\KGL(L).c_1^\KGL(L')
\end{align*}
\item Let $k$ be a field of exponential characteristic $p$.
 Let $\MGL[1/p]$ be the absolute ring spectrum obtained from
 $\MGL$ by inverting $p$.
Then, as a consequence of the Theorem of Hopkins-Morel-Hoyois
 (cf. \cite[Th. 7.12]{Hoy}), we know that
 the formal group law of the absolute oriented ring spectrum
 $\MGL[1/p]$, considered over the category of all $k$-schemes,
 is isomorphic to the universal formal group law.

More precisely, if $(L,F_{univ})$ denotes the Lazard ring
 equipped with its canonical formal group law,
 according to \emph{loc. cit.}, Prop. 8.2,
 there exists an isomorphism of formal group laws:
$$
(L[1/p],F_{univ}) \rightarrow (\MGL_{(2,1)*}(k)[1/p],F_\MGL).
$$
\end{enumerate}
\end{ex}

\begin{prop} \label{prop:FGL}
Consider the notations of the previous definition.
\begin{enumerate}
\item For any vector bundle $E/X$
 and any integer $i>0$, the class $c_i(E)$ is nilpotent in $\E^{**}(X)$.
\item  For any line bundles $L_1,L_2$ over $X$,
$$
c_1(L_1 \otimes L_2)=F_X\big(c_1(L_1),c_1(L_2)\big) \in \E^{2,1}(X).
$$
\end{enumerate}
\end{prop}
\begin{proof}
Point (1) follows from the hypothesis that $X$ is noetherian
 according to the proof of \cite[3.8(1)]{Deg8}.
 Point (2) is tautological by definition of the
 first Chern class $c_1$ and of the formal group law $F_X$. 
\end{proof}

\begin{rem}
Note that unlike in \cite[Prop. 3.8]{Deg8},
 to prove Point (2),
 we do not need that $X$ admits an ample line bundle. 
 This is because we consider cohomology theories that
 are representable in $\SH$: in fact,
 any line bundle $L/X$ can be represented
 by a map $S^0 \rightarrow \PP^\infty_X$
 in $\H(X)$ according to the theorem of Morel-Voevodsky.
 (If $L/X$ is not generated by its sections, 
 this map cannot be lifted in the category of schemes.)
\end{rem}

\subsection{Thom classes and $\MGL$-modules}

\begin{num} \label{num:Thom_space}
Let $(\E,c)$ be an absolute oriented ring spectrum 
 (see Def. \ref{df:orientation}).

Let $E/X$ be a vector bundle of rank $n$, $\PP(E)$ 
(resp. $\bar E$, $E^\times$) be the associated projective bundle
 (resp. projective completion, complement of the zero section).
 Recall from \cite[\S 3, 2.16]{MV} that one defines the Thom space
 of $E/X$ as the pointed sheaf
$$
\Th(E)=E/E^\times.
$$
According to \emph{loc. cit.}, Prop. 2.17, we get a canonical $\AA^1$-equivalence of simplicial sheaves
 $\Th(E) \simeq \bar E/\PP(E)$,
 the right hand side being the cokernel of the canonical embedding
  $\nu:\PP(E) \rightarrow \bar E$
 in the category of sheaves, equipped with its obvious base point.
 Thus we get a homotopy cofiber sequence for the $\AA^1$-local model structure:
$$
\PP(E) \xrightarrow \nu \bar E \xrightarrow \pi \Th(E)
$$
which induces a long exact sequence
$$
\hdots \rightarrow \E^{**}(\Th(E)) \xrightarrow{\pi^*} \E^{**}(\bar E)
 \xrightarrow{\nu^*} \E^{**}(\PP(E)) \rightarrow \hdots
$$
According to Theorem \ref{thm:PBF}, $\nu^*$ is a split epimorphism
 of free $\E^{**}(X)$-modules of respective ranks $n$ and $n-1$.
 Thus $E^{**}(\Th(X))$
 is a free $\E^{**}(X)$-module of rank $1$, isomorphic to $\ker(\nu^*)$.
\end{num}
\begin{df} \label{df:Thom}
Consider the notations and assumptions above.

We define the \emph{Thom class} of $E/X$ as
 the following element of $\E^{2n,n}(\bar E)$:
$$
\thom(E)=\sum_{i=0}^n p^*(c_i(E)).\big(-c_1(\lambda)\big)^{n-i}
$$
where $p$ is the canonical projection.
We define the \emph{refined Thom class} $\rthom(E)$ of $E$
 as the unique element of $\E^{2n,n}(\Th(E))$ such that
$$
\pi^*\big(\rthom(E)\big)=\thom(E).
$$
When the base of the vector bundle $E$ is not clear,
 we indicate it as follows: $\PP_X(E)$, $\Th_X(E)$, $\thom(E/X)$, $\rthom(E/X)$.

Note that $\E^{**}(\Th(E))=\E^{**}_X(E)$ (see Definition \ref{df:coh_support}).
 According to the preceding paragraph, the canonical map:
\begin{equation} \label{eq:thom_iso}
\E^{**}(X) \rightarrow \E^{**}_X(E), \lambda \mapsto \lambda.\rthom(E)
\end{equation}
defined in Paragraph \ref{num:products} is an isomorphism
 called the \emph{Thom isomorphism}.
\end{df}
We deduce from formulas (a), (b) of \ref{num:prop_Chern} 
 that Thom classes are compatible with base change
 and invariant under isomorphisms of vector bundles.

\begin{rem} \label{rem:non_constant_rk}
Recall that in general the rank of a vector bundle $E/X$
 is Zariski locally constant on $X$. In other
  words, it is a function $r:\pi_0(X) \rightarrow \NN$.
On the other hand, $\E$-cohomology of $X$ is additive.
 Moreover, we can give sense to the formula defining the refined Thom class
 of $\E$ without requiring $E/X$ is of (constant) rank equal to $n$.
 Then $\rthom(E)$ is an element of $\E^{**}(\Th(E))$,
  which still induces a Thom isomorphism as above.
Then, we can define the bidegree of this class
 as the function $\pi_0(X) \rightarrow \ZZ^2$
 which to a connected component $X_i$ of $X$ associates
 the couple $(2r(X_i),r(X_i))$.
The Thom isomorphism is then Zariski locally homogeneous on $X$
 with the same bidegree.
\end{rem}

\begin{ex} \label{ex:thom&euler_quotient_bdl}
Recall the universal quotient bundle $\xi$ on $\bar E$ is defined by 
the exact sequence
$$
0 \rightarrow \lambda \rightarrow p^{-1}(E \oplus 1) \rightarrow \xi \rightarrow 0.
$$
Thus the Whitney sum formula \ref{num:prop_Chern}(c)
gives the following formula:
\begin{equation}\label{eq:compute_Thom}
\thom(E)=c_n(\xi).
\end{equation}
\end{ex}

\begin{num} \label{num:orientation_MGL}
The cobordism spectrum $\MGL$
 is given by the sequence of Thom spaces $\Th(\gamma_n/\mathrm{BGl}_n)$ for $n>0$
 where $\gamma_n$
 is the tautological rank $n$ vector bundle over the classifying
 space of $\mathrm{GL}_n$.
As $\Th(\gamma_1)=\mathrm B \GG=\PP^\infty$,
 $\MGL_S$ is canonically oriented.
 We denote by $c^\MGL$ this orientation.
 
Given an absolute oriented ring spectrum $(\E,c)$,
 the Thom class defined previously allows to define a morphism 
 $\varphi_c:\MGL \rightarrow \E$. This is the key observation
 of the following proposition (see \cite[4.3]{Vez}):
\end{num}
\begin{prop} \label{prop:orientation&MGL-modules}
Let $\E$ be an absolute ring spectrum. Then the following sets correspond
 bijectively:
\begin{enumerate}
\item[(i)] orientations $c$ of $\E$;
\item[(ii)] maps of absolute ring spectra $\varphi:\MGL \rightarrow \E$.
\end{enumerate}
by the following maps:
\begin{align*}
(ii) &\rightarrow (i), \varphi \mapsto \varphi_*(c^\MGL),
 \ \varphi_* \text{ induced map in cohomology,} \\
(i) &\rightarrow (ii), c \mapsto \varphi_c.
\end{align*}
\end{prop}

\begin{num}
A module over the ring spectrum $\MGL_S$ is an $\MGL_S$-module in $\SH(S)$
 in the classical sense: a spectrum $\E$ over $S$ equipped
 with a multiplication map $\gamma_\E:\MGL_S \wedge \E \rightarrow \E$
 satisfying the usual identities -- see \cite{Mcl}.

Given two $\MGL_S$-modules $\E$ and $\F$,
 a morphism of $\MGL_S$-modules is a morphism $f:\E \rightarrow \F$ in $\SH(S)$
 such that the following diagram is commutative:
$$
\xymatrix@R=12pt@C=28pt{
\MGL_S \wedge \E\ar^{1 \wedge f}[r]\ar_{\gamma_\E}[d]
 & \MGL_S \wedge \F\ar^{\gamma_\F}[d] \\
\E\ar^f[r] & \F.
}
$$
We will denote by $\wmod \MGL_S$ the additive category of $\MGL_S$-modules.

Given any spectrum $\E$, $\MGL_S \wedge \E$ has an obvious structure
 of an $\MGL_S$-module. The assignment $L_\MGL^w:\E \mapsto \MGL_S \wedge \E$ 
 defines a functor left adjoint to the inclusion functor $\mathcal O^w_\MGL$
 and we get an adjunction of categories:
\begin{equation} \label{eq:SH_wmod}
L_\MGL^w:\SH(S) \leftrightarrows \wmod \MGL_S:\mathcal O_\MGL^w.
\end{equation}
The category $\wmod \MGL_S$ is not well behaved: it has no triangulated
 monoidal structure. In \cite[Ex. 2.12(2)]{Deg8} and \cite[Sec. 2.2]{Deg10},
 we introduced the homotopy category of \emph{strict} $\MGL_S$-modules.
 We will denote it by $\smod \MGL_S$ and call it
  the \emph{homotopy category of $\MGL_S$-modules}.
This is an enrichment of the category of $\MGL_S$-modules: one defines
a canonical commutative diagram of functors:
\begin{equation}
\begin{split}\label{eq:weak_MGL-mod_diag}
\xymatrix@=14pt{
& \SH(S)\ar^{L^w_\MGL}[rd]\ar_{L_\MGL}[ld] & \\
\smod \MGL_S\ar|{\mathcal O'_\MGL}[rr] && \wmod \MGL_S
}
\end{split}
\end{equation}
where $L_\MGL$ is a triangulated functor
 and $\mathcal O'_\MGL$ is a conservative functor.
Therefore, using the previous commutative diagram and
 Proposition \ref{prop:orientation&MGL-modules},
 we get the following result.
\end{num}

\begin{prop}
Let $(\E,c)$ be an absolute oriented ring spectrum.

Then for any scheme $S$,
 the functor $\varphi_\E=\Hom_{\SH(S)}(-,\E)$ induces a canonical
 functor $\tilde \varphi_\E$ which fits in the following commutative 
 diagram:
\begin{equation*}
\xymatrix@R=4pt@C=40pt{
\SH(S)^{op}\ar^/8pt/{\varphi_\E}[rd]\ar_{(L_\MGL)^{op}}[dd] &  \\
& \mathscr Ab \\
(\smod \MGL_S)^{op}\ar_/10pt/{\tilde \varphi_\E}[ru] &
}
\end{equation*}
\end{prop}

\begin{num} \label{num:apply_Deg8}
Let $S$ be a scheme. Recall the considerations of \cite[\S 2.3.2]{Deg8}.
For any cartesian square of smooth $S$-schemes
$$
\xymatrix@=12pt{
W\ar[d]\ar[r]\ar@{}|\Delta[rd] & V\ar[d] \\
U\ar[r] & X
}
$$
made of immersions,
 we denote by $\frac{X/U}{V/W}$ the colimit of $\Delta$
 in the category of sheaves over $S$.\footnote{As all maps in this diagram
 are cofibrations, this is the homotopy colimit of the diagram in the
 $\AA^1$-local model category of simplicial sheaves. It measures
 the obstruction for the square $\Delta$ to be homotopy cocartesian.}

We denote by $\MGL_S\left(\frac{X/U}{V/W}\right)$ the image of the infinite suspension
 spectrum $\Sus\!\left(\frac{X/U}{V/W}\right)$ by the functor $L_\MGL$ defined above.
When $V=W=\varnothing$ (resp. $U=V=W=\varnothing$), we simply denote this
object by $\MGL_S\left(X/U\right)$ (resp. $\MGL_S(X)$)
 -- according to the notation of Par. \ref{num:i_*_monoidal}.

As explained in \cite[2.3.2, Ex. 2.12(2)]{Deg8},
 the category $\smod \MGL_S$ together with the canonical orientation
 $c^\MGL$ of $\MGL_S$ satisfies all the axioms of \cite[\S 2.1]{Deg8}.
 Therefore we can apply the results of \emph{loc. cit.}
 to that category.

Then, according to the previous proposition, for any smooth closed $S$-pair $(X,Z)$,
 we get:
$$
\E_Z^{n,m}(X)=\Hom_{\SH(S)}(X/X-Z,\E(m)[n])
 =\tilde \varphi_\E\left(\MGL_S\left(X/X-Z\right)(-m)[-n]\right).
$$
The important thing for us is that,
 given smooth closed $S$-pairs $(X,Z)$ and $(Y,T)$ and a couple of integers $(n,m)$,
 any morphism of $\MGL_S$-modules of the form:
$$
\MGL_S(X/X-Z) \rightarrow \MGL_S(Y/Y-T)(m)[n]
$$
induces a canonical homogeneous morphism of bigraded abelian groups
$$
\E^{**}_T(Y) \rightarrow \E_Z^{*+n,*+m}(X).
$$
Obviously, this association is compatible with composition.
\end{num}

\begin{ex} \label{ex:Gysin&MGL-mod}
Let $f:Y \rightarrow X$ be a projective $S$-morphism between smooth
 $S$-schemes. Assume $f$ has dimension $d$.
 Then, applying \cite[5.12]{Deg8}, we get the Gysin morphism
$$
f^*:\MGL_S(X) \rightarrow \MGL_S(Y)(-d)[-2d]
$$
which in turn induces a push-forward in cohomology:
$$
f_*:\E^{*,*}(Y) \rightarrow \E^{*-d,*-2d}(X).
$$
The purpose of the next sections is to generalize this
 pushforward in the case of regular schemes.
\end{ex}

\begin{rem}
In fact, all the orientation theory exposed here for ring spectra
 can be done without ring structure by replacing
 an orientation $c$ by a structure of a module over $\MGL_S$,
  on a given spectrum $\E$.
This is the content of \emph{loc. cit.} in the case
$\base=\sm_S$.\footnote{The case of weak $\MGL$-modules can be deduced
from the case of strict $\MGL$-modules which is treated in \emph{loc. cit.}
using the commutative diagram \eqref{eq:weak_MGL-mod_diag}.}

Then, instead of having Chern classes in $\E$ cohomology,
 we get an \emph{action} of Chern classes,
 through the $\MGL^{**}$-module structure of $\E^{**}$
  -- see also Example \ref{ex:RR_source_MGL}.
Similarly, to anticipate what follows,
 we do not get fundamental classes in $\E$-cohomology
 but an action of fundamental classes with coefficients in $\MGL$.
 Note this is enough to get the purity isomorphism \eqref{eq:purity}
 and therefore, Gysin long exact sequence (as in Def. \ref{df:Gysin_triangle}),
 Gysin morphisms for projective lci morphisms (as in \ref{df:Gysin}).

However, the general case $\base=\reg$ is not covered by \cite{Deg8}.
In fact, in this case,
 Conjecture \ref{conj:purityMGL} (of section \ref{sec:abs_purity})
 seems very natural.
\end{rem}

\subsection{Fundamental classes}
\label{sec:fdl_class}

In this section and the following one,
 we fix an absolute oriented ring spectrum $(\E,c)$.
\begin{df} \label{df:fdl_class}
Let $(X,Z)$ be a regular closed pair of codimension $n$
 and normal bundle $N_ZX$. 
 Let $i:Z \rightarrow X$ be the corresponding embedding.
 
When $(X,Z)$ is $\E$-pure
 we define the \emph{refined fundamental class} of $Z$ in $X$
  (with coefficients in $\E$)
 as the image of the Thom class $\rthom(N_ZX)$ by the
 isomorphisms (see Definition \ref{df:abs_pur}):
\begin{equation} \label{eq:deformation_iso}
\E^{2n,n}(\Th(N_ZX))
 \xrightarrow{(\sigma_0^*)^{-1}} \E^{2n,n}_{\AA^1_Z}(D_ZX)
 \xrightarrow{\ \sigma_1^*\ } \E^{2n,n}_Z(X).
\end{equation}
We denote it by $\rfdl_X(Z)$.

We also define the \emph{fundamental class} of $Z$ in $X$
 as the following class of $\E^{2n,n}(X)$:
$$
\fdl_X(Z):=i_!\left(\rfdl_X(Z)\right).
$$
\end{df}
Because $\sigma_0^*$, $\sigma_1^*$ are isomorphisms of $\E^{**}(Z)$-modules,
 and according to the Thom isomorphism \eqref{eq:thom_iso},
 the class $\rfdl_X(Z)$ is in fact a base of the $\E^{**}(Z)$-bigraded module
 $\E^{**}_Z(X)$. In other words, we get an isomorphism
\begin{equation} \label{eq:purity}
\pur{X,Z}:\E^{**}(Z) \rightarrow \E^{**}_Z(X), z \mapsto z.\rfdl_X(Z)
\end{equation}
called the \emph{purity isomorphism}.

\begin{rem} \label{rem:non_constant_codim}
For the reader who does not want to restrict to pure codimension
 pairs, let us indicate that in general, the codimension
 of a regular closed pair $(X,Z)$ is Zariski local
 in $Z$, as well as the rank of its normal bundle $N_ZX$.
Then the convention of \ref{rem:non_constant_rk}
 will apply both to the refined fundamental class
  and the purity isomorphism associated with $(X,Z)$.
\end{rem}

\begin{ex}
Let $E/X$ be a vector bundle.
Then the closed $X$-pair $(E,X)$ corresponding to the zero section
 is $\E$-pure -- see Example \ref{ex:abs_pure}(1).
 Then we get from the above definition the
 equality of classes in $\E^{**}(\Th(E/X))$:
\begin{equation} \label{rfdl_class&vector_bdl}
\rfdl_E(X)=\rthom(E/X).
\end{equation}
Indeed, one can identify the deformation space
 $D_XE$ (resp. normal bundle $N_XE$)
 with the affine line $\AA^1_E$ (resp. the $X$-vector bundle $E$ itself)
 and the  maps $\sigma_1, \sigma_0$ of \eqref{eq:deformation}
 with respectively the unit and zero section of $\AA^1_E$
 (see \cite[Rem. 4.2]{Deg8}) so that the above identification
 follows from the previous definition.
Moreover, if $\bar E/X$ denotes the projective completion of $E/X$,
 so that $E \subset \bar E$ is an open subscheme,
 one gets by Nisnevich excision an  identification
 $\E^{**}_X(\bar E)=\E^{**}_X(E)$ through which
 $\rfdl_{\bar E}(X)$ corresponds to the class $\rfdl_E(X)=\rthom(E/X)$.
 In particular, one gets the equality of (non refined) fundamental classes:
\begin{equation} \label{fdl_class&proj_bdl}
\fdl_{\bar E}(X)=\thom(E/X).
\end{equation}
\end{ex}

\begin{rem} \label{rem:purity&MGL-mod}
Consider a smooth closed $S$-pair $(X,Z)$ of codimension $n$.
Recall from Example \ref{ex:abs_pure}(1) that $(X,Z)$ is $\E$-pure.
We have defined in \cite[4.6]{Deg8} a purity isomorphism
 in the homotopy category of $\MGL_S$-modules:
$$
\MGL_S(X/X-Z) \xrightarrow \sim \MGL_S(Z)(n)[2n].
$$
According to the above example and \cite[4.3, 4.4]{Deg8},
 we get that the induced morphism in $\E$-cohomology
 (see \ref{num:apply_Deg8}) coincides exactly with
 the purity isomorphism $\pur{X,Z}$ introduced above.
In particular, the class $\rfdl_X(Z)$ in $\E^{**}_Z(X)$
 introduced here coincides with that obtained from \cite[4.14]{Deg8}
 by considering $\MGL_S$-modules.\footnote{The terminology
 here differs slightly from that of \emph{loc. cit.}:
 we use the adjective "refined" instead of "localized" which
 seems more classical.}
\end{rem}

\subsection{Intersection theory}
\label{sec:inter_theory}

\begin{df}
Consider a morphism $\Delta:(Y,T) \rightarrow (X,Z)$ of closed pairs
 (Par. \ref{num:closed_pairs})
 such that $(X,Z)$ and $(Y,T)$ are $\E$-pure.

Then, according to Definition \ref{df:fdl_class},
 there exists a unique class $e_\Delta$ in
 $\E^{**}(T)$ such that
$$
\Delta^*(\rfdl_X(Z))=e_\Delta.\rfdl_Y(T).
$$
We call $e_\Delta$ the \emph{defect} of the morphism $\Delta$.\footnote{
One could call this class the \emph{defect of transversality}.
In fact, there are two kinds of possible defect: excess of intersection,
 ramification.}
\end{df}

The following result generalizes \cite[Prop. 4.16]{Deg8}:
\begin{thm} \label{thm:pullback&fdl_class}
With the notations of the above definition,
 assume $\Delta=(f,g)$ is cartesian.
Then it induces a monomorphism $\nu:N_TY \rightarrow g^{-1}(N_ZX)$
 of vector bundles over $T$. We put 
\begin{equation}\label{eq:excess_bdl}
\xi=g^{-1}(N_ZX)/N_TY
\end{equation}
and denote by $e$ the rank of this vector bundle.

Then $e_\Delta=c_e(\xi)$.
\end{thm}
\begin{rem}
One can apply this theorem in two cases:
\begin{itemize}
\item Given a base scheme $S$, $\Delta$ is a morphism of smooth closed $S$-pairs.
This case is already known from \cite[4.16]{Deg8}.
\item $\E$ is $\reg$-absolutely pure and the schemes $X$, $Z$, $Y$, $T=f^{-1}(Z)$
 are all regular.
\end{itemize}
\end{rem}
\begin{proof}
Recall the deformation diagram \eqref{eq:deformation} is functorial
 with respect to the morphism $\Delta$, because it is assumed to be cartesian.
 Thus, going back to the definition of fundamental classes
 (\ref{df:fdl_class}), using the deformation diagrams respectively
 for the closed pairs $(X,Z)$ and $(Y,T)$
 we are reduced to compare Thom classes for the corresponding normal bundles,
 as summarized in the following cartesian square:
$$
\xymatrix@C=40pt@R=14pt{
T\ar_g[dd]\ar@{^(->}[r] & N_TY\ar^\nu[d] \\
& g^{-1}(N_ZX)\ar^{g'}[d] \\
Z\ar@{^(->}[r] & N_ZX,
}
$$
where $\nu$ is a monomorphism of vector bundles over $T$.
 This now follows from \cite[Lem. 4.18]{Deg8}
 (using the considerations of Paragraph \ref{num:apply_Deg8}).
\end{proof}

In terms of fundamental classes, we get:
\begin{cor} \label{cor:pullback&fdl_class}
Let $(f,g):(Y,T) \rightarrow (X,Z)$ be a cartesian morphism of $\E$-pure closed pairs.
 Let $f^*:\E^{**}_Z(X) \rightarrow \E^{**}_T(Y)$
 be the morphism defined in Paragraph \ref{num:pullback_support}.

Then, using Definition \eqref{eq:excess_bdl}, we get the following formula
 in $\E^{**}_T(Y)$:
$$
f^*(\rfdl_X(Z))=c_e(\xi).\rfdl_Y(T).
$$
In particular when $f$ is transversal to $Z$,
\begin{equation} \label{eq:transversal_pullback}
f^*(\rfdl_X(Z))=\rfdl_Y(T).
\end{equation}
\end{cor}

\begin{rem}
According to formulas (E3) and (E7) of Proposition \ref{prop:ppty_product},
 we get in the assumptions of the above corollary the even more usual formula
 in $\E^{**}(Y)$:
$$
f^*(\fdl_X(Z))=c_e(\xi).\fdl_Y(T).
$$
Applying this formula in the case where $f$ is the zero section $s:X \rightarrow E$
 of a vector bundle of rank $n$, we get the relation
$$
s^*(\fdl_E(X))=c_n(E).\footnote{In other words,
 given a vector bundle $E/X$
 the fundamental class of its zero section
 coincides up to homotopy with its \emph{Euler class}.
 This fact justifies our choice of Thom class in Definition \ref{df:Thom}.}
$$
These two relations give the following classical trick,
 which will be used later, to compute fundamental classes.
\end{rem}
\begin{cor} \label{cor:fdl&Chern_transversal_section}
Let $(X,Z)$ be an $\E$-pure closed pair of codimension $n$
 corresponding to an immersion $i:Z \rightarrow X$.
Assume there exists a vector bundle $E/X$ with a section $s$ transversal
 to the zero section $s_0$
 and such that $s_0^{-1}(s)=i$.
 Then
$$
\fdl_X(Z)=c_n(E).
$$
\end{cor}

An important point of intersection theory
 is the associativity of the intersection product.
 In our setting, this can be expressed nicely using the refined product
 of Paragraph \ref{num:products}. We begin with a particular case of Theorem \ref{thm:product_fdl}
 which will be the crucial step.
\begin{prop}\label{prop:product&thom}
Consider an exact sequence $(\sigma)$ of vector bundles over a	scheme $X$:
$$
0 \rightarrow E' \xrightarrow{\nu} E \rightarrow E'' \rightarrow 0.
$$
Then the following relation holds in $\E^{**}(\Th(E))$:
$$
\rthom(E/X)=\rthom(E'/X).\rthom(E/E')
$$
using the pairing $\E^{**}_X(E') \otimes \E^{**}_{E'}(E) \rightarrow \E^{**}_X(E)$
 (Paragraph \ref{num:products}).
\end{prop}
\begin{proof}
Note that according to formula \eqref{rfdl_class&vector_bdl},
 we have to prove the equality:
\begin{equation} \label{eq:product&thom:pf1}
\rfdl_E(X)=\rfdl_{E'}(X).\rfdl_E(E').
\end{equation}

As in the proof of \cite[4.1.1]{Riou}, we can find a torsor $T$
 over the $X$-vector bundle $\uHom(E'',E')$ such that
 the sequence $(\sigma)$ splits over $T$.
 The morphism $T \rightarrow X$ is an $\AA^1$-weak equivalence ; thus,
 by compatibility of the Thom class with base change, 
 we can assume the sequence $(\sigma)$ splits.

Let us consider $P'$ (resp. $P''$) the projective completion
 of $E'/X$ (resp. $E''/X$), $P=P' \times_X P''$,
 and the following commutative diagram 
 made of cartesian squares:
$$
\xymatrix@R=15pt@C=24pt{
X\ar[r]\ar@{=}[d] & E'\ar[r]\ar^{j'}[d] & E\ar^j[d] \\
X\ar^{s'}[r]\ar@/_5pt/_s[rr] & P'\ar^{s''}[r] & P
}
$$
where $j$ (resp. $s$, $s'$, $s''$) stands for the natural open immersion
 (resp. canonical sections coming from the obvious zero sections).
Remark that if we apply the functor $j^*:\E^{**}_X(P) \rightarrow \E^{**}_X(E)$
 to the following equality:
\begin{equation} \label{eq:product&thom:pf2}
\rfdl_{P}(X)=\rfdl_{P'}(X).\rfdl_{P}(P'),
\end{equation}
we get \eqref{eq:product&thom:pf1}
 -- use formula (E5) of Prop. \ref{prop:ppty_product}
  and formula \eqref{eq:transversal_pullback}.
Thus we are reduced to prove that equality \eqref{eq:product&thom:pf2}
 holds in $\E^{**}_X(P)$.

According to the projective bundle formula,
 the morphism $s_!:\E^{**}_X(P) \rightarrow \E^{**}(P)$ is
 a split monomorphism (see the end of Par. \ref{num:Thom_space} for more details).
 Thus it is sufficient to prove that \eqref{eq:product&thom:pf2} holds
 after applying $s_!$. That means we have to prove $\fdl_P(X)$
 is equal to:
\begin{align*}
s_!\big(\rfdl_{P'}(X).\rfdl_{P}(P')\big)&=
s''_!s'_!\big(\rfdl_{P'}(X).\rfdl_{P}(P')\big) \\
& \stackrel{(E6)}=s''_!\big(\fdl_{P'}(X).\rfdl_{P}(P')\big)
=s''_!\left(s^{\prime\prime*}[\fdl_{P}(P'')].\rfdl_{P}(P')\right) \\
& \stackrel{(E7)}=\fdl_{P}(P'').\fdl_{P}(P')
\end{align*}
using once again the properties of the refined product enumerated
 in Proposition \ref{prop:ppty_product}.
 This last check follows from a direct computation
 -- see \cite[Lem. 4.25]{Deg8}.
\end{proof}

\begin{rem} \label{rem:virtual_Thom_classes}
Consider a scheme $X$. Given vector bundles $E'$ and $E''$ over $X$,
 we get using the product defined in Paragraph \ref{num:products} a canonical
 map of K\"unneth type:
$$
\E^{**}(\Th(E')) \otimes_{\E^{**}(X)} \E^{**}(\Th(E'')
 \rightarrow \E^{**}(\Th(E' \oplus E'')).
$$
According to the above proposition, $\rthom(E' \oplus E'')=\rthom(E').\rthom(E'')$.
In other words, the above map is an isomorphism (of bigraded $\E^{**}(X)$-modules).

Now, given an exact sequence $(\sigma)$ of vector bundles as in the above proposition,
 Riou in the proof of \cite[4.1.1]{Riou} shows
 there exists a canonical isomorphism in $\SH(X)$
$$
\Sus \Th(E') \wedge \Sus \Th(E'')\xrightarrow{\epsilon_\sigma} \Sus \Th(E).
$$
Using this isomorphism, one defines a canonical product of bigraded $\E^{**}(X)$-modules:
$$
\xymatrix@C=22pt@R=0pt{
{\E}^{**}(\Th(E')) \otimes_{{\E}^{**}(X)} {\E}^{**}(\Th(E''))\ar[r]
 & {\E}^{**}(\Th(E') \wedge \Th'(E''))\ar^-{(\epsilon_\sigma^*)^{-1}}[r]
 & {\E}^{**}(\Th(E)). \\
\qquad\qquad\qquad\qquad\qquad\qquad a \otimes b\ar@{|->}[rr] && \langle a,b\rangle_\sigma \qquad
}
$$
Then the relation obtained in the previous proposition can be rewritten
 as follows:
$$
\rthom(E)=\langle\rthom(E'),\rthom(E'')\rangle_\sigma.
$$
Thus, the map $\langle -,-\rangle_\sigma$ is an isomorphism.

Recall that Deligne defines in \cite[4.12]{Del}
 the category of \emph{virtual vector bundle}
 denoted by $\underline K(X)$. As remarked by Riou,
  the isomorphisms of type $\epsilon_\sigma$ show
  that
 the functor $\Sus \Th_X$ induces a canonical functor
$$
\Th_X:\underline K(X) \rightarrow \SH(X).
$$
In fact the preceding considerations show that,
 for any virtual vector bundle $\xi$ in $\underline K(X)$,
 the bigraded $\E^{**}(X)$-module $\E^{**}(\Th(\xi))$
 is free of rank $1$. Moreover, it admits a canonical trivialization
 which we will denote $\rthom(\xi)$: if $\xi=[E]-[E']$ for two
 vector bundles $E$ and $E'$, one puts:
$$
\rthom(\xi):=\rthom(E).\rthom(E')^*.
$$
A conceptual way of stating this result:
 the canonical functor $\E^{**} \circ \tilde \Th_X$
 induces a functor of Picard categories:
$$
\big(\underline K(X)\big)^{op} \rightarrow \Pic\big(\E^{**}(X)\big)
$$
where the right hand side category is the Picard category
 of bigraded virtual line bundles over $\E^{**}(X)$.
\end{rem}

\begin{thm} \label{thm:product_fdl}
Consider regular closed immersions $Z \xrightarrow k Y \xrightarrow i X$
 in $\base$
 and assume $\E$ is absolutely pure.

Then we get the following equality in $\E^{**}_Z(X)$:
$$
\rfdl_X(Z)=\rfdl_Z(Y).\rfdl_X(Y),
$$
using the pairing $\E^{**}_Z(Y) \otimes \E^{**}_{Y}(X) \rightarrow \E^{**}_Z(X)$
 (Paragraph \ref{num:products}) for the right hand side.
\end{thm}
\begin{proof}
Recall from Paragraph \ref{num:def_space} the deformation space $D_ZX$ 
 associated with the closed pair $(X,Z)$.
 Following \cite[\S 10]{Ros},
 we define the double deformation space by the following formula:
$$
D(X,Y,Z)=D(D_Z(X),D_Z(X)|_Y).
$$
This scheme is flat over $\AA^2$.
 Its fiber over $(1,1)$ is $X$ while its fiber over $(0,0)$
 is the normal bundle $E:=N\big(N_ZX,N_ZY\big)$.
Let us put $D'=D(Y,Y,Z)$ and $E'=N_ZY$. Then we get the following diagram
made of cartesian squares:
$$
\xymatrix@R=8pt@C=20pt{
Z\ar^{s'}[r]\ar_{s_1}[d] & E'\ar^\nu[r]\ar[d] & E\ar^{d_0}[d] \\
\AA^2_Z\ar[r] & D'\ar[r] & D \\
Z\ar^k[r]\ar^{s_0}[u] & Y\ar^i[r]\ar[u] & X\ar_{d_1}[u]
}
$$
where the first (resp. third) row is seen as the $(0,0)$-fiber
 (resp. $(1,1)$-fiber) of the second row, with respect to its canonical
 projection to $\AA^2$. Because this projection is flat,
 all the squares in this diagram are transversal.

Thus, one can apply Corollary \ref{cor:pullback&fdl_class}
 to the morphism of closed pairs $(X,Z) \xrightarrow{d_1} (D,\AA^2_Z)$
 (resp. $(E,Z) \xrightarrow{d_0} (D,\AA^2_Z)$):
 because $s_1$ (resp. $s_0$) is a strong $\AA^1$-homotopy equivalence,
 and using the absolute purity property,
 the pullback morphism $d_1^*:\E^{**}_{\AA^2_Z}(D) \rightarrow \E^{**}_Z(X)$
 (resp. $d_0^*:\E^{**}_{\AA^2_Z}(D) \rightarrow \E^{**}_Z(E)$)
 is an isomorphism.

Applying again Corollary \ref{cor:pullback&fdl_class},
 we deduce that to prove the theorem, it is enough to prove the relation
$$
\rfdl_E(Z)=\rfdl_{E'}(Z).\rfdl_E(E').
$$
In view of formula \eqref{rfdl_class&vector_bdl}, this is precisely 
 Proposition \ref{prop:product&thom}.
\end{proof}

\begin{rem}
In the case where $\base$ is the category of smooth $S$-schemes,
 for a base scheme $S$,
 the above theorem is very close to \cite[4.30]{Deg8}.
 The proof given here is considerably simpler,
 as we use the refined product of \ref{num:products} -- it uses
 the localization property of $\SH$, a very strong result.
\end{rem}
\section{Gysin morphisms and localization long exact sequence}

\label{sec:Gysin}

In this section, we fix an absolute oriented ring spectrum $(\E,c)$
 -- recall Def. \ref{df:orientation}.

\subsection{Residues and the case of closed immersions}

\begin{num} \label{num:loc_exact_seq}
Consider a closed immersion $i:Z \rightarrow X$ 
 with complementary open immersion $j:U \rightarrow X$.
 Property (A4) of Paragraph \ref{num:functoriality_properties_SH}
  implies the existence of the 
 \emph{localization long exact sequence}:
\begin{equation} \label{eq:loc_exact_seq}
\cdots \E^{n,m}_Z(X) \xrightarrow{i_!} \E^{n,m}(X) \xrightarrow{j^*} \E^{n,m}(U)
 \xrightarrow{\delta_{X,Z}} \E^{n+1,m}_Z(X) \cdots
\end{equation}
If we assume that the corresponding closed pair $(X,Z)$
 is regular of codimension $c$ and $\E$-pure, 
 the associated purity isomorphism
 \eqref{eq:purity} induces a long exact sequence in cohomology
 from the preceding one:
\begin{equation} \label{eq:Gysin_exact_seq}
\cdots \E^{n-2c,m-c}(Z) \xrightarrow{i_*} \E^{n,m}(X)
 \xrightarrow{j^*} \E^{n,m}(U)
 \xrightarrow{\partial_{X,Z}} \E^{n-2c+1,m-c}(Z) \cdots
\end{equation}
\end{num}
\begin{df} \label{df:Gysin_triangle}
Under the notations above, assuming $(X,Z)$ is $\E$-pure of codimension $c$,
 we call \eqref{eq:Gysin_exact_seq} (resp. $i_*$, $\partial_{X,Z}$)
 the Gysin long exact sequence (resp. Gysin morphism, residue morphism)
 associated with the closed pair $(X,Z)$ (or with the immersion $i$).
\end{df}
Note that in terms of the refined fundamental class (Def. \ref{df:fdl_class}),
 the Gysin morphism can be described by the following formula
 for a cohomology class $z \in \E^{**}(Z)$:
\begin{equation} \label{eq:Gysin&rfdl}
i_*(z):=i_!\big(z.\rfdl_X(Z)\big).
\end{equation}
Moreover, the residue of a cohomology class $u \in \E^{**}(U)$ is
 uniquely determined by the following property:
\begin{equation} \label{eq:residue&rfdl}
\delta_{X,Z}(u)=\partial_{X,Z}(u).\rfdl_X(Z).
\end{equation}

\begin{rem} \label{rem:Gysin_depends_orient}
Note that the Gysin morphism $i_*$ defined above do depend
 on the orientation $c$ on the absolute ring spectrum $\E$,
 fixed at the beginning of the section --- this dependence
 is responsible for the Riemann-Roch theorem \ref{thm:RR1}.
 Usually, there will be no possible confusion about the chosen
 orientation, but in case there might, we will denote by $i_*^c$
 the corresponding Gysin morphism (see in particular
 \ref{num:orient_in_notations}).
\end{rem}

\begin{prop} \label{prop:basic_Gysin_closed}
\begin{enumerate}
\item Assume $\E$ is absolutely pure.
Then for any $\E$-pure closed immersions
$$
T \xrightarrow k Z \xrightarrow i X,
$$
 the following equality holds: $i_*k_*=(ik)_*$.
\item Consider an $\E$-pure closed immersion $i:Z \rightarrow X$.
 For any couple $(x,z)$ in $\E^{**}(X) \times \E^{**}(Z)$,
$$
i_*\big(i^*(x).z\big)=x.i_*(z).
$$
\item Consider a cartesian square 
$$
\xymatrix@=10pt{
T\ar^k[r]\ar_g[d] & Y\ar^f[d] \\
Z\ar_i[r] & X
}
$$
in $\base$ such that $i$ and $k$ are $\E$-pure closed immersions.
 Let $h:(Y-T) \rightarrow (X-Z)$ be the morphism induced by $f$.
Let $\xi$ be the vector bundle over $T$ defined by formula \eqref{eq:excess_bdl}
 with respect to the morphism of closed pairs $(f,g)$.
 Let $e$ be the rank of $\xi$.

Then the following formulas hold:
\begin{align}
f^*i_*(z)=k_*\big(c_e(\xi).g^*(z)\big). \\
\partial_{Y,T} h^*(u)=c_e(\xi).g^*\partial_{X,Z}(u).
\end{align}
\end{enumerate}
\end{prop}
\begin{proof}
Taking care of formulas \eqref{eq:Gysin&rfdl}
 and \eqref{eq:residue&rfdl},
 points (1), (2), (3) are respectively consequences 
 of Th. \ref{thm:product_fdl}, Proposition \ref{prop:ppty_product}(E4)+(E7),
 Corollary \ref{cor:pullback&fdl_class}.
\end{proof}

\begin{rem} \label{rem:Gysin_closed_im&Deg8}
Let $i:Z \rightarrow X$ be a closed immersion between smooth $S$-schemes.
According to  point (1) of Example \ref{ex:abs_pure}, the closed pair $(X,Z)$
 is $\E$-pure.
Then both Example \ref{ex:Gysin&MGL-mod}
 and the preceding definition gives a pushforward
 of the form
$$
i_*:\E^{**}(Z) \rightarrow \E^{**}(X).
$$
Remark \ref{rem:purity&MGL-mod} shows that these two constructions
 coincide.
\end{rem}

As we will see in the next examples, residues are connected to the
 \emph{tame residue symbols} of Milnor in K-theory.
 This connection can be reduced to the following interesting
 property:
\begin{prop}\label{prop:specialisation}
Let $X$ be a regular scheme and $Z$ be a principal divisor
 of $X$ parametrized by a regular function $\pi:X \rightarrow \AA^1_k$.
Let us denote by $i:Z \rightarrow X$ ($j:U \rightarrow X$)
 be the corresponding closed immersion (resp. complementary open immersion).
 Let us consider the following composite map:
$$
\gamma_\pi:\E^{n,m}(U) \xrightarrow{(1)} 
 \E^{n,m}(U) \oplus \E^{n+1,m+1}(U) \simeq \E^{n+1,m+1}(\GG \times U)
 \xrightarrow{(2)} \E^{n+1,m+1}(U)
$$
where $(1)$ is the obvious inclusion
 and $(2)$ is the pullback by the graph of 
 the induced map $\pi|_U:U \rightarrow \GG$.

Then if $(X,Z)$ is $\E$-pure, the following relation holds:
$$
\partial_{X,Z} \circ \gamma_\pi \circ j^*=i^*.
$$
\end{prop}
\begin{proof}
The proof is essentially the same as the one explained in
 \cite[Prop. 2.6.5]{Deg5}. Let us indicate the main steps.

According to the definition of the purity isomorphism
 through the deformation diagram \eqref{eq:deformation},
 and the fact the canonical morphism $\sigma_1:X \rightarrow D_ZX$
 is a split monomorphism (cf \emph{loc. cit.}), we restrict to the
 case of $(N_ZX,Z)$ where $N_ZX$ is the normal bundle of $Z$ in $X$,
 while $\pi$ can be considered as a global trivialization of
 this line bundle.

Up to isomorphism, we thus are reduced to the case of $(\AA^1_Z,Z)$
 and the canonical parametrization $t$ of the affine line (given by the identity).
 Then, the problem can be unfolded to a trivial identity given at the level
 of schemes (cf. the end of the proof of \emph{loc. cit.}).
\end{proof}

\begin{ex} \label{ex:tame_residue}
Let $R$ be a discrete valuation ring with valuation $v$,
 fraction field $K$ and residue field $k$.
 
Recall first that the \emph{tame residue symbol}
 on Milnor K-theory:
\begin{equation} \label{eq:tame_residue}
\partial_v:K_n^M(K) \rightarrow K_{n-1}^M(k)
\end{equation}
is uniquely characterized by the following properties,
where $\pi$ is a given uniformizing parameter:
\begin{enumerate}
\item $\partial_v( \{\pi\})=1$ ;
\item for units $u_1,...,u_r \in R^\times$,
 $\partial_v(\{\pi,u_1,\hdots,u_r\})=\{\bar u_1,\hdots,\bar u_r\}$;
\item for units $u_1,...,u_r \in R^\times$,
 $\partial_v(\{u_1,\hdots,u_r\})=0$.
\end{enumerate}

If we assume that $(\spec R,\spec k)$ is $\E$-pure, 
 we get an abstract residue map with coefficients in $\E$:
$$
\partial_v^\E:\E^{n,m}(K) \rightarrow \E^{n-1,m-1}(k).
$$
When $R$ is of equal characteristics,
 we can apply this construction to the case of Voevodsky's Einlenberg-MacLane motivic
 ring spectrum $\HH_\ZZ$ -- $S$ is the spectrum of the prime field of $R$.
 Then, the induced map
$$
\partial_v:K_n^M(K) \simeq \HH_{\mathcal M}^{n,n}(K)
 \rightarrow \HH_{\mathcal M}^{n-1,n-1}(k)\simeq K_{n-1}^M(k)
$$
coincides with \eqref{eq:tame_residue}. In fact,
 relations (1) and (2) follows from the previous proposition
 while relation (3) follows from the fact $\partial_v \circ j^*=0$
 according to the Gysin long exact sequence \eqref{eq:Gysin_exact_seq}.

In the unequal characteristics case,
 we get the same result with rational coefficients by considering
 Beilinson motivic cohomology ring spectrum $\HB$. 
\end{ex}

\begin{rem}
Other applications of the previous proposition will be given
 in Section \ref{sec:application_residual}.
\end{rem}

\subsection{Projective lci morphisms}

\begin{num} \label{num:Gysin4projection}
Let $X$ be a scheme, put $P=\PP^n_X$ and consider the canonical projection
 $p:P \rightarrow X$.
 Recall we have introduced in Example \ref{ex:Gysin&MGL-mod}
 the Gysin morphism: $p_*:\E^{*,*}(\PP^n_X) \rightarrow \E^{*-2n,*-n}(X)$
 based on the construction of \cite{Deg8}.
 We give an alternative construction which uses the point of view
 considered in the present setting.
 It is based on the following facts which follow directly from the projective
 bundle theorem (\ref{thm:PBF}). Consider the bigraded ring $A=\E^{**}(X)$.
 In the paragraph and definition that follow, we work in the category 
 of bigraded $A$-modules and refer to them simply as $A$-modules:
\begin{itemize}
\item The bigraded group $\E^{**}(P)$ is a free $A$-module of finite rank.
\item The K\"unneth map
$$
\E^{**}(P) \otimes_A \E^{**}(P) \rightarrow \E^{**}(P \times_X P),
 (x,y) \mapsto p_1^*(x)\cupp p_2^*(y)
$$
is an isomorphism of $A$-modules.
\end{itemize}
From the first point, we deduce that the $A$-module $\E^{**}(P)$ is dualizable.
Let $\E^{**}(P)^\vee$ be its $A$-dual and 
$$
ev:\E^{**}(P)^\vee \otimes_A \E^{**}(P) \rightarrow A
$$
be the evaluation map.
Let $\delta:P \rightarrow P \times_X P$ be the diagonal immersion of $P/X$.
Using the Gysin morphism and the second point, we get a co-pairing of $A$-modules:
\begin{equation} \label{eq:copairing_duality_Pn}
\epsilon_P:A=\E^{**}(X) \xrightarrow{p^*} \E^{**}(P)
 \xrightarrow{\delta_*} \E^{**}(P \times_X P)
 \simeq \E^{**}(P) \otimes_A \E^{**}(P).
\end{equation}
We claim this is a duality co-pairing in the sense that the induced map
\begin{equation} \label{eq:duality_Pn}
\mathcal D_P:\E^{**}(P)^{\vee}
 \xrightarrow{1 \otimes \epsilon_P}
  \E^{**}(P)^{\vee} \otimes_A \E^{**}(P) \otimes_A \E^{**}(P)
 \xrightarrow{ev \otimes 1} \E^{**}(P)
\end{equation}
is an isomorphism. Indeed, the computation of the matrix of $\mathcal D_P$
 in the base given by the projective bundle theorem is precisely the same as 
 that of \cite[Lem. 5.5]{Deg8}: it is a triangular lower matrix
 with the identity  diagonal. 
 Note that $\mathcal D_P$ is homogeneous of degree $(2n,n)$.
\end{num}
\begin{df} \label{df:Gysin_Pn}
Using the above notations, one defines the Gysin morphism associated with $p$
as the following composite morphism:
$$
p_*:\E^{**}(P) \xrightarrow{\mathcal D_P^{-1}} \E^{**}(P)^\vee
 \xrightarrow{(p^*)^\vee} \E^{**}(X)^\vee=\E^{**}(X),
$$
homogeneous of degree $(-2n,-n)$.
\end{df}
In other words, $p_*$ is the transpose of $p^*$ with respect
 to the duality given by the co-pairing \eqref{eq:copairing_duality_Pn}.
 In view of \cite[Prop. 5.26]{Deg8}, the above definition coincides
 with that of Example \ref{ex:Gysin&MGL-mod}.

\begin{rem}\label{rem:compute_Gysin_proj_bdl}
Here is a practical way to determine the Gysin map $p_*$ constructed in the
 preceding definition. Let us again denote by $A$ the cohomology ring
 $\E^{**}(X)$.

 First, one determines the map $\epsilon_P$ --- see \eqref{eq:copairing_duality_Pn}.
 According to the projective bundle theorem, the $\E$-cohomology of
 $P \times_X P$ is a free $A$-module with basis $c^i.d^j$, $0 \leq i,j \leq n$,
 where $c$ (resp. $d$) is the first Chern class of the canonical
 line bundle on the first (resp. second) factor of $P \times_X P$ (see
 Theorem \ref{thm:PBF}).
 The $A$-linear map $\epsilon_P$ is uniquely determined by its value
 on the multiplicative identity element $1$ of $A$,
 which can be a priori written as:
$$
\epsilon_P(1)=\delta_*(1)
 =\sum_{0 \leq i,j \leq n} \pi^*\big(\eta_{ij}^{(n)}\big).c^id^j.
$$
where $\pi$ is the projection map of $P \times_X P/X$.

Second, the map $\mathcal D_P$ --- see \eqref{eq:duality_Pn} ---
 is an $A$-linear map between rank $n+1$ free $A$-modules.
 Then the matrix $M$ of $\mathcal D_P$ with respect to the basis 
 on the source and targets coming from the projective bundle theorem applied to
 $P/X$ can be easily computed with the previous notation as:
$$
M:=\left(\eta_{ij}^{(n)}\right)_{\{0 \leq i,j \leq n\}}.
$$
From what was said in the paragraph preceding the above definition,
 the matrix $M$ is invertible
 --- it is in fact symmetric, lower triangular with $1$ on the anti-diagonal;
 see Example \ref{ex:Myschenko} for a computation of this matrix in terms
 of the formal group law associated with $(\E,c)$.

According to the preceding definition, the Gysin map $p_*$ is then
 completely determined by the inverse matrix $M^{-1}$.
 In fact, here is a way to compute the image of an element $x \in \E^{**}(P)$
 under the Gysin map $p_*$:
\begin{itemize}
\item We write:
$$
x=\sum_{0 \leq i \leq n} p^*(x_i).c^i
$$
where $c$ is the canonical line bundle on $P$.
 Let $X$ be the element of $A^{\{0,...,n\}}$ whose value at the $i$-th place
 is $x_i$.
\item The map $(p^*)^\vee$ just extract the coefficient of the element
 $(c^0)^\vee$ with respect to the decomposition of elements of 
 $\E^{**}(P)$ in the dual base $\big((c^i)^\vee\big)$. Therefore, 
 we get that $p_*(x)$ is the $0$-th coefficient of $M^{-1}.X$, 
 seen as an element of $A^n$.
\end{itemize}
Note that, given the cofactors formula to compute the inverse of a matrix,
 we get in particular the following formula:
\begin{equation} \label{eq:pre_Myschenko}
p_*(1)=\frac 1 {\det(M)}.\det\left(\eta_{ij}^{(n)}\right)_{\{1 \leq i,j \leq n\}}.
\end{equation}
This will be made precise in Example \ref{ex:Myschenko}
\end{rem}

\begin{num}\label{num:basic_Gysin_projection}
The Gysin morphism introduced above satisfies the following properties:
\begin{enumerate}
\item For any integers $n,m \geq 0$, considering the canonical projections
$$
\xymatrix@=10pt{
\PP^n_X \times_X \PP^m_X\ar_{p'}[d]\ar^-{q'}[r]
 & \PP^n_X\ar^p[d] \\
\PP^m_X\ar^q[r] & X
}
$$
one has:
$p_*q'_*=q_*p'_*$.
\item For any cartesian square
$$
\xymatrix@=10pt{
\PP^n_Z\ar_{q}[d]\ar^-{l}[r]
 & \PP^n_X\ar^p[d] \\
Z\ar^i[r] & X
}
$$
where $i$ is an $\E$-pure closed immersion and $p$ is the canonical
 projection: $i_*q_*=p_*l_*$.
\item For any integer $n \geq 0$ and any section $s$ of the
 projection $p:\PP^n_X \rightarrow X$:
  $p_*s_*=1$.
\end{enumerate}
For the proof of these properties,
 we use Remark \ref{rem:Gysin_closed_im&Deg8} and
 refer the reader to \cite{Deg8}, respectively
 Lemmas 5.8, 5.9 and 5.10.
Then the following lemma follows formally from these properties
 -- see \cite[Lem. 5.11]{Deg8}:
\end{num}

\begin{lm} \label{lm:Gysin}
Assume $\E$ is absolutely pure and consider a commutative diagram~:
$$
\xymatrix@C=14pt@R=-4pt{
& {}\PP^n_X\ar^p[rd] & \\
Y\ar^k[ru]\ar_i[rd] & & X \\
& {}\PP^m_X\ar_q[ru] &
}
$$
where $i$ (resp. $k$) is a closed immersion and $p$ (resp. $q$) is the canonical projection.
Then, using the above definitions, $p_*k_*=q_*i_*$.
\end{lm}

\begin{num} \label{num:lci_morphisms}
Recall that an $X$-scheme $Y$ is said to be a \emph{local complete intersection}
 if it admits locally an immersion into an affine space $\AA^n_X$
 -- see \cite[VIII, 1.1]{SGA6}. We will say abusively \emph{lci} instead
 of local complete intersection. \\
 Given a projective lci morphism $f:Y \rightarrow X$ in $\base$,
 it admits a factorization
 of the form 

\begin{equation} \label{eq:generic_fact_lci}
Y \xrightarrow i \PP^n_X \xrightarrow p X
\end{equation}
 where $i$ is a regular closed immersion and $p$ is
 the canonical projection.
 This follows from our convention on projective morphisms
  and \cite[1.2]{SGA6}.

Now, assume that $\E$ is absolutely pure.
From the preceding lemma, the composite morphism:
$$
\E^{**}(Y) \xrightarrow{i_*} \E^{**}(\PP^n_X)
 \xrightarrow{p_*} \E^{**}(X)
$$
obtained using Definition \ref{df:Gysin_triangle}
 and Paragraph \ref{num:Gysin4projection} is independent
 of the choice of the factorization of $f$.
\end{num}
\begin{df}\label{df:Gysin}
Consider the above assumptions and notations.

We put $f_*=p_*i_*$ and call it the Gysin morphism
 associated with $f$.
\end{df}
Note that if $f$ has dimension $d$, then $f_*$
 is homogeneous of bidegree $(-2d,-d)$.

\begin{rem}
The warning of Remark \ref{rem:Gysin_depends_orient} can be equally
 applied to the preceding Gysin morphism. In fact, it is proved 
 in Theorem \ref{thm:compare_Gysin} that the preceding Gysin morphisms
 are uniquely determined by the chosen orientation $c$ on $\E$.
\end{rem}

\begin{ex} \label{ex:gysin_H&KGL}
\begin{enumerate}
\item \textit{Motivic cohomology and Chow groups}.-- Let $k$ be a perfect field and consider the $k$-absolute ring spectrum
 $\HH_\ZZ$ of Example \ref{ex:abs_ring_spectra_orientation}.
 It is well known that for any integer $n \geq 0$
 and any smooth $k$-scheme:
$$
\HH^{2n,n}_\cM(X,\ZZ) \simeq CH^n(X).
$$
In fact, we know that the Nisnevich sheaf $H^n(\ZZ(n))$ is the
 unramified Milnor K-theory $\mathcal K_n^M$,
 and from the hypercohomology spectral sequence 
 associated with the homotopy $t$-structure along with some known
 vanishing\footnote{namely: $H^q(\ZZ(n))=0$ if $q>n$ and
 $H^p(X,H^q(\ZZ(n)))=0$ if $p>n$, which follows from the Gersten
 resolution satisfied by the sheaf $H^q(\ZZ(n))$.},
 one gets:
$$
\HH^{2n,n}_\cM(X,\ZZ)
\simeq H^n_{\mathrm{Zar}}(X,\mathcal K_*^M)\simeq CH^n(X).
$$
Then, from \cite[Prop. 3.16]{Deg9}, we get that
 the Gysin morphism associated with $\HH_\ZZ$ and
  a projective morphism in the previous definition coincides in degree $(2n,n)$
  with the usual pushforward on Chow groups.\footnote{Let
  us indicate also that one can go from the case of
 smooth $k$-schemes to arbitrary regular schemes
 of equal characteristics using Popescu's desingularization
 theorem together with Lemma \ref{lm:coh&lim}.}
 (See also Example \ref{ex:comput_Gysin}(2)).
\item \textit{The spectrum $\KGL$ and algebraic K-theory}.--
 the Gysin morphism on $\KGL$ associated with a projective
 morphism of regular schemes agrees with the usual functoriality:
 see Example \ref{ex:comput_Gysin}(3).
\end{enumerate}
\end{ex}

\begin{prop} \label{prop:basic_Gysin1}
Assume $\E$ is absolutely pure.
\begin{enumerate}
\item For any composable projective lci morphisms $f$ and $g$: $f_*g_*=(fg)_*$.
\item For any projective lci morphism $f:Y \rightarrow X$
 and any couple $(x,y)$ in $\E^{**}(X) \times \E^{**}(Y)$,
$$
f_*\big(f^*(x).y\big)=x.f_*(y).
$$
\end{enumerate}
\end{prop}
\begin{proof}
Point (1) follows formally from the preceding facts
 (for details, see the proof of \cite[Prop. 5.14]{Deg8},
  using: relations (1), (2) of Paragraph \ref{num:basic_Gysin_projection},
  point (1) of Proposition \ref{prop:basic_Gysin_closed} and
  the preceding lemma).

Point (2) follows from relation (2) of Proposition \ref{prop:basic_Gysin_closed}
 for the case of a closed immersion together with \cite[Cor. 5.18]{Deg8}
 for the case of the projection of $\PP^n_X/X$.
\end{proof}

\begin{num} \label{num:excess_square}
 Consider a cartesian square 
$$
\xymatrix@=10pt{
Y'\ar^q[r]\ar_g[d]\ar@{}|\Delta[rd] & X'\ar^f[d] \\
Y\ar_p[r] & X
}
$$
in $\base$ such that $p$ is a projective lci morphism.

Choose an $X$-embedding $\nu$ of $Y$ into a suitable projective bundle $\PP^n_X$.
The preceding square induces a cartesian morphism of closed pairs
 $\Theta:(\PP^n_{X'},Y') \rightarrow (\PP^n_X,Y)$ and we denote by $\xi$
 the vector bundle over $Y'$ defined by formula \eqref{eq:excess_bdl}
 with respect to $\Theta$.
 According to \cite[Prop. 6.6(c)]{Ful} the vector bundle $\xi$ 
 is independent up to isomorphism of the choice of $\nu$
 and we call it the \emph{excess intersection bundle} associated with 
 the square $\Delta$.
\end{num}
\begin{prop} \label{prop:basic_Gysin2}
Assume $\E$ is absolutely pure and consider the above notations.
Let $e$ be the rank of $\xi$ over $Y'$.

Then for any $y \in \E^{**}(Y)$, one has:
 $f^*p_*(y)=q_*\big(c_e(\xi).g^*(y)\big)$.
\end{prop}
This follows easily from the definition of $\xi$ and point (3)
 of Proposition \ref{prop:basic_Gysin_closed}.

As usual, one says the square $\Delta$ is transversal 
 if $e=0$ -- \emph{i.e.} for any point $a$ in $Y'$,
 the dimension of $q$ at $a$ is equal to the dimension of
 $p$ at $g(a)$. In this case, the preceding formula
 reads simply: $f^*p_*=q_*g^*$.

As a companion of the Gysin morphism,
 one gets another kind of characteristic classes,
 the cobordism classes. The definition obtained here
 is valid in a slightly better generality than that of
 \cite[\textsection 5.27]{Deg8}.
\begin{df}\label{df:cobordism}
Assume $\E$ is an absolutely pure oriented ring spectrum.

Given any base scheme $S$ in $\base$
 and any projective $S$-scheme $X$ in $\base$ with structural
 morphism $p$,
 we define the \emph{cobordism class} of $X/S$ as the following element
 of $\E^{**}(S)$:
$$
[X/S]:=p_*(1).
$$
\end{df}
Note that if $X/S$ has constant relative dimension $d$,
 the class $[X/S]$ has cohomological degree $(-2d,-d)$.

\begin{ex}\label{ex:Myschenko}
Let us now recall the explicit computation of the cobordism class
 of the linear projective $n$-plane obtained in \cite[5.31]{Deg8}.
 In fact, we will actually correct a sign mistake \emph{loc. cit.} that was
 indicated to us by the referee.

According to Remark \ref{rem:compute_Gysin_proj_bdl},
 the key point is to compute the \emph{fundamental class of the diagonal}
 $\delta$ of the $X$-scheme $P=\PP^n_X$; 
 in other words, the element $\delta_*(1)$, with
 the notation of Definition \ref{df:Gysin_triangle}.

The main point used to compute that class is the classical fact
 that the diagonal subscheme $\delta(P)$ of $P \times_X P$ is the
 zero locus of a canonical section of the vector bundle
$$
\lambda_1^\vee \otimes \xi_2=\uHom(\lambda_1,\xi_2)
$$
where $\lambda_1$ (resp. $\xi_2$) is the canonical line (resp. quotient) 
 bundle on the first (resp. second) factor of $P \times_X P$,
 seen as a vector bundle over $P \times_X P$
 (see \cite[Ex. 5.29]{Deg8}).
 Therefore according to Corollary \ref{cor:fdl&Chern_transversal_section},
 one gets:
$$
\delta_*(1)=c_n\big(\lambda_1^\vee \otimes \xi_2\big).
$$
Then, if $(a_{ij})_{i,j \in \NN}$ denotes the coefficients of the formal
 group law associated with the orientation $c$ of $\E$, one gets the following
 formula\footnote{We refer the reader to \cite[proof of 5.30]{Deg8} for
 a proof.}:
$$
\delta_*(1)=\sum_{0\leq i,j\leq n} \pi^*(a_{1,i+j-n}).c^i.d^j.
$$
Therefore, with the notations of Remark \ref{rem:compute_Gysin_proj_bdl},
 the matrix $M$ corresponding to $\delta_*(1)$ is:
$$
\left(
\raisebox{0.5\depth}{
\xymatrix@=0.1ex{
0\ar@{.}[rr]\ar@{.}[dd]
 & & 0\ar@{.}[lldd]
 & 1\ar@{-}[llldddd] \\
&&& a_{1,1}\ar@{-}[llddd]\ar@{.}[ddd] \\
0 & & & \\
& & & \\
1 & a_{1,1}\ar@{.}[rr] & & a_{1,n}
}
}
\right).
$$
Note a basic computation gives that the determinant of this matrix is
 $(-1)^{\lfloor(n+1)/2\rfloor}$. Therefore, using formula
 \eqref{eq:pre_Myschenko} together with the preceding computation,
 we get the following expression of the cobordism class of the projective
 $n$-plane:
\begin{equation} \label{eq:Myschenko}
[\PP^n/S]=(-1)^{\lfloor(n+1)/2\rfloor}.\left|
\raisebox{0.5\depth}{
\xymatrix@=0.1ex{
0\ar@{.}[rr]\ar@{.}[dd]
 & & 0\ar@{.}[lldd] & 1\ar@{-}[lllddd]
 & a_{1,1}\ar@{-}[lllldddd] \\
&&&& a_{1,2}\ar@{-}[lllddd]\ar@{.}[ddd] \\
0 & & & & \\
1 & & & & \\
a_{1,1} & a_{1,2}\ar@{.}[rrr] & & & a_{1,n}
}
}
\right|.
\end{equation}
Beware that the sign written here corrects the sign in \cite[Cor. 5.31]{Deg8}.
 In fact, this formula coincides with its analogue in topology
 given by the classical Myschenko formula. 

Given an arbitrary projective $S$-scheme $X$,
 we get the following method to compute the cobordism
 class of $X/S$: choose an embedding $i:X \rightarrow \PP^N_S$;
 compute the fundamental class of $i$ in terms of the canonical
 basis $(c^i, 0 \leq i \leq N)$ of $\E^{**}(\PP^n_S)$
 as a free  $\E^{**}(S)$-module, where $c$ is the first Chern
 class of dual canonical invertible bundle:
$$
\eta_{\PP^n_S}(X)=\sum_{i=0}^N x_i.c^i, x_i \in \E^{**}(S).
$$
Then the projection formula yields the following expression:
$$
[X/S]=\sum_{i=0}^N x_i.[\PP^{N-i}/S].
$$
\end{ex}

\begin{rem}
Using the Riemann-Roch formula, we will reprove Quillen formula 
 which computes the cobordism class of
 an arbitrary projective bundle
 (see Ex. \ref{ex:Quillen_cobord_formula})
\end{rem}

\subsection{Uniqueness}

In our setting, we can extend the uniqueness statement
 of \cite[4.1.4]{Panin2}, to obtain the following characterization
 of the Gysin morphisms we have introduced.
\begin{thm}\label{thm:compare_Gysin}
Assume $\E$ is an absolutely pure ring spectrum over $\base$.
 For any scheme $S$, we will denote by $\eta_S$
 the class in $\E ^{2,1}(\PP^1_S)$ obtained by $\PP^1$-suspension
 of the unit $S^0 \rightarrow \E$ of the ring spectrum structure.

Suppose given for any projective lci morphism $f:Y \rightarrow X$
 a morphism of bigraded abelian groups (non necessarily homogeneous):
$$
f_\star:\E^{**}(Y) \rightarrow \E^{**}(X)
$$
such that the following properties hold:
\begin{enumerate}
\item $(fg)_\star=f_\star g_\star$;
\item $f_\star$ is $\E^{**}(X)$-linear: $f_\star(f^*(x).y)=x.f_\star(y)$;
\item for any square $\Delta$ as in \ref{num:excess_square},
 which is in addition transversal, $f^*p_\star=q_\star g^*$;
\item if $i$ is a closed immersion with complementary open immersion $i'$,
 $\mathrm{Im}(i_\star)=\mathrm{Ker}(i^{\prime*})$.
\end{enumerate}
For any integer $n>0$, let $\lambda_n$ be the canonical line bundle
 on $\PP^n_S$ and $s_n$ be its zero section. Let us put:
$$
c_{n,S}=s_n^*s_{n\star}(1).
$$
Then the following conditions are equivalent:
\begin{enumerate}
\item[(i)] the sequences $c_S=(c_{n,S})_{n >0}$ indexed by a scheme $S$
 form an absolute orientation of $\E$;
\item[(ii)] for any $n>0$, $c_{n,S}$ has bidegree $(2,1)$ and $c_{1,S}=\eta_S$
 in $\E^{2,1}(\PP^1_S)$.
\end{enumerate}
When these equivalent conditions are fulfilled, we get in addition
 for any scheme $S$ which admits an ample family of line bundles:
\begin{enumerate}
\item[(5)] for any section $s$ of a line bundle $L/S$,
 $s^*s_\star(1)=c_1(L)$ where the right hand side is the first Chern
 class associated with the orientation $c$ of $\E$ 
 given by the above condition (i).
\end{enumerate}
Finally, if we assume that condition (5) holds for any scheme $S$
 in $\base$, we get:
\begin{itemize}
\item for any projective lci morphism $f$,
 one has: $f_\star=f_*$, the last morphism
 being the Gysin morphism defined in \ref{df:Gysin}
 with respect to $(\E,c)$.
\end{itemize}
\end{thm}

\begin{rem}
To summarize: the Gysin morphism of an absolutely pure
 oriented ring spectrum $(\E,c)$ is uniquely characterized by properties
 (1)-(5).
\end{rem}

\begin{proof}
The equivalence between (i) and (ii) is obvious from Definition \ref{df:orientation}
 and the fact that Point (3) implies: $\iota_n^*(c_{n+1}^S)=c_n^S$.

Let us assume these equivalent conditions are satisfied.
To prove (5), using homotopy invariance, we can assume $s$ is the zero section
 of $L/S$. Because $S$ admits an ample family of line bundles,
 we can use Jouanolou's trick which reduces us to the case where
 $S$ is affine.
 Then $L$ is generated by a finite number of its sections and one
 can find an immersion $i:S \rightarrow \PP^n_S$ such that $L=i^{-1}(\lambda_n)$.
 Using (3), one reduces to the case where $S$ is $\PP^n_S$, $L=\lambda_n$
 which holds by definition.

Let us now prove the final point. Given any projective lci morphism $f$,
 we put $f_\flat=f_*-f_\star$. We have proved previously that
 the Gysin morphisms $f_*$ satisfies properties (1)-(4) stated above.
 In particular, $f_\flat$ satisfies properties (1)-(3) and in the situation of (4),
 we get:
\begin{enumerate}
\item[(4')]  $\mathrm{Im}(i_\flat) \subset \mathrm{Ker}(i^{\prime *})$.
\end{enumerate}
\underline{Case of closed immersions}: assuming $f=i:Z \rightarrow S$ is a closed immersion,
 we show $i_\flat=0$.
Using (3) for $?_\flat$,
 the deformation diagram \eqref{eq:deformation} induces a commutative diagram:
$$
\xymatrix@R=14pt{
\E^{**}(Z)\ar_{i_\flat}[d] & \E^{**}(\AA^1_Z)\ar^{s_0^*}[r]\ar_{s_1^*}[l]\ar^{k_\flat}[d]
 & \E^{**}(Z)\ar^{s_\flat}[d] \\
\E^{**}(X) & \E^{**}(D_ZX)\ar^{\sigma_0^*}[r]\ar_{\sigma_1^*}[l] & \E^{**}(N_ZX).
}
$$
We use the following lemma:
\begin{lm}
Consider the notations above and let $k':(D_ZX-\AA^1_Z) \rightarrow D_ZX$ be the complementary
 open immersion to $k$. Then the morphism: $(k^{\prime*},\sigma_0^*)$ is injective.
\end{lm}
Indeed, let $x \in \E^{**}(D_ZX)$ being an element of the kernel.
 Because of property (4), we get a cohomology class $y$ such that $x=k_\star(y)$.
 On the other hand, because of (3), $\sigma_0^* k_\star=s_\star s_0^*$.
 But $s_\star$ is a split monomorphism (because of (1)), and $s_0^*$ is an isomorphism.
 Thus we deduce $y=0$ and this concludes.

With the help of this lemma, and the fact $k^{\prime*}k_\flat=0$ from (4') above,
 we see that it is sufficient to prove $\sigma_0^*k_\flat=0$. Thus,
 we are reduced to show that $s_\flat=0$, where again $s$ is the zero
 section of the vector bundle $E=N_ZX$ over $Z$.

Using the splitting principle applied to the vector bundle $E/Z$ 
 and property (1), we reduce to the case where $s$
 is the zero section of a line bundle $p:L \rightarrow Z$.
 According to (2), $s_\flat(x)=s_\flat(1).p^*(x)$. The fact $s_\flat(1)=0$
 is nothing else than assumption (5).

\bigskip

\underline{Case of a projective bundle}: let $p:\PP^n_S \rightarrow S$
 be the canonical projection, and let us show $p_\flat=0$.
 Let $\delta$ be the diagonal immersion of $\PP^n_S/S$.
 Let us consider the bigraded ring $A=\E^{**}(S)$.
 Below any $A$-module is assumed to be a bigraded $A$-module,
  but we do not require morphisms of bigraded $A$-modules are homogeneous.
  The symbol $\otimes_A$ means the tensor product of bigraded $A$-modules.
 Then,
  according to the projective bundle theorem, $M=\E^{**}(\PP^n_S)$
  is a free $A$-module of finite rank. 
 Thus it is a rigid object of the category of $A$-modules.
 On the other hand, one gets using (2) the following morphisms
 of $A$-modules:
\begin{align*}
p_\star \delta^*: & \; M  \otimes_A M \longrightarrow A \\
\delta_\star p^*: & \; A \longrightarrow M \otimes_A M.
\end{align*}
Using (1), (2) and (3), we get that the following composite is the identity
 morphism:
$$
M \xrightarrow{M \otimes_A \delta_\star p^*} M  \otimes_A M \otimes_A M
 \xrightarrow{p_\star \delta^* \otimes_A M} M.
$$
The same is true when replacing $\star$ by $*$.
 Thus, as we have already proved that $\delta_\star=\delta_*$,
 we deduce that:
$$
(p_\flat \delta^* \otimes_A M) \circ (M \otimes_A \delta_* p^*)=0.
$$
As $\delta_*$ and $p^*$ are split monomorphisms, we deduce that
 $p_\flat \delta^* \otimes_A M=0$.
 This allows to conclude because $M$ is faithfully flat over $A$ and
 $\delta^*$ is a split epimorphism.
\end{proof}

\begin{ex} \label{ex:comput_Gysin}
As said before, the case where $\base$ is the category of
 smooth $k$-schemes was already obtained by Panin and Smirnov,
 later published by Panin in \cite[4.1.4]{Panin2}.\footnote{See \cite{PS}
 for the original work of Panin and Smirnov.
 I thank the referee for attracting my attention on the paper of Smirnov \cite{Smir1}.}
\begin{enumerate}
\item It is worthwhile to mention that this result applies in particular to any
 Mixed Weil cohomology theory: example \ref{ex:abs_ring_spectra}(1),
 where $\base$ is the category of regular $k$-schemes.
 Thus, there is only one way to define pushforwards on any such cohomology
 satisfying conditions (1)-(5) -- and in particular compatible with a
 well behaved first Chern class.
\item It can also be applied to the usual Chow groups:
 in \cite{Deg8}, we have shown that the category of stable motivic complexes
 $DM(k)$ is endowed with a $t$-structure, called the homotopy $t$-structure\footnote{extending
 the homotopy $t$-structure of Voevodsky that he defined on $DM_-^{eff}(k)$.}
 and whose heart is the category of cycle modules defined by Rost.
 In particular,
 the cycle module corresponding to Milnor K-theory defines an object
 in $DM(k)$
 which is nothing else than the $0$-th cohomology object of the unit
 in $DM(k)$
 (see \emph{loc. cit.}, Th. 5.11). Using the forgetful functor
$$
DM(k) \rightarrow SH(k)
$$
we can see this object as a ring spectrum in $SH(k)$: by definition,
 it corresponds
 to the unramified Milnor-K-theory sheaf $\mathcal K_*^M$,
 seen over the Nisnevich site of smooth $k$-schemes,
 and represents the Chow group in $SH(k)$, with its product structure.
 Note also this is a direct factor of the motivic Eilenberg-MacLane spectrum
 $\mathbf H_{\mathcal M,k}^\ZZ$ representing motivic cohomology (cf. \ref{ex:abs_ring_spectra}):
 it corresponds to cut-out the groups of bidegree $(2*,*)$. 
 As usual, one extends $\mathcal K_k$
 to a $k$-absolute ring spectrum $\mathcal K$ by taking pullbacks
 within the fibred
 category $SH$.

Then the preceding theorem applies to the latter ring spectrum $\mathcal K$,
 when $\base$ is the category of regular $k$-schemes
  (because of Lemma \ref{lm:coh&lim} and the fact that usual Chow groups
   commute with limits of regular $k$-schemes).
 Thus, it says that prescribing the orientation on $CH^*$,
 there exists a unique way of defining pushforwards satisfying properties
 (1)-(5).
\item More interestingly, the proposition applies to $\KGL$, the
 absolute ring spectrum representing $K$-theory when $\base$ is the category
 of regular schemes (Example \ref{ex:abs_ring_spectra}(3)).
 Recall from \cite{QuiK} the corresponding cohomology has well defined
 pushforward that satisfies conditions (1)-(4) of the proposition.
Moreover, with the choice of orientation of Example \ref{ex:abs_ring_spectra_orientation}(4),
 condition (5) of the above Theorem has been proved by Thomason in \cite[Th. 3.1]{Thomason}.
Thus,  the isomorphism \eqref{eq:iso_KGL&Ktheory}
  is covariantly functorial with respect to projective
   morphisms between regular schemes, where on the source we consider
  the pushforward defined above with respect to the absolute
  oriented ring spectrum ($\KGL$,$c^\KGL$)
  and on the aim we consider Quillen pushforward.
 
 Note also that we need Thomason excess intersection formula only
 in the case of line bundles:
 the machinery developed here gives a proof in the general case.
 Note also that if we make the ``hypoth\`ese paresseuse''
 (referred to by Thomason in the end of the introduction
 of \emph{loc. cit.}) that schemes, in addition to being regular,
 admit an ample line bundle, 
 we even get a new proof of the excess intersection formula for higher K-theory
 -- because (5) is then automatically satisfied by our choice of orientation
 and we can compare our Gysin morphism with the one defined by Quillen.
\end{enumerate}
\end{ex}
\section{Riemann-Roch formulas}
\label{sec:RR}

\begin{num} \label{num:convention_RR}
In this section,
 we will consider two absolute oriented ring spectra $(\E,c)$, $(\F,d)$
 and a morphism of absolute ring spectra (cf Definition \ref{df:abs_ring_sp}):
$$
\varphi:\E \rightarrow \F.
$$
When considering the constructions of orientation theory as described
 previously for $\E$ (resp. $\F$), we will will put an index
 $\E$ (resp. $\F$) in the notation. However, when no confusion is possible,
 we drop this index.

In the particular case $\E=\F$,
 we will also put an upper-index $c$ (resp. $d$) in the notation 
 of orientation theory (Chern classes, Gysin morphisms, ...) with
 respect to the orientation $c$ (resp. $d$).
\end{num}

\subsection{Todd classes}

\begin{num} \label{num:todd}
Let $S$ be a scheme (in $\base$).
 We deduce from $\varphi$ a morphism of graded rings:
$$
\E^{**}(\PP^\infty_S)
 \xrightarrow{\varphi_{\PP^\infty_S}} \F^{**}(\PP^\infty_S)
$$
using the notations of Paragraph \ref{num:Pinfty}.
According to the projective bundle theorem (\ref{thm:PBF}),
 this corresponds to a morphism of ring:
$$
\E^{**}(S)[[u]] \rightarrow \F^{**}(S)[[t]]
$$
and we denote by $\Psi_S(t)$ the image of $u$ by this map.
In other words,
 this formal power series is characterized by the relation:
\begin{equation} \label{eq:ass_power_series}
\varphi_{\PP^\infty_S}\left(c\right)=\Psi_S(d).
\end{equation}
Note that the restriction of $\varphi_{\PP^\infty_S}(c)$
 to $\PP^0_S$ (resp. $\PP^1_S$) is $0$ (resp. the multiplicative identity
 element of the ring $\F^{**}(\PP^1_S)$) because
 $c$ is an orientation and $\varphi$ is a morphism of
 ring spectra.
 Thus, as $d$ is also an orientation of $\F$,
 we can write $\Psi_S(t)$ as:
$$
\Psi_S(t)=t+\sum_{i>1} \alpha_i^S.t^i
$$
where $\alpha_i^S \in \F^{2-2i,1-i}(S)$.
Obviously, the power series $\Psi_S(t)/t$ is invertible.

We will also consider the commutative monoid $\mathcal M(S)$
 generated by the isomorphism classes of vector bundles over $S$ 
 modulo the relations $[E]=[E']+[E'']$ coming from exact sequences
$$
0 \rightarrow E' \rightarrow E \rightarrow E'' \rightarrow 0.
$$
Then $\mathcal M$ is a presheaf of monoids on $\base$
 whose presheaf of
 abelian groups is the functor $K_0$.

Note that $\F^{00}(S)$, equipped with cup-product,
 is a commutative monoid. We will denote by
 $\F^{00\times}(S)$ the group made by its invertible elements.
\end{num}
\begin{prop} \label{prop:todd}
There exists a unique natural transformation of presheaves
 of monoids over $\base$
$$
\td_\varphi:\mathcal M \rightarrow \F^{00}
$$
such that for any line bundle $L$ over a scheme $X$,
\begin{equation} \label{eq:td_line_bdl}
\td_\varphi(L)=\frac t {\Psi_S(t)}(t=d_1(L)).
\end{equation}
Moreover, it induces a natural transformation 
 of presheaves of abelian groups:
$$
\td_\varphi:K_0 \rightarrow \F^{00\times}.
$$
\end{prop}
\begin{proof}
The proof is very classical:
 the uniqueness statement follows from the \emph{splitting principle}
 while the existence statement follows from the use of \emph{Chern roots}.
 Note also that the relation \eqref{eq:td_line_bdl} is well defined
 because $c_1(L)$ is nilpotent (see Proposition \ref{prop:FGL}).
The final assertion follows from the fact $t/\Psi_S(t)$ is an invertible
 formal power series.
\end{proof}

\begin{rem}
According to the construction of the first Chern classes
 for the oriented ring spectra $(\E,c)$ and $(\F,d)$
 together with Relations \eqref{eq:ass_power_series} and \eqref{eq:td_line_bdl},
 we get for any line bundle $L/S$ the following identity in $\F^{2,1}(S)$:
\begin{equation} \label{eq:todd&chern}
\varphi_S\big(c_1(L)\big)
 =\td_\varphi(-L) \cupp d_1(L).
\end{equation}
\end{rem}

\begin{df} \label{df:todd}
Consider the context and notations of the previous
 proposition.

Given any virtual vector bundle $e$ over $X$,
 the element $\td_\varphi(e) \in \F^{00}(X)$
 is called the \emph{Todd class} of $e$ over $X$ associated
 with the morphism of ring spectra $\varphi$.
\end{df}

\subsection{The case of closed immersions}

\begin{num}
Consider a regular closed immersion $i:Z \rightarrow X$, $U=X-Z$.
As by assumption $\E$ (resp. $\F$) is absolutely pure,
 we can consider the associated refined fundamental class
 $\rfdl^\E_X(Z)$ in $\E^{**}_Z(X)$
 (resp. $\rfdl^\F_X(Z)$ in $\F^{**}_Z(X)$)
 -- Definition \ref{df:fdl_class}.
 The morphism $\varphi$ induces a map in relative cohomology:
$$
\E^{**}_Z(X) \xrightarrow{\varphi_{X,Z}} \F^{**}_Z(X)
$$
According to the definition of the purity isomorphism \eqref{eq:purity},
 we deduce there exists a unique class $\tau_\varphi(X,Z) \in \F^{00}(Z)$
 such that
\begin{equation} \label{eq:tauto_todd}
\varphi_{X,Z}\big(\rfdl^\E_X(Z)\big)=\tau_\varphi(X,Z).\rfdl^\F_X(Z).
\end{equation}
This relation together with the definition of 
 the localization long exact sequence \eqref{eq:Gysin_exact_seq} 
 immediately gives the following commutative diagram:
\begin{equation}\label{eq:RR_closed0}
\begin{split}
\xymatrix{
\E^{**}(U)\ar^{\partial_{X,Z}}[r]\ar_{\varphi_U}[d]
 & \E^{**}(Z)\ar^{i_*}[r]\ar|{\tau_\varphi(X,Z).\varphi_Z}[d]
 & \E^{**}(X)\ar^{\varphi_X}[d] \\
\F^{**}(U)\ar^{\partial_{X,Z}}[r]
 & \F^{**}(Z)\ar^{i_*}[r]
 & \F^{**}(X).
}
\end{split}
\end{equation}
The commutativity of the square on the right hand side
 is a tautological Riemann-Roch formula.
 Our proof of the actual Riemann-Roch theorem in fact lies
 in the computation of the class involved in this formula:
\end{num}
\begin{lm}\label{lm:RR_rfdl}
In the above assumptions, the following relation holds:
 $$\tau_\varphi(X,Z)=\td_\varphi(-N_ZX).$$
\end{lm}
\begin{proof}
As Relation \eqref{eq:tauto_todd} characterizes uniquely
 the class $\tau_\varphi(X,Z)$,
 whatever the regular closed pair $(X,Z)$ is,
 the deformation diagram \eqref{eq:deformation_iso}
 gives the relation:
$\tau_\varphi(X,Z)=\tau_\varphi(N_ZX,Z)$.

Thus, we are reduced to prove that for any scheme $X$
 and any vector bundle $E/X$,
$$
\td_\varphi(-E)=\tau_\varphi(E,X).
$$
 Using again the characterizing relation \eqref{eq:tauto_todd},
 which involves refined Thom classes according to \eqref{rfdl_class&vector_bdl},
 we deduce:
\begin{itemize}
\item for any morphism $f:Y \rightarrow X$ and any vector bundle $E/X$,
 $f^*\tau_\varphi(E,X)=\tau_\varphi(f^{-1}(E),Y)$.
\item for any scheme $X$
 and any exact sequence of vector bundles over $X$,
$$
0 \rightarrow E' \rightarrow E \rightarrow E'' \rightarrow 0,
$$
one has: $\tau_\varphi(E,X)=\tau_\varphi(E',X)\tau_\varphi(E'',X)$.
\end{itemize}
More precisely: for the first relation, one uses
 the compatibility of the refined Thom class with pullback
 and for the second, one applies Proposition \ref{prop:product&thom}.

In other words, one gets a morphism of monoids
$$
M(X) \rightarrow F^{00}(X), [E] \mapsto \tau_\varphi(E,X)
$$
which is contravariantly natural in $X$.
According to the uniqueness statement of Proposition \ref{prop:todd},
 we are reduced to the case of line bundles.

Let $L/X$ be a line bundle and $P$ be its projective completion.
Let $\xi$ be the universal quotient bundle on $P$.
According to Relation \eqref{eq:todd&chern}, one obtains:
$$
\varphi_P\big(c(\xi)\big)
 =\td_\varphi(-\xi) \cupp d(\xi).
$$
Let $s$ be the section of $P/X$ induced by the zero section
 of $L/X$.
 Using the commutativity of the right square of \eqref{eq:RR_closed0}
 when the closed pair $(X,Z)$ is $(P,X)$, we obtain:
$$
\varphi_P\big(s_*(1)\big)
 =s_*(\tau_\varphi(P,X))=s_*(\tau_\varphi(L,X))
$$
According to \eqref{eq:Gysin&rfdl}, \eqref{fdl_class&proj_bdl}
 and \eqref{eq:compute_Thom}: $s_*(1)=c(\xi)$.
 Thus we obtain:
$$
\td_\varphi(-\xi) \cupp s_*(1)=s_*(\tau_\varphi(L,X)).
$$
Let $p:P \rightarrow X$ be the canonical projection.
 Of course, $p \circ s=Id$. Thus,
 it is sufficient to apply $p_*$ to the preceding relation to conclude:
$$
p_*\big(\td_\varphi(-\xi) \cupp s_*(1)\big)
\stackrel{(*)}=p_*s_*\big(s^*\td_\varphi(-\xi)\big)
=s^*(\td_\varphi(-\xi))=\td_\varphi(-L),
$$
where $(*)$ is the projection formula \ref{prop:basic_Gysin_closed}(2).
\end{proof}

Finally the following Riemann-Roch formulas are consequences
 of the commutative diagram \eqref{eq:RR_closed0} in conjunction
 with the preceding lemma.
\begin{thm} \label{thm:RR1}
Under the assumptions of the previous paragraph,
 the following formulas hold:
\begin{align*}
\varphi_X(i_*(z))&=i_*\big(\td_\varphi(-N_ZX) \cupp \varphi_Z(z)\big), \\
\partial_{X,Z}\big(\varphi_U(u)\big)&=  \td_\varphi(-N_ZX) \cupp \varphi_Z\big(\partial_{X,Z}(u)\big).
\end{align*}
\end{thm}

\subsection{The general case}

\begin{num} \label{num:virtual_tangent}
Consider a projective lci morphism $f:Y \rightarrow X$
 (Par. \ref{num:lci_morphisms}).
 Given a factorization \eqref{eq:generic_fact_lci} of $f$,
 we define the \emph{virtual tangent bundle} of $f$ as
 the element of $K_0(Y)$:
$$
\tau_f:=[i^{-1}T_p]-[N_i]
$$
where $T_p$ is the tangent space of $p$
 and $N_i$ is the normal bundle of $i$.
Standard considerations show this definition is independent
 of the chosen factorization (see \cite[VIII, 2.2]{SGA6}). 
\end{num}
\begin{thm} \label{thm:projective_RR}
Consider the notations above.
Then for any element $y \in \E^{**}(Y)$,
 the following formula holds:
$$
\varphi_X(f_*(y))=f_*\big(\td_\varphi(\tau_f) \cupp \varphi_Y(y)\big).
$$
\end{thm}
\begin{proof}
According to the multiplicativity of the Todd class and the compatibility
 of the Gysin morphism with respect to composition,
 this formula can be divided in two cases according to a factorization
 \eqref{eq:generic_fact_lci} of $f$.
 The case of a closed immersion was treated above (\ref{thm:RR1})
 and it remains to consider the case of the canonical projection
 $p$ of a projective space $P=\PP^n_X/X$.

We first treat the case where $\varphi_{\PP^\infty_X}(c)=d$.
 This implies that $\td_\varphi=1$, so we have to prove that
 for $y \in \E^{**}(P)$, one has:
\begin{equation}\label{eq:pre_P^n_RR}
\varphi_X(p_*(y))=p_*\big(\varphi_Y(y)\big).
\end{equation}
 Consider the notations of Paragraph \ref{num:Gysin4projection}
  with respect to both $\E$ and $\F$.
 The Riemann-Roch formula for the diagonal embedding $\delta$ gives the formula:
$$
\varphi_{P \times_X P} \delta_*=\delta_*\varphi_P.
$$
This is sufficient to conclude in view of Definition \ref{df:Gysin_Pn}.
 Let us give more details with the considerations of Remark
 \ref{rem:compute_Gysin_proj_bdl}.
The preceding formula implies that the map
 $\varphi_X:\E^{**}(X) \rightarrow \F^{**}(X)$
 sends each of the coefficients of the matrix $M^\E$
 of $\delta_*(1) \in \E^{**}(P \times_X P)$
 to the corresponding coefficient of the matrix $M^\F$
 of $\delta_*(1) \in \F^{**}(P \times_X P)$.
 Therefore, the same result holds for the inverse matrices
 of $M^\E$ and $M^\F$.
 Thus finally, Formula \eqref{eq:pre_P^n_RR}
 comes from the algorithm to compute $p_*$ given at the end of
 Remark \ref{rem:compute_Gysin_proj_bdl}.

It remains now to consider the case where $\E=\F$, $\varphi=1_\E$,
 for which we adopt the convention of the last paragraph of
 \ref{num:convention_RR}.
 Thus, we have to prove
\begin{equation} \label{eq:RR3}
\forall \alpha \in \E^{**}(P), \ 
p_*^{c}(\alpha)=p_*^{d}(\td_\varphi(T_p).\alpha).
\end{equation}
 We consider again the notations of Paragraph \ref{num:Gysin4projection}
  with respect to both $(\E,c)$ and $(\E,d)$.
 Let $\pi:P \times_X P \rightarrow P$ be the projection to
 the second factor.
 The Riemann-Roch formula for $\delta$ can be read in this case as:
\begin{align*}
\delta_*^{c}(\alpha)
 =\delta^{d}_*\big(\td_\varphi(-N_\delta) \cupp \alpha\big)
 =\delta^{d}_*\big(\delta^*\pi^*(\td_\varphi(-N_\delta)) \cupp \alpha\big)
 =\pi^*(\td_\varphi(-N_\delta)) \cupp \delta^{d}_*(\alpha)
\end{align*}
-- we use Prop. \ref{prop:basic_Gysin1} for the last equality.
Recall from Paragraph \ref{num:Gysin4projection} we have a K\"unneth
 isomorphism $\E^{**}(P \times_X P)=\E^{**}(P) \otimes_A \E^{**}(P)$.
 Through this identification, and the fact $N_\delta=T_p$,
 the class $\pi^*(\td_\varphi(-N_\delta))$ corresponds to
  $1 \otimes_A \td_\varphi(-T_p)$.
Thus the preceding equality leads to the commutativity of the
 following diagrams:
$$
\xymatrix@R=14pt{
A\ar^-{\epsilon_P^{c}}[r]\ar@{=}[dd]
 & \E^{**}(P) \otimes_{A} \E^{**}(P)\ar^{1 \otimes_A \td_\varphi(T_p)}[dd]
 &&
\E^{**}(P)^\vee\ar^-{\mathcal D_P^{c}}[r]\ar@{=}[dd]
 & \E^{**}(P)\ar^{\td_\varphi(T_p)}[dd] \\
 &&\Rightarrow&& \\
A\ar^-{\epsilon_P^{d}}[r]
 & \E^{**}(P) \otimes_{A} \E^{**}(P)
 &&
\E^{**}(P)^\vee\ar^-{\mathcal D_P^{d}}[r]
 & \E^{**}(P).
}
$$
Then formula \eqref{eq:RR3} finally follows from Definition \ref{df:Gysin_Pn}.
\end{proof}

\section{Examples and applications}
\label{sec:examples}

\begin{num}\label{num:orient_in_notations}
In this section, we will adopt the following notations.
 Orientations of an absolute ring spectrum 
 (Def.  \ref{df:orientation})
 will be denoted by letters $c$ (or $c'$, $d$, ...).
 Then the corresponding Chern classes (Def. \ref{df:Chern_classes})
  with be denoted by the same letter with a lower index: $c_1$, $c_2$,...
 The corresponding Gysin morphisms (resp. fundamental classes, residue
 morphisms) will be denoted by $p_*^c$ (resp.$\bar \eta^c_X(Z)$,
 $\partial^c_{X,Z}$).

All spectra in this section will be
 $\base$-absolutely pure ring spectra,
  so we simply say spectra for $\base$-absolutely pure ring spectra.
 Recall that an orientation $c$ of such a spectrum $\E$
 is a family of orientation $c_S$ of $\E_S$ indexed
 by schemes $S$ of $\base$, and stable by pullbacks
 (Def.~\ref{df:orientation}).
 As seen above, such a collection of orientations
 gives rise to a collection $F$ of formal group laws
 $F_S$ on the ring $\E^{**}(S)$, which will simply be
 called the formal group law associated with $(\E,c)$
  -- or just $c$ when $\E$ is clear.
\end{num}

\subsection{Principle of computation}
\label{sec:principle_comput}

\begin{df}
Let $(\E,c)$ and $(\F,d)$ be oriented spectra.

A \emph{pseudo-morphism} of oriented spectra
 $\varphi:(\E,c) \rightarrow (\F,d)$
 is simply a morphism of ring spectra.
We will say $\varphi$ is a \emph{morphism} of oriented spectra
 if for any scheme $S$, one has:
$$
\varphi_{\PP^\infty_S}(c)=d
$$
in $\tilde \F^{2,1}(\PP^\infty_S)$, the reduced cohomology
 of the pointed scheme $\PP^\infty_S$ pointed by the infinite point
 of $\PP^1_S$. \\
We will also say that a pseudo-morphism $\varphi$ is \emph{identical} if,
 as a morphism of ring spectra, it is the identity.
\end{df}

Note that one immediately deduces
 from the construction of Chern classes
 and from the Riemann-Roch formulas \ref{thm:RR1}
 and \ref{thm:projective_RR} the following result:
\begin{prop}
Let $\varphi:(\E,c) \rightarrow (\F,d)$ be a morphism of oriented
 spectra.
\begin{enumerate}
\item For any vector bundle $E/S$ and any integer $n \geq 0$,
 one has:
 $\varphi_S(c_n(E))=d_n(E)$.
\item For any projective lci morphism $f:Y \rightarrow X$
 one has: $\varphi_X \circ f^c_*=f_*^d \circ \varphi_Y$.
\item For any regular closed immersion $i:Z \rightarrow X$,
 $U=X-Z$, one has:
 $\varphi_Z \circ \partial^c_{X,Z}=\partial^d_{X,Z} \circ \varphi_{U}$.
\end{enumerate}
\end{prop}

\begin{ex} \label{ex:RR_source_MGL}
Let $(\E,c)$ be an oriented spectrum.
 According to Proposition \ref{prop:orientation&MGL-modules},
 the choice of $c$ uniquely corresponds to a morphism of
 oriented spectra:
$$
\varphi:\big(\MGL,c^\MGL\big) \rightarrow (\E,c)
$$
where $c^\MGL$ is the canonical orientation of $\MGL$.

Then, according to the preceding proposition,
 the morphism $\varphi$ is compatible with all the constructions
 of orientation theory given in this paper.
\end{ex}

\begin{num} \label{num:comput_todd1}
Given any pseudo-morphism as in the above definition,
 we obviously get
 that $c'=\varphi_{\PP^\infty_S}(c)$ is an orientation
 of $\F$. In particular, the pseudo-morphism $\varphi$ admits a
 canonical factorization:
\begin{equation}\label{eq:facto_pseudo-morph}
(\E,c) \xrightarrow{\tilde \varphi} (\F,c')
 \xrightarrow{\psi} (\F,c)
\end{equation}
such that $\tilde \varphi$ (resp. $\psi$) is a morphism
 (resp. identical pseudo-morphism) of oriented spectra.

According to Definition \ref{df:todd},
 the Todd class of $\varphi$ is equal to the Todd class of
 $\psi$. 
 Thus the computation of the Todd class of a pseudo-morphism
 can always be reduced to the case of an identical
 pseudo-morphism of the aim,
 which corresponds to the effect of 
 \emph{changing the orientation}.
We will describe this case in more details below
 (see in particular \ref{num:basic_change}).
\end{num}

\begin{rem}
According to Proposition \ref{prop:orientation&MGL-modules},
 oriented (ring) spectra correspond to $\MGL$-algebras.
 In the light of this analogy, pseudo-morphisms
 (resp. morphisms) of oriented spectra corresponds 
 to morphisms of rings (resp. $\MGL$-algebras)
 between $\MGL$-algebras.
 In light of this analogy, the preceding factorization
 can be understood as follows: $\tilde \varphi$ corresponds
 to a morphism of $\MGL$-algebras while $\psi$ corresponds
 to a change of the $\MGL$-algebra structure on $\F$.
\end{rem}

\begin{num}
Recall
 that given a ring $R$ and formal group laws $F$, $G$
 with coefficients in $R$, a morphism $\Phi:(R,F) \rightarrow (R,G)$
 of formal group laws is a power series
$\Phi(t) \in R[[t]]$ of positive valuation such that
$$
\Phi(F(x,y))=G(\Phi(x),\Phi(y))
$$
in $R[[x,y]]$. Such a morphism is an isomorphism if and only
 if $\Phi$ admits a composition inverse -- equivalently,
  $\Phi'(0)$ is invertible in $R$.
 It is called a \emph{strict} isomorphism if $\Phi'(0)=1$.

According to the conventions of this section,
 a morphism of formal group laws arising from 
 orientations of an (absolute) ring spectrum will be
 a family of morphisms indexed by schemes of $\base$
 and stable by pullbacks.
\end{num}
\begin{prop} \label{prop:chg_orient&iso_fgl}
Let $\E$ be a ring spectrum
 and $c$ be an orientation of $\E$ with associated 
 formal group law $F_c$.
Consider the following sets:
\begin{enumerate}
\item the orientations $d$ of $\E$ (equivalently:
 the identical pseudo-morphisms of ring spectra with
 source $(\E,c)$);
\item the strict isomorphisms $\Phi$ of formal group laws
 with source $F_c$ such that for any scheme $S$,
 $\Phi_S(t)$ can be written as a power series
 of the form:
$$
t+\sum_{i>1} \alpha^S_i.t^i
$$
where $\alpha^S_i$ is an element of $\E^{2-2i,1-i}(S)$.
\end{enumerate}
Then the map
$$
(2) \xrightarrow{(*)} (1), \Phi \mapsto \Phi(c)
$$
is a well defined bijection.
\end{prop}
\begin{proof}
We prove the map $(*)$ is well defined.
 First note that, because of the condition on
 the degree of the coefficients of $\Phi_S$,
 the cohomology class $\Phi_S(c)$ is of degree $(2,1)$.
 The fact it is an orientation of $\E_S$ simply follows
  from the form of $\Phi_S(t)$ which implies:
	 $\Phi_S(0)=0$ and $\Phi_S'(0)=1$.

To prove that $(*)$ is a bijection, we construct
 its inverse.
 Let $d$ be an orientation of $\E$ and $S$ be a base scheme.
 Then $d_S$ is an element of the bigraded algebra
$$
\E^{**}(\PP^\infty_S) \simeq \E^{**}(S)[[c]].
$$
This means there is a unique power series
$$
\Phi_S(t)=\sum_{i \geq 0} \alpha_i^S.t^i
$$
with coefficients in $\E^{**}(S)$ such that $d_S=\Phi_S(c_S)$.
Because $d_S$ has bidegree $(2,1)$, we deduce that
 necessarily, $\alpha_i^S$ has bidegree $(2-2i,1-i)$.
 Moreover, because $d_S$ is an orientation of $\E_S$,
 we deduce that $\alpha_0^S=0$ and $\alpha_1^S=1$.
 By uniqueness of $\Phi_S$, we deduce that the family
 $\Phi=(\Phi_S)_{S \in \base}$ is stable by pullbacks.

Let us consider the Segre embedding (Par.~\ref{num:FGL}):
$$
\sigma:\PP^\infty_S \rightarrow \PP^\infty_S \times_S \PP^\infty_S.
$$
By definition of the formal group law $F_c$
 (resp. $F_{d}$) associated with $c$ (resp. $d$)
 we obtain (dropping the reference to the base $S$):
\begin{align*}
\sigma^*(c)=F_c(c',c''), 
\text{ resp. } \sigma^*(d)=F_d(d',d'').
\end{align*}
where, on the left hand side $c', c''$ (resp. $d', d''$) corresponds
 to the first Chern class associated with the orientation $c$ (resp. $d$)
 of the canonical line bundle on the first (resp. second) factor
 of $\PP^\infty_S \times_S \PP^\infty_S$.
 Thus, we obtain that $\Phi_S$ is a strict isomorphism
 from $F_c$ to $F_d$:
$$
\Phi(F_c(c',c''))=\Phi(\sigma^*(c))=\sigma^*(\Phi(c))
 =\sigma^*(d)
 =F_d(d',d'')=F_d(\Phi(c'),\Phi(c'')).
$$
Therefore $\Phi_S(t)$ is a strict isomorphism
 $(\E^{**}(S),F) \rightarrow (\E^{**}(S),F_d)$
 of formal group laws.

The uniqueness of $\Phi_S$ shows that we have indeed constructed
 an inverse map to $(*)$.
\end{proof}

\begin{rem}
This proposition is essentially an elaboration
 of the construction done in Paragraph \ref{num:todd}.
 Note to be more precise that, given an identical
 pseudo-morphism $\varphi:(\E,c) \rightarrow (\E,d)$,
 the power series $\Phi(t)$ obtained by applying the previous
 proposition is the composition inverse of the power series
 $\Psi(t)$ obtained in \ref{num:todd}. Thus, in fact,
 we get reciprocal strict isomorphisms of formal group laws
 on $\E^{**}$:
\begin{equation} \label{eq:associated_isos_fgl}
\Phi:F_c \leftrightarrows F_d:\Psi
\end{equation}
From a categorical perspective, having fixed
 an oriented spectrum $(\E,c)$,
 the proposition defines a covariant functor $(*)$
 from the category of orientations on $\E$ with
 morphisms the identical pseudo-morphisms
 to the category of ($\base$-families) of formals group laws
 on the bigraded ring $\E^{**}$.

The functor $(*)$ induces a bijection of groupoids
 when one restricts morphisms on the target category 
 to strict isomorphisms satisfying the condition on
 the degrees of point (2) above.
\end{rem}

\begin{rem}
According to this proposition,
 two orientations on a given (orientable ring) spectrum $\E$
 necessarily yields (strictly) isomorphic formal group laws.
 Thus, there is a uniquely defined isomorphism class
 of formal group law associated with an orientable
 ring spectrum $\E$ and we can safely qualify such
 a ring spectrum as being additive, multiplicative, etc.
 Note that it is not known whether any formal group law
 can be realized as the formal group law 
 associated with an orientable ring spectrum.
 However, interesting new examples are provided
  in \cite[Th. A]{LYZ}.
\end{rem}

\begin{cor} \label{cor:uniqueness_orientation}
Let $\E$ be an absolute ring spectrum satisfying the following
 property:
\begin{itemize}
\item[(Ann)] For any integer $i<0$,
 and any scheme $S$, $\E^{2i,i}(S)=0$.
\end{itemize}
Then the following assertions hold:
\begin{enumerate}
\item[(i)] If an orientation exists on $\E$ over a scheme $S$,
 it is unique. Moreover, if $\E$ is rational,\footnote{\emph{i.e.} $\E$ is isomorphic
 to its rationalization.} then the formal
  group law associated with an orientation on $\E$ is necessarily
  additive.
\item[(ii)] Assume $\E$ is oriented with orientation $c$
 and $(\F,d)$ is an oriented spectrum. 
Then, any morphism of ring spectra $\varphi:\F \rightarrow \E$
 automatically satisfies $\varphi(d)=c$.
\end{enumerate}
\end{cor}

\begin{ex} \label{ex:uniqueness_orientation}
Assumption (Ann) is most common in algebraic geometry. It is fulfilled in
 any of the examples (1), (2), (5) of \ref{ex:abs_ring_spectra_orientation}.

Over an algebraically closed field $S=\spec k$,
 or more generally any scheme $S$ such that $-1$ is a sum of squares
 in all residue fields, we get from a theorem of Morel
 (cf \cite[16.2.14]{CD3})
 that any $S$-absolute rational ring spectrum satisfying assumption
 (Ann) is uniquely oriented and the corresponding orientation is
  additive.
\end{ex}

\subsection{Change of orientation}\label{sec:change_orientation}

\begin{num}  \label{num:basic_change}
Let $\E$ be  a ring spectrum
 equipped with two orientations $c$, $d$
 whose respective associated formal group laws are $F_c$ and $F_d$.
 This is pictured by the following identical pseudo-morphism
 of oriented ring spectra:
$$
\varphi:(\E,c) \rightarrow (\E,d).
$$
 Let us consider the reciprocal isomorphism of formal group laws
 \eqref{eq:associated_isos_fgl} uniquely associated with $\varphi$
 in Proposition \ref{prop:chg_orient&iso_fgl}.
 Recall the convention: $d=\Phi(c)$, $c=\Psi(d)$.

Then, from Proposition \ref{prop:todd},
 we get two expressions of the Todd class
 of a line bundle $L/S$ associated with $\varphi$:
$$
\td_{\varphi}(L)=\left.\frac t {\Psi_S(t)}\right|_{t=d_1(L)}
 =\left.\frac {\Phi_S(t)} t\right|_{t=c_1(L)}.
$$
Another way of saying this is the relation:
$$
d_1(L)=\td_\varphi(L).c_1(L).
$$
\end{num}

\begin{num} \label{num:log&exp}
In general, it is not easy to compute the strict isomorphism
 associated with a change of orientations on a ring spectrum.
However, let us recall that for any $\QQ$-algebra $A$,
 given any formal group law $F$ with coefficients in $A$,
 there exists a unique
 strict isomorphism between $F$ and the additive formal group law,
 called the \emph{logarithm of $F$}:
$$
\log_F:F \rightarrow F_{add}.
$$
There exists a well known formula for the logarithm.
 First, one attach to $F$ a formal differential form
$$
\omega_F(x)
=\left(\left.\frac{\partial F(x,y)} y\right|_{y=0}\right)^{-1}.dx
$$
and then one defines the logarithm as the primitive of
 this differential form:
\begin{equation} \label{eq:formula_log_F}
\log_F(t)=\int \omega_F(x).
\end{equation}
The composition inverse for this power series
 is usually called the \emph{exponential of $F$} and denoted
 by:
$$
\exp_F:F_{add} \rightarrow F.
$$
As a consequence, we get the well-known fact
 that any two formal group laws
 with coefficients in $A$ are uniquely strictly isomorphic.

This gives the following formula for computing
 the Todd class associated with a change of orientations.
One determines the formal group laws $F_c$ and $F_d$ associated
 with $c$ and $d$. Then the strict isomorphisms
 $\Phi:F_c \leftrightarrows F_d:\Psi$ are given
 by the power series:
$$
\Phi(t)=\exp_{F_d} \circ \log_{F_c}(t),
 \Psi(t)=\exp_{F_c} \circ \log_{F_d}(t),
$$
\end{num}

\begin{ex}
It is easy to compute the Todd class using the splitting principle.
Let us consider the abstract case of a rational oriented ring spectrum 
 $(\E,c)$ with an abstract formal group law
$$
F(x,y)=\sum_{i,j} a_{ij}x^iy^j,
$$
and assume $\bar c=\log_F(c)$ is the canonical orientation corresponding
 to the additive formal group law (cf. the preceding examples).

It is convenient to denote by $p_i=[\PP^i]$ the cobordism class
 of the projective space of dimension $i$ -- whose expression in terms
 of the coefficients $a_{1,i}$ can be found in Example \ref{ex:Myschenko}.
 Indeed, according to the Myschenko formula,
 one gets:
$$
\log_F(t)=\sum_{i=0}^\infty \frac{p_i}{i+1}.t^{i+1}.
$$
Then,
 associated with the identical pseudo-morphism $(\E,c) \rightarrow (\E,\bar c)$,
 the Todd class of an arbitrary vector bundle $E/S$ of rank $n$ 
 can be expressed as a polynomial of the $p_i$ and of the 
 Chern classes $\bar c_i:=\bar c_i(E)$ for $1 \leq i \leq n$.
 The first few terms are as follows:
\begin{align*}
&\td(E)=1 + \left(\frac{1}{2} p_{1}\right) \bar c_{1}
+ \left(\frac{3}{4} p_{1}^{2} - \frac{2}{3} p_{2}\right) \bar c_{2}  \\
& + \left(-\frac{1}{4} p_{1}^{2} + \frac{1}{3} p_{2}\right) \bar c_{1}^{2} +
\left(-\frac{7}{8} p_{1}^{3} + \frac{5}{3} p_{1} p_{2} - \frac{3}{4} p_{3}\right)
\bar c_{1} \bar c_{2} \\
&\qquad\qquad\qquad\qquad\quad+ \left(-\frac{13}{16} p_{1}^{4} + 2 p_{1}^{2} p_{2} -
 \frac{5}{4} p_{1} p_{3} - \frac{1}{3} p_{2}^{2} + \frac{2}{5} p_{4}\right)
\bar c_{2}^{2} \\
& + \left(\frac{1}{4} p_{1}^{3} - \frac{1}{2} p_{1} p_{2} +
\frac{1}{4} p_{3}\right) \bar c_{1}^{3} + \left(\frac{11}{8} p_{1}^{4} - \frac{43}{12}
p_{1}^{2} p_{2} + \frac{17}{8} p_{1} p_{3} + \frac{8}{9} p_{2}^{2} -
\frac{4}{5} p_{4}\right) \bar c_{1}^{2} \bar c_{2} + \cdots
\end{align*}
Note we can also give an expression of the Todd class
 in function of the original Chern classes $c_i= c_i(E)$:
\begin{align*}
&\td(E)=1 + \left(\frac{1}{2} p_{1}\right) c_{1} 
 + \left(\frac{1}{4} p_{1}^{2} - \frac{2}{3} p_{2}\right) c_{2} \\
&+ \left(\frac{1}{3} p_{2}\right) c_{1}^{2}
+ \left(\frac{1}{6} p_{1} p_{2}- \frac{3}{4} p_{3}\right) c_{1} c_{2}
+ \left(-\frac{1}{4} p_{1} p_{3}
 +\frac{1}{9} p_{2}^{2} + \frac{2}{5} p_{4}\right) c_{2}^{2} \\
& + \left(\frac{1}{4}
p_{3}\right) c_{1}^{3} + \left(\frac{1}{8} p_{1} p_{3} - \frac{4}{5} p_{4}\right)
c_{1}^{2} c_{2} + \left(-\frac{3}{10} p_{1} p_{4} + \frac{1}{12} p_{2} p_{3}\right) c_{1} c_{2}^{2} \\
&+ \left(-\frac{2}{15} p_{2} p_{4} + \frac{1}{16} p_{3}^{2}\right) c_{2}^{3} + \cdots
\end{align*}
\end{ex}

As an illustration of the Riemann-Roch formula
 in the case of a change of orientations,
 we give the following simple proof of a formula due to Quillen
 (in complex cobordism, \cite{Qui}):
\begin{thm}\label{thm:Quillen_formula}
Let $E$ be a vector bundle of rank $n+1$ over a scheme $S$,
 $P=\PP(E)$ be the associated projective bundle,
 $p:P \rightarrow S$ be the canonical projection and
 $\lambda=\mathcal O(1)$ be the canonical dual line bundle on $P$.

Put $A=\MGL^{**}(S)$
 and let $c$ be the canonical orientation of $\MGL$ over $S$.
Let $F(x,y) \in A[[x,y]]$ be the formal group law associated with $(\MGL,c)$.

Then for any polynomial $\phi(t) \in A[t]$, the following formula holds in
 $A \otimes_\ZZ \QQ$:
$$
p_*\big(\phi(c_1(\lambda))\big)=\res_t\left(\frac{\phi(t).\omega_F(t)}{\prod_i F(t,l_i)}\right)
$$
where $l_i$ are the Chern roots of $E$ with respect to the orientation $c$.
\end{thm}
\begin{proof}~[D., Levine, Vishik]\footnote{I thank
 a lot M.~Levine and A.~Vishik
 for helping me to finish this proof.}
Let $\log_F(t)$ (resp. $\exp_F(t)$)
  be the logarithm (resp. exponential) associated with the formal group law $F$
  -- see above.
Then the class $\bar c=\log_F(c)$ is an orientation of $\MGL_\QQ$
 whose formal group law is additive.
 Let us denote by $\bar p_*$ the Gysin morphism associated
 with $p$ with respect to $\bar c$.

Then the Riemann-Roch formula for the identical pseudo-morphism
$$
(\MGL_\QQ,c) \rightarrow (\MGL_\QQ,\bar c)
$$
 and the morphism $p$ reads:
\begin{equation} \label{eq:thm:Quillen}
p_*\big(\phi(d)\big)=\bar p_*\big( P(d).\td(T_p) \big)
\end{equation}
where we have put $d=c_1(\lambda)$ and $T_p$ denotes the tangent bundle
 of $P/S$.
The following lemma is a reformulation of a well known formula
 in the (classical) theory of Chern classes:
\begin{lm}
Let $(\E,\bar c)$ be an additive oriented ring spectrum
 over $S$.
Consider the total Chern class of $E$:
$$
\bar c_t(E)=\sum_{i \geq 0} \bar c_i(E).t^i.
$$
as an invertible power series
 with coefficients in $\E^{**}(S)$. Put $\bar d=\bar c_1(\lambda)$.

Then for any power-series $\psi(t)$ with coefficients in $\E^{**}(S)$,
 one has the following equality in $\E^{**}(S)$:
$$
\bar p_*(\psi(\bar d))
=\res_t\left(\frac{\psi(t).dt}{t^{n+1}.\bar c_{t^{-1}}(E)}\right)
$$
where $\res_t$ stands for the residue of the indicated Laurent
power-series.\footnote{The formula inside $\res_t$
 is a Laurent power series because the Chern classes of $E$ 
 are nilpotent (\ref{prop:FGL}).}
\end{lm}
Let us write
 $\psi(t)=\sum_i \alpha_i.t^i$ and
  ${\bar c_t(E)}^{-1}=\sum_{i\geq 0} \bar c_i(-E).t^i$.
Recall that for any element $\alpha \in \E^{**}(P)$,
 $\bar p_*(\alpha)$ is the coefficient of $\bar d^n$ in the
 decomposition of $\alpha$ within the $A$-basis $(\bar d^i)_{0 \leq i \leq n}$
 of $\E^{**}(P)$.
 A classical computation according to the defining relation of Chern classes
 \eqref{eq:Chern} in the additive case gives us:
$$
\bar p_*(\bar d^i)=\bar c_{i-n}(-E).
$$
In particular, one obtains:
$$
\bar p_*\left(\sum_i \alpha_i.\bar d^i\right)=\sum_i \alpha_i.\bar c_{i-n}(-E).
$$
To end the proof of the lemma, one has only to realize that the right hand side
 is the coefficient of $t^n$ in the following Laurent power-series:
$$
\left(\sum_i \alpha_i.t^i\right).c_{t^{-1}}(-E).
$$

Let us now compute the right hand side of \eqref{eq:thm:Quillen}.
 From the exact sequence of vector bundles over $P$:
$$
0 \rightarrow \mathcal O_{P} \rightarrow \lambda \otimes p^{-1}(E)
 \rightarrow T_p \rightarrow 0
$$
we get $\td(T_p)=\td(\lambda \otimes p^{-1}(E))$.
By assumption,
 the classes $\bar l_i=\log_F(l_i)$ are the Chern roots of $p^{-1}(E)$
 with respect to the additive orientation $\bar c$. Thus, by
  definition of the Todd class
  (see also the formula of \ref{num:basic_change} with $\Psi_S=\exp_F$),
  one gets:
$$
\td(\lambda \otimes p^{-1}(E))
 =\prod_{i \in I} \frac t {\exp_F(t)}\left(t=\bar d+\bar l_i\right).
$$
Note that because $\bar l_i$ are Chern roots of $E$ with respect to $\bar c$,
 one has:
$$
\prod_{i \in I} (\bar d+\bar l_i)
 =\sum_i \bar c_{n+1-i}(E).\bar d^{i}
 =\left.\big(t^{n+1}.c_{t^{-1}}(E)\big)\right|_{t=\bar d}.
$$
Because $F(x,y)=\exp_F(\log_F(x)+\log_F(y))$,
 one also has:
$$
\exp_F(\bar d+\bar l_i)=F(d,l_i)=F\big(\exp_F(\bar d),l_i\big).
$$
Thus to compute the right hand side of \eqref{eq:thm:Quillen}
 using the formula of the preceding lemma, we are led to
 introduce the following power series:
$$
\psi(t)
=\frac{\phi(\exp_F(t)).t^{n+1}.c_{t^{-1}}(E)}{\prod_{i \in I} F\big(\exp_F(t),l_i\big)}.
$$
Applying the preceding lemma:
$$
\bar p_*\big( \phi(d).\td(T_p) \big)=
 \bar p_*(\psi(\bar d))=\res_t
  \left(\frac{\phi(\exp_F(t)).dt}{\prod_{i \in I} F\big(\exp_F(t),l_i\big)}\right).
$$
Computing this residue with the change of variables $x=\exp_F(t)$,
 one gets the desired result:
 by \eqref{eq:formula_log_F}, $d\log_F(x)=\omega_F(x)$.
\end{proof}

\begin{rem}\label{rem:Quillen_formula}
\begin{enumerate}
\item In view of Remark \ref{ex:RR_source_MGL},
 the preceding formula is universal: It is valid  without any change
 for any $S$-absolute oriented ring spectrum $\E$ with associated
 formal group law $F$. Anyway, it is clear in the above proof that
 one can faithfully replace $\MGL$ by $\E$.
\item The above proof is particularly simple but it works only
 with rational coefficients whereas Quillen formula is stated in \cite{Qui}
 for complex cobordism with integral coefficients.
 However, one can at least deduce from this proof the formula
  with integral coefficients in characteristic $0$: indeed, in this case,
  one reduces to $S=\spec{\QQ}$ and we know from \cite{Lev} that
  $A=\MGL^{**}(\QQ)$ is the Lazard ring and thus has no torsion.
 One gets back the usual Quillen formula by using
  the complex realization functor.
\item As a particular case of the preceding formula, one gets
 the classical Myschenko formula computing the cobordism
 class of $\PP^n_S$ for any integer $n$: take $E=\AA^{n+1}_S$
 and $P(t)=1$. May be it is worth to summarize the proof in this case:
\begin{align*}
p_*(1)&=\bar p_*(\td(T_{\PP^n}))=\bar p_*(\td(\lambda \otimes \AA^{n+1})) \\
&=\bar p_*\left(\left(\frac{\bar d}{\exp{\bar d}}\right)^{n+1}\right)
=\res_t\left(\frac {dt} {\exp_F(t)^{n+1}}\right)
=\res_x\left(\frac {\omega_F(x)} {x^{n+1}}\right).
\end{align*}
Note also that this computation,
 including the change of variable $x=\exp_F(t)$, 
 was used by Borel and Serre in \cite{BS} to prove
 the classical Grothendieck-Riemann-Roch formula
 for $p:\PP^n_S \rightarrow S$
 -- replacing $\MGL$ by $\KGL$: as explained below (\ref{num:chern_character}),
  the Chern  character corresponds to changing the natural orientation on $\KGL$
  to the additive orientation.
\item The general formula of Quillen has been proved integrally
 in the context of oriented cohomology theories in \cite{Shi}
 and  \cite [App. B]{Vis}. The proofs given here are equally valid
 in our more general context.
\end{enumerate}
\end{rem}

\begin{ex}\label{ex:Quillen_cobord_formula}
An interesting particular case of Quillen formula is the following
 computation of the cobordism class of a projective bundle $P/S$
 associated with a vector bundle $E/S$:
$$
[P/S]=\res_t\left(\frac{\omega_F(t)}{\prod_i F(t,l_i)}\right)
$$
where $l_i$ are the Chern roots of $E$.
\end{ex}

\begin{ex}
We end-up this series of illustration of the (generalized) Riemann-Roch
 formula by explaining how one can also compare the Chern classes 
 arising from different orientations.

Let us consider again the general setting of Paragraph \ref{num:basic_change}.
In general, given a vector bundle $E/S$ of rank $n$,
 we can use the splitting principle to compute 
 the two different type of Chern classes of $E$
 associated respectively with the orientations $d$ and $c$.
 For any integer $0 \leq r \leq n$, one gets:
$$
d_r(E)=d_r(\oplus_{i=1}^n L_i)
=\sum_{\alpha} d_1(L_{\alpha_1})...d_1(L_{\alpha_r})
=\sum_{\alpha} \Phi(c_1(L_{\alpha_1}))...\Phi(c_1(L_{\alpha_r})),
$$
where $\alpha$ runs over the $r$-uple of integers $(\alpha_1,\hdots,\alpha_r)$
 such that
$$
1 \leq \alpha_1 < \alpha_2 < \hdots < \alpha_r \leq n.
$$
Once again, one has to express the right hand side in terms
 of the elementary symmetric functions in the $c_1(L_i)$ to obtain
 an expression in terms of $c_i(E)$.

In general, the formulas are pretty complicated.
However, one can remark there is a simple formula when $r=n$. Then,
 one simply obtains:
$$
d_n(E)=\td(E).c_n(E),
$$
generalizing the case of a line bundle (Par. \ref{num:basic_change}).

Another case that can be computed is the case $r=1$.
 Assume the strict isomorphism of formal group law $\Phi$
 such that $d=\Phi(c)$ has been written as:
$$
\Phi(x)=x+\sum_{i>1} \alpha_i.x^i.
$$
Recall the following classical relation (determinantal form
 of the Newton's identity) between the power sum
 symmetric polynomials $p_i$ and the elementary symmetric polynomials
 $e_j$:
$$
p_i=\left|
\raisebox{0.5\depth}{
\xymatrix@=0.1ex{
e_1 & 1\ar@{-}[rrrddd] & 0\ar@{.}[rr]\ar@{.}[rrdd] && 0\ar@{.}[dd]\\
2e_2\ar@{.}[ddd] & e_1\ar@{-}[rrrddd] & & & \\
& e_2\ar@{.}[dd]\ar@{-}[rrdd] & & & 0 \\
& & & & 1 \\
i.e_i & e_{i-1}\ar@{.}[rr] & & e_2 &  e_1 \\
}
}
\right|
$$
Using the splitting principle together with this relation, one obtains
 the following relation:
$$
d_1(E)=c_1(E)+\sum_{i>1} \alpha_i.
\left|
\raisebox{0.5\depth}{
\xymatrix@=0.1ex{
c_1(E) & 1\ar@{-}[rrrddd] & 0\ar@{.}[rr]\ar@{.}[rrdd] && 0\ar@{.}[dd]\\
2c_2(E)\ar@{.}[ddd] & c_1(E)\ar@{-}[rrrddd] & & & \\
& c_2(E)\ar@{.}[dd]\ar@{-}[rrdd] & & & 0 \\
& & & & 1 \\
i.c_i(E) & c_{i-1}(E)\ar@{.}[rr] & & c_2(E) &  c_1(E) \\
}
}
\right|.
$$
Note the sum above is finite because Chern classes are always
 nilpotent (cf. Prop. \ref{prop:FGL}(1)).
\end{ex}

\begin{ex}
As a last illustration of the change of orientation principle,
 let us consider an arbitrary oriented ring spectrum $(\E,c)$.
 One defines a new orientation on $E$ by putting for any
 line bundle $L$ over a scheme $S$:
$$
c'_1(L)=-c_1(L^\vee),
$$
where $L^\vee$ denotes the dual of $L/S$.
Indeed, if one denotes by $\mu(t)$ the formal inverse associated with
 the formal group law of $(\E,c)$, then the strict isomorphism
 corresponding to this change of orientation is the composition inverse
 of the power series $-\mu(t)$.

Using the splitting principle, one easily checks that for any vector bundle
 $E/S$, the following relation holds:
\begin{equation} \label{eq:each_Chern_dual}
c'_i(E)=(-1)^i.c_i(E^\vee).
\end{equation}
Let $\lambda$ be the canonical line bundle on $\PP(E)$.
 Then, from this formula and the relation \eqref{eq:Chern}
 defining Chern classes with respect
 to the orientation $c'$, one gets the following formula\footnote{This
  formula was suggested to me by Alberto Navarro.}
  expressing Chern
 classes in terms of the dual canonical line bundle $\lambda^\vee$:
\begin{equation} \label{eq:Chern_dual}
\sum_{i=0}^n (-1)^i p^*\big(c_i(E^\vee)\big).c_1\big(\lambda^\vee\big)^{n-i}=0.
\end{equation}
\end{ex}

\subsection{Universal formulas and the Chern character}

In the next proposition,
 we work over a fixed base scheme $S$.
Recall we have seen in Example \ref{ex:uniqueness_orientation}
 that there exists a unique orientation $c^{\HB}$ on the Beilinson
 motivic cohomology ring spectrum $\HBx S$ whose formal
 group law is necessarily additive. Using a fundamental result
 of \cite{CD2}, we get slightly more:
\begin{prop}\label{prop:HB_univ}
Let $(\E,c)$ be an additive oriented ring spectrum with rational coefficients
 over $S$.

Then there exists a unique morphism of absolute ring spectra
$$
\sigma:\HBx S \rightarrow \E
$$
which is moreover a morphism of oriented ring spectra:
 $\sigma(c^{\HB})=c$.

In particular, the morphism $\sigma$ is compatible with Gysin morphisms,
 residues, Chern classes and fundamental classes as constructed in the
 preceding sections.
\end{prop}
\begin{proof}
The existence and the uniqueness of $\sigma$ follow from
 \cite[14.2.16]{CD3}. The fact $\sigma$ is necessarily
 a morphism of oriented ring spectra follows from
 Proposition \ref{prop:chg_orient&iso_fgl}
 and the fact that the unique strict automorphism of
 a formal group law on a rational ring is the identity
 (fact recalled in Paragraph \ref{num:log&exp}).
The last assertion is then an application
 of the Riemann-Roch formula (or its extension for Chern classes
 and fundamental classes) to $\sigma$.
\end{proof}

\begin{rem}
Another way of stating the previous proposition
 is that $\HB$ is the universal absolute
 \emph{orientable} rational ring spectrum,
 -- whereas $\MGL$ is the universal absolute \emph{oriented}
  ring spectrum.

Thus, the proposition (as well as its integral counterpart stated in a few
 paragraphs) answers a \emph{desideratum} raised in the introduction 
 (cf. second paragraph) of Beilinson's fundamental work \cite{Bei}.
 
The morphism $\sigma$ could be called the \emph{higher cycle class} morphism:
 in degree $(2n,n)$, it gives a morphism from the Chow group of $n$-codimensional
 cycles.\footnote{Note that,
 even when $\E$ is the spectrum
 representing Deligne-Beilinson cohomology (say over $\QQ$, see \cite{HS}), the map
 $\sigma$ is not exactly the \emph{regulator map} as described
 in \cite[\textsection 3.3]{Sou2}.}
 Over a field and in the integral case (see below),
 this fits well with the fact motivic cohomology (of regular schemes)
 is given by Bloch's higher Chow groups.
\end{rem}

\begin{num} \label{num:chern_character}
Recall that the decomposition of algebraic K-theory
 according to the eigenvalues of the Adams operation has
 been lifted by J.~Riou to the stable homotopy category
 of schemes resulting in a canonical isomorphism of ring spectra:
$$
\ch_t:\KGL_\QQ \rightarrow \oplus_{i \in \ZZ} \HB(i)[2i].
$$
This is essentially the results of \cite{Riou}
 as explained in \cite{CD3}, Lemma 14.1.4 for the existence of this isomorphism
 and Corollary 14.2.17 for the fact it is a morphism of ring spectra.
By definition, for any regular scheme $X$ and any integer $n$,
 this isomorphism induces the canonical decomposition
$$
K_n(X)_\QQ \rightarrow \oplus_{i \geq 0} \mathrm{Gr}^i_\gamma K_n(X)_\QQ
$$
for the $\gamma$-filtration on Quillen K-theory with rational coefficients
 (cf. \cite{Sou}). According to \cite[\textsection 7]{Sou},
 the projector on the $i$-th graded
 part of this decomposition is given by the $i$-th Chern character.
 Thus the morphism of ring spectra $\ch_t$ lifts 
 to the category of spectra
 the usual Chern character
 in higher $K$-theory with values in Beilinson motivic
 cohomology (see again \cite{Sou}):
 for any pair $(n,r) \in \NN^2$ and any regular scheme,
 one gets the usual (higher) Chern character:
$$
\ch_{r,n}:K_r(X) \rightarrow \HB^{2n-r,n}(X).
$$
For $r=0$, this coincide withs the Chern character of \cite{SGA6}.
In particular, it is uniquely characterized by its value
 on the class of a line bundle $L/X$:
$$
\ch_{0,n}\big([L]\big)=\frac 1 {n!} c_1(L)^n.
$$

Using the principle explained in paragraphs \ref{num:comput_todd1}
 and \ref{num:basic_change}, we can easily determine
 the Todd class associated with the morphism of spectra $\ch_t$.
 Indeed, the formal group law associated with $c^\KGL$ is
 $F_\KGL(x,y)=x+y-\beta.xy$ (Ex. \ref{ex:abs_ring_spectra_FGL}).
 By definition of $\beta$, we get that $\ch_t(\beta)=1$.
 We deduce that the formal group law
 associated with the orientation $\ch_t(c^\KGL)$ on the aim
 is the following multiplicative formal group law $F_{mult}(x,y)=x+y-xy$.
 On the other hand $c^\HB$ is additive,
 thus the pair of reciprocal strict isomorphisms associated
 with $\ch_t$ is the logarithm/exponential of the formal group
 law $F_{mult}$.
 Formula \eqref{eq:formula_log_F} easily yields:
\begin{align*}
& \Phi(t)=\log_{F_{mult}}(t)=-\ln(1-t),\\
& \Psi(t)=\exp_{F_{mult}}(t)=1-e^{-t}.
\end{align*}
Thus, from \ref{num:basic_change},
 we obtain that the Todd class associated with $\ch_t$
 is defined by the power series $\frac t {1-e^{-t}}$
  (Prop. \ref{prop:todd}).
 In particular, it coincides exactly with the usual Todd class
$$
\td:K_0 \rightarrow  \mathrm{Gr}^*_\gamma K_0(X)_\QQ
 \simeq \HB^{2*,*}(X)
$$
 of \cite{SGA6}.

In the end, taking care about the isomorphisms
 \eqref{eq:iso_KGL&Ktheory}, Theorem \ref{thm:projective_RR}
 yields the following higher arithmetic Grothendieck-Riemann-Roch formula:
\end{num}
\begin{prop} \label{prop:classical_GRR}
Consider the preceding notations.
Let $f:Y \rightarrow X$ be a projective morphism
 between regular schemes and $\tau_f \in K_0(Y)$ be its
 virtual tangent bundle (see \ref{num:virtual_tangent}).

Then for any integer $r \geq 0$,
 the following diagram is commutative:
$$
\xymatrix@R=18pt@C=24pt{
K_r(Y)_\QQ\ar_{\td(\tau_f).\ch_t}[d]\ar^{f_*}[r]
 & K_r(X)_\QQ\ar^{\ch_t}[d] \\
\oplus_{n \geq 0}\HB^{2n-r,n}(Y)\ar^{f_*}[r]
 & \oplus_{n\geq 0} \HB^{2n-r,n}(X). \\
}
$$
In other words, for any integer $n \geq 0$
 and any element $y \in K_r(Y)$, one has:
$$
\ch_{r,n}\big(f_*(y)\big)
 =f_*\left( \sum_{i+j=n} \td_i(\tau_f).\ch_{r,j}(y)\right). 
$$
\end{prop}

\begin{rem}
Taking into account Proposition \ref{prop:HB_univ},
 the preceding formula is universal.
 Indeed, given any additive oriented rational ring spectrum $\E$,
  the Chern character and the Todd class with values in $\E$
  are induced by those with values in Beilinson motivic cohomology.
 Then we get the following commutative diagram:
$$
\xymatrix@R=18pt@C=58pt{
K_r(Y)_\QQ\ar^-{\td(\tau_f).\ch_t}[r]\ar_{f_*}[d]
 & \oplus_{n \geq 0}\HB^{2n-r,n}(Y)\ar^{f_*}[d]\ar^{\sum \sigma^{2n-r,n}_Y}[r] 
 & \oplus_{n \geq 0}\E^{2n-r,n}(Y)\ar^{f_*}[d] \\
K_r(X)_\QQ\ar^-{\ch_t}[r]
 & \oplus_{n \geq 0} \HB^{2n-r,n}(X)\ar^{\sum \sigma^{2n-r,n}_X}[r] 
 & \oplus_{n \geq 0}\E^{2n-r,n}(X),
}
$$
\end{rem}

Note that the map $\ch_t$ exists over any scheme.
On the other hand, it is known from \cite{Cis} that
 the spectrum $\KGL_\QQ$ represents Weibel homotopy invariant
 K-theory $KH_*$. Therefore, because $\KGL_\QQ$
 and $\HB$ are $\sm_S$-absolutely pure for any scheme $S$,
 the preceding proposition admits the following version,
 which is finer when $S$ is singular:
\begin{prop} \label{prop:H_classical_GRR}
Let $S$ be any scheme and $f:Y \rightarrow X$ be a projective morphism
 between smooth $S$-schemes with virtual tangent bundle $\tau_f \in K_0(Y)$.

Then for any integer $r \geq 0$,
 the following diagram is commutative:
$$
\xymatrix@R=18pt@C=24pt{
KH_r(Y)_\QQ\ar_{\td(\tau_f).\ch_t}[d]\ar^{f_*}[r]
 & KH_r(X)_\QQ\ar^{\ch_t}[d] \\
\oplus_{n \geq 0}\HB^{2n-r,n}(Y)\ar^{f_*}[r]
 & \oplus_{n\geq 0} \HB^{2n-r,n}(X). \\
}
$$
\end{prop}

\begin{rem}
The morphism $f_*$ in the above proposition stands \emph{a priori}
 for the Gysin
 morphisms we have defined. However,
 recall from \cite[3.16.4]{TT}, that K-theory is covariant
 with respect to proper maps of finite Tor-dimension,
 thus a fortiori with respect to the morphism $f$ in the above
 proposition. The definition of $KH_*$ makes it clear that
 this covariant functoriality extends to $KH_*$ in such
 a way that the following diagram commutes:
$$
\xymatrix@R=10pt@C=16pt{
K_r(Y)\ar[d]\ar^{f_*}[r]
 & K_r(X)\ar[d] \\
KH_r(Y)\ar^{f_*}[r]
 & KH_r(X).
}
$$
Applying Theorem~\ref{thm:compare_Gysin} as in
 Example~\ref{ex:comput_Gysin}(3), we get that this last
 morphism coincides with our Gysin morphism.
 Therefore, because of the preceding commutative square,
 the above Riemann-Roch formula is even true when
 replacing $KH_*$ by Thomason-Trobaugh algebraic
 K-theory.\footnote{Recall it coincides with Quillen's definition
 whenever $X$ admits an ample line bundle.}
\end{rem}

\begin{num}
Finally, we show how to get an almost integral (weaker) version
 of Proposition \ref{prop:HB_univ}
 using the recent Hopkins-Morel-Hoyois theorem (cf. \cite{Hoy}).
 Let us restate this theorem for the sake of notations.
Consider a field $k$ of exponential characteristic $p$.
The canonical orientation on $\MGL$ induces a formal group law
 on its ring of coefficients, say over $k$.
 By universality of the Lazard ring, we get a canonical morphism of ring:
$$
L \longrightarrow \MGL_{**}:=\MGL^{-*,-*}(\spec k)
$$
which in fact induces the following isomorphism (cf. \cite[8.2]{Hoy})
 of graded rings:
\begin{equation}\label{eq:MGL&Lazard}
L[1/p] \longrightarrow \MGL_{(2,1)*}[1/p]
\end{equation}
where $\MGL_{(2,1)*}$ denotes the elements of $\MGL_{**}$ of
 bidegree $(2n,n)$ for an integer $n \in \ZZ$.
 The Lazard ring is in fact a graded polynomial algebra
 $L=\ZZ[b_1,b_2,\hdots]$ where $b_i$ has degree $i$.
 In particular, any element $b_i$ corresponds to an element of $\MGL_{2i,i}$.
 Following \cite{Hoy}, paragraph after Corollary 6.9, one can define
 the quotient spectrum
$$
\MGL/\{b_i, i > 0\}
$$
by killing the elements $b_i$.
 In fact, this amounts to kill the coefficients $a_{ij}$ of the formal group
 law of $\MGL$ --- which, after inverting $p$, corresponds to the universal
 formal group law on the Lazard ring through the isomorphism \eqref{eq:MGL&Lazard}.

As the formal group law associated with the canonical orientation
 on the motivic Eilenberg Mac Lane spectrum $\HH_{\ZZ}$ is additive,
 one gets a canonical map of spectra:
\begin{equation}\label{eq:HMH_iso}
\big(\MGL/\{b_i, i>0\}\big) \rightarrow \HH_{\ZZ}
\end{equation}
which becomes an isomorphism after inverting $p$ (cf. \cite[7.12]{Hoy}).
As a corollary, we get:
\end{num}
\begin{prop}
Let $k$ be a field of exponential characteristic $p$
 and $\E$ be a $\ZZ[1/p]$-linear $k$-absolute ring spectrum.

Then the following conditions are equivalent:
\begin{itemize}
\item $\E$ admits an additive orientation $c$.
\item There exists a morphism of absolute ring $k$-spectra:
$$
\sigma:\HH_\ZZ \rightarrow \E.
$$
\end{itemize}
Moreover, when these conditions are fulfilled,
 the additive orientation $c$ on $\E$ is unique
 and is the image under $\sigma$ of the canonical
 orientation on $\HH_\ZZ$.
 In particular, $\sigma$ is compatible with Gysin morphisms,
 residues, Chern classes and fundamental classes as constructed 
 in the preceding sections.
\end{prop}
\begin{proof}
The fact (i) implies (ii) is obvious: the orientation
 of $\HH_\ZZ$ induces an orientation of $\E$, which is
 necessarily additive.

Reciprocally. According to
 Prop. \ref{prop:orientation&MGL-modules}, the orientation $c$ of $\E$
 corresponds to a morphism of ring spectra over $\spec k$
$$
\varphi:\MGL \rightarrow \E.
$$
This map induces a morphism of formal group law.
Moreover, according to the commutative diagram:
$$
\xymatrix@C=20pt@R=10pt{
\MGL^{**}(\PP^\infty)\ar[r]\ar_{\varphi_{*}}[d]
 & \MGL^{**}(\PP^\infty \times \PP^\infty)\ar^{\varphi_{*}}[d] \\
\E^{**}(\PP^\infty)\ar[r]
 & \E^{**}(\PP^\infty \times \PP^\infty),
}
$$
the definition of the formal group law associated with
 an oriented ring spectrum and the fact $(\E,c)$ is additive,
 we obtain that
 $\varphi_{*}:\MGL_{**}(k) \rightarrow \E_{**}(k)$
 sends all the elements $a_{ij}$ to $0$
 as soon as $(i,j) \neq (1,0), (0,1)$.
 Thus, $\varphi$ induces the morphism of spectra
 $\sigma$ using the isomorphism induced by \eqref{eq:HMH_iso}
 and the assumption that $\E$ is $\ZZ[1/p]$-linear.

Having fixed an additive orientation $c$ of $\E$,
 from \cite[Th. 4.3]{Vez}, the orientations $d$ on $\E$ are in one-to-one
 correspondence with the isomorphisms of formal group laws on $\E_{(2,1)*}(k)$
 from the additive formal group to the formal group law $F_d$ associated
 with $d$. As there is only one automorphism of the additive formal group
 law, the uniqueness statement follows.
 The remaining assertion is clear from what we have seen.
\end{proof}

\begin{rem}
Another way of stating the preceding proposition is that
 given a representable $\ZZ[|1/p]$-linear
 cohomology theory $\E^{**}(-)$ over (smooth) $k$-schemes,
 there is only one possible way to define Chern classes which are additive. 
 Moreover the other structures on the cohomology, that is
 residues, Gysin morphisms and fundamental classes, are unique.
 Similarly, the cycle class map, from classical Chow groups, is unique.
\end{rem}

\subsection{Residues and symbols}

\begin{num}\label{num:res&symbols}
Let us consider one of the following two cases:
\begin{itemize}
\item $S=\spec{\ZZ}$, $\Lambda=\QQ$;
\item $S=\spec{k}$, $k$ a field of exponential characteristic $p$,
 $\Lambda=\ZZ[1/p]$.
\end{itemize}
Let us denote simply by $\HH_\Lambda$ either the Beilinson motivic cohomology
 spectrum in the first case or the Eilenberg-Mac Lane motivic ring
 spectrum with coefficients in $\Lambda$ in the second case.

According to the preceding section,
 any $S$-absolute oriented ring spectrum $\E$ with additive formal
 group law has a structure of
 $\HH_\Lambda$-algebra with structural morphism $\sigma$.

In particular, given any field $K$ over $S$,
 we get a canonical symbol map:
$$
K_n^M(K) \otimes_\ZZ \Lambda
 \simeq \HH_\Lambda^{n,n}(K) \xrightarrow{\sigma_*}
  \E^{n,n}(K), \{f_1,\hdots,f_n\} \mapsto \{f_1,\hdots,f_n\}_\E.
$$
If $K$ admits a discrete valuation $v$ such that its ring of integers 
 is an $S$-scheme, with residue field $k$, we have obtained
 (Ex. \ref{ex:tame_residue})
 a residue map with coefficients in $\E$: 
$$
\partial_v^\E:\E^{n,m}(K) \rightarrow  \E^{n-1,m-1}(k).
$$
As an easy application of the Riemann-Roch residual theorem,
 we get the following computation of this residue:
 given any elements $f_1,...,f_n$ in $K^\times$,
$$
\partial_v^\E(\{f_1,\hdots,f_n\}_\E)
=\sigma_*\big(\partial_v(\{f_1,\hdots,f_n\})\big)
$$
where $\partial_v$ is Milnor tame residue symbols
 (cf Example \ref{ex:tame_residue}).
In fact, this follows from Theorem \ref{thm:RR1}
 applied to $\sigma$: the Todd class is equal to $1$ in this case as
 $\sigma$ is a morphism of oriented ring spectra.\footnote{One could
 also obtained this computation directly as in Example
 \ref{ex:tame_residue} using Proposition \ref{prop:specialisation}.}
\end{num}

\begin{ex}\label{ex:residues_Weil}
\begin{enumerate}
\item {\it De Rham cohomology}.-- Assume $S=\spec k$ where $k$ is a field of characteristic $0$
 and consider $\E_{dR}$ the $k$-absolute ring spectrum
 representing De Rham cohomology
  (cf. \ref{ex:abs_ring_spectra}(1)).
 Recall the twist on that cohomology is just given
 by the tensor product with the $1$-dimensional $k$-vector space
 $k(1):=H^1_{dR}(\GG)^\vee$.

Note that $H^1_{dR}(\GG)=k.d\log(t)$,
 where $\GG=\spec{k[t,t^{-1}]}$.
 The choice of the generator $d\log(t)$ determines
 an isomorphism:
$$
\E^{n,i}_{dR}(X) \simeq H^n_{dR}(X/k)
$$
functorial in any smooth $k$-scheme $X$.
As already mentioned, the fact $\E_{dR}$ is a $k$-absolute
 ring spectrum extends De Rham cohomology to any $k$-scheme.
 In the particular case of an extension field $K/k$,
 the choice of $d\log(t)$ gives a canonical isomorphism:
$$
\E^{n,i}_{dR}(K) \simeq H^n_{dR}(K/k).
$$
Through this isomorphism (case $i=n$),
 the \emph{de Rham symbol} associated with a family of units
 $f_1,...,f_n$ in an extension field $K/k$ is given by the classical formula:
$$
\{f_1,\hdots,f_n\}_{dR}=d\log(f_1) \wedge \hdots \wedge d\log(f_n)
 \in H^n_{dR}(K/k).
$$
Given now a discrete valuation ring $(K,v)$ over $k$,
 according to
 the previous paragraph one gets a residue map:
$$
\partial_v^{dR}:H^1_{dR}(K/k) \simeq \E^{1,1}_{dR}(K) \rightarrow
 \E^{0,0}_{dR}(\kappa(v))\simeq k.
$$
According to the previous computation of residues on symbols,
 one gets:
\begin{equation}
\partial^{dR}_v(d\log(f))=v(f). 
\end{equation}
\item {\it Rigid cohomology}.-- 
Let $V$ be a complete discrete valuation ring with fraction field $E$
 and residue field $k$. Let $\E_{rig}$
 be the $k$-absolute ring spectrum representing rigid cohomology $H_{rig}(-/E)$.

The situation is analogous to the previous one though the constructions are
 less concrete. By definition, $H^1_{rig}(\GG/E)$ is the rational part of 
 the first cohomology
 group of the weakly complete De Rham complex associated with
 the weakly complete $V$-algebra $V\{t,t^{-1}\}$. In particular,
 it is generated by the differential form $d\log(t)$ of $V\{t,t^{-1}\}$.
 For smooth (or even singular) $k$-scheme $X$,
  the choice of this differential gives a canonical isomorphism
$$
\E^{n,i}_{rig}(X) \simeq H^n_{rig}(X/E)
$$
Given a unit $f:X \rightarrow \GG$,
 this isomorphism in the case $n=i=1$ sends the symbol $\{f\}_{rig}$
 to the element $d\log(f):=f^*(d\log(t))$.
 Similarly, $\{f_1,\hdots,f_n\}_{rig}$ corresponds
 to $d\log(f_1) \wedge \hdots \wedge d\log(f_n)$.

Given an extension field $K/k$, one gets:
$$
\E^{n,i}_{rig}(K)=\ilim{A \subset K} H^n_{rig}(\spec A/E)=:H^n_{rig}(K/E)
$$
where $A$ runs over the sub-rings of $K$ which are smooth
 of finite type over the inseparable closure of $k$ in $K$.
 And for any discrete valuation $v$ on $K$,
 we get a tame residue symbol:
$$
\partial_v^{rig}:H^1_{rig}(K/E)
 \rightarrow H^0_{rig}(K/E)\simeq E
$$
satisfying the expected property on symbols.
\end{enumerate}
\end{ex}

\begin{rem}
Symbols in differential calculus have a beautiful history
 starting from van der Kallen formula (\cite{vdK}) and going to
 the Bloch-Kato conjecture modulo $p$ (\cite[2.1]{BK}). 
\end{rem}

\begin{num}
Consider again the notations of point (1) of the preceding example.
Let $C$ be a proper connected regular curve over $k$
 with function field $K$. Let us fix a closed point
 of $C$, in other words a discrete valuation $v$ on $K$
 trivial on $k$.
 In \cite{Tate}, Tate gives a purely algebraic definition of
 a residue map: $\partial^{Tate}_v:H^1_{dR}(K/k) \rightarrow k$.
 In view of the preceding example, the reader should not be
 surprised by the following comparison result:
\end{num}
\begin{prop} \label{prop:compare_res_Tate}
Using the above notations, $\partial_v^{dR}=\partial_v^{Tate}$.
\end{prop}
\begin{proof}
Both definitions are invariant under completion with respect to 
 the valuation $v$
 (see Prop. \ref{prop:inv_completion} for $\partial_v^{dR}$).
 Thus we can replace $K$ by $\hat K \simeq k((t))$ and
  we are reduced to identify the two residues on
  differential forms of the form $\omega=f(t)dt$ when 
  $f(t)$ is a power series with coefficients in $k$.
  By continuity and additivity of residues, we can assume
  $\omega=t^i.dt$.

The case $i\geq 0$ is easy because then $\omega$ can be extended
 to the valuation ring of $k[[t]]$ and therefore its image by $\partial_v^{dR}$ 
 is $0$ according to the Gysin long exact sequence \eqref{eq:Gysin_exact_seq}.
 The case $i=-1$ follows from the residual Riemann-Roch formula
 as explained in the previous example.
 For the remaining case, we use the reciprocity formula for $\PP^1_k$
 and $\partial_v^{dR}$ (according to \cite[5.2.1]{Deg5} and \cite[2.2]{Ros}).
 According to this formula, we get:
$$
\sum_{x \in \PP^1_{k(0)}} \partial_x(\omega)=0.
$$
Because $\omega=t^i.dt$, the above equality gives:
$$
\partial_0^{dR}(t^i.dt)=-\partial_\infty^{dR}(t^i.dt)
\stackrel{(1)}=-\partial_\infty^{dR}(-u^{-(i+2)}.du)
\stackrel{(2)}= 0
$$
where equality (1) is given by the
 substitution $u=t^{-1}$, and (2) is true
 because $i+2 \leq 0$. This concludes.
\end{proof}

\begin{rem}
The above proof also gives a tool to compute residues in rigid cohomology.
Given any cohomology class $\omega$ in $H^1_{rig}(K/E)$,
 we can extend it to the completion of $K$ thus it corresponds to an element
 in $H^1_{rig}(k((t))/E)$, in other words an overconverging
 differential $\hat \omega$ form over $V((t))$.
 Then $\partial_v(\omega)$ is the residue in the usual sense:
 write $\omega=f(t).dt$ with $f(t) \in V\{t\}$, then
$$
\partial_v(\omega)=\mathrm{res}_t(f(t))
$$
seen as an element of $K$. More generally,
 our construction of residues
 should be linked with that of \cite[VII, 1.2]{Ber}.
\end{rem}

\begin{ex} \label{ex:coniveau}
\textit{Functoriality of coniveau spectral sequences}.--
Assume that $\base$ is the category of excellent regular schemes.
Let $\E$ be any absolute oriented ring spectrum.

Using the method of Bloch-Ogus in \cite{BO}
 and Gysin long
 exact sequences of the form \eqref{eq:Gysin_exact_seq},
 we get for any regular excellent scheme $X$ and any integer $n \in \ZZ$,
 a \emph{coniveau} spectral sequence of the form:
$$
E_1^{p,q}=\bigoplus_{x \in X^{(n)}} \E^{q-p,n-p}(\kappa(x))
 \Rightarrow \E^{p+q,n}(X)
$$
which converges to the coniveau filtration on $\E^{**}(X)$.

Recall that this spectral sequence can be defined 
 using the exact couple:
$$
D^{p,q}=\ilim{Z^*} \E^{p+q,n}(X-Z^{p+1}), 
 E^{p,q}=\ilim{Z^*} \E^{q-p,n-p}(Z^p-Z^{p+1})
$$
where the limit is taken over the sequences $(Z^p)_{p  \in \NN}$
 such that $Z^p$ is a closed subscheme of codimension $\geq p$ in $X$
 satisfying the condition that $(Z^p-Z^{p+1})$ is regular. 
 Indeed, using the fact $X$ is excellent, 
 we obtain that the corresponding set, ordered
 by term-wise inclusion, is filtering.

Then the maps of the exact couples are given by considering
 localisation long exact in sequences in cohomology --- which come
 from Gysin triangles:
$$
\xymatrix@R=20pt@C=40pt{
E^{p,q+1} \leftrightsquigarrow \E^{q-p+1,n-p}(Z^p-Z^{p+1})
& \E^{q-p,n-p}(Z^p-Z^{p+1}) \leftrightsquigarrow E^{p,q}
 \ar^{i_{p*}}[d] \\
D^{p,q} \leftrightsquigarrow \E^{p+q,n}(X-Z^{p+1})
  \ar^{\partial_{X-Z^{p+1},Z^p-Z^{p+1}}}[u]
& \E^{p+q,n}(X-Z^{p}) \leftrightsquigarrow D^{p-1,q+1}
 \ar^{j_p^*}[l]
}
$$
(see also the presentation of \cite{Deg13}).
According to this presentation, the fact that this spectral sequence,
 and especially the differentials in the $E_1$-term, is functorial
 with respect to any pseudo-morphism $\E \rightarrow \F$ of oriented ring spectra
 follows from the Riemann-Roch formula applied to the morphisms in the above diagram
 and the fact that it is enough to consider sequences $Z^*$
 as above and such that the normal bundle of $Z^p-Z^{p+1}$ in
 $X-Z^{p+1}$ is trivial.

Note that this phenomena was already observed in a particular case
 in \cite[Th. 3.9]{Gil}.  
\end{ex}

\subsection{Residual Riemann-Roch formula} \label{sec:application_residual}

\begin{num}
So far, we have only worked out the trivial form
 of the residual Riemann-Roch formula,
 when the Todd class involved is $1$. Let us express the general residual
 Riemann-Roch formula in the case of the usual Chern character,
 as introduced in Paragraph \ref{num:chern_character}.
\end{num}
\begin{thm}
Consider a closed regular pair $(X,Z)$ of codimension $c$.
Let $N_ZX$ be the normal bundle of $Z$ in $X$ and put $U=X-Z$.

Then the following diagram is commutative:
$$
\xymatrix@R=20pt@C=60pt{
K_r(U)\ar^{\partial_{X,Z}}[r]\ar_{\qquad\qquad\ch_{r,n}^U}[d]
& K_{r-1}(Z)\ar^{\sum_{i+j=n-c} \td_i(-N_ZX).\ch_{r-1,j}^Z}[d] \\
\HB^{2n-r,n}(U)\ar^-{\partial_{X,Z}}[r] & \HB^{2(n-c)-r+1,n-c}(Z)
}
$$
\end{thm}

\begin{rem}
Once again, using the universality of Beilinson motivic cohomology
 among the absolute oriented cohomology with additive formal group law
 (Prop. \ref{prop:HB_univ}),
 the preceding formula gives also a formula for the classical
 (mixed Weil) cohomologies.
\end{rem}

\begin{ex} \underline{Let us first consider the case $r=1$}.
Then one gets an explicit description of the residue morphism for K-theory,
 when $X=\spec A$ and $Z=\spec{A/I}$,
 assuming $A$ and $A/I$ are regular rings.

Indeed, one knows that $K_1(A_I)=\mathrm{GL}(A_I)^{ab}$, the abelianization
 of the group of invertible matrices of arbitrary dimensions.
 Assume we are given an endomorphism $u:A^r \rightarrow A^r$
  such that $u \otimes_A A_I$ is an automorphism of $A_I^r$.
 We will denote by $[u]$ the class of this isomorphism in
  $K_1(A)$.
 By assumption, $u$ is a monomorphism whose cokernel is supported on $I$.
 We denote by $[\coker(u)]$ the class of the corresponding 
 (finitely presented) $A/I$-module in $K_0(A/I)$.
 With these notations, one has the following formula:
$$
\partial_{X,Z}([u])=[\coker(u)].
$$

Assume furthermore $n=c=1$.
Recall that the following part of the higher Chern character
 $\ch_{1,1}:K_1(A) \rightarrow \HB^{1,1}(A)=A^\times \otimes \QQ$
 sends any matrix of $GL(A)$ to its determinant.

Assume $Z$ is connected. By assumption, $I=(\pi)$
 for a prime divisor $\pi$: $A_I$ is a discrete valuation ring.
 We let $v_\pi$ denote its valuation.
 Then, giving the above notations, the residual Riemann-Roch formula 
 lands in $\HB^{0,0}(Z)=\QQ$ and reads:
$$
\boxed{
v_\pi \big(\det(u)\big)=\mathrm{rk}_{A/I}\big([\coker(u)]\big).
}
$$
Note that in fact, it is an integral formula as all the members are integers.

\noindent \underline{Consider now the case $n>1$, $r$ arbitrary}:
According to the coniveau spectral sequence of \ref{ex:coniveau} applied for
 $\E=\HB$, we get that the group $\HB^{2n-1,n}(U)$ is the cohomology, tensored with $\QQ$,
 in the middle of the following complex:
$$
\bigoplus_{y \in U^{(n-2)}} K_2^M(\kappa(y))
 \longrightarrow \bigoplus_{x \in U^{(n-1)}} \kappa(x)^\times
 \xrightarrow{\ \mathrm{div}\ } Z^n(U)
$$
where the last map is the usual divisor class map (computation of Quillen).
Thus, any element $f$ of $\HB^{2n-1,n}(U)$ can be described as
 the class of a finite sum:
$$
\sum_{x \in U^{(n-1)}} f_{x}
$$
where $f_{x}$ is a unit of $\kappa(x)$, which is
 the identity for almost all $x$,
 and such that the following $n$-codimensional cycle of $U$ is zero:
$$
\sum_{x \in U^{(n-1)}} \mathrm{div}_U(f_{x})=0.
$$
Moreover, using this description of the group $\HB^{2n-1,n}(U)$,
 the residue map 
 $\partial_{X,Z}:\HB^{2n-1,n}(U) \rightarrow CH^{n-c}(Z)_\QQ$
 can be described as follows:
$$
\partial_{X,Z}(f)=\sum_{x \in U^{(n-1)}} \mathrm{div}_X(f_{x})
$$
where $\mathrm{div}_X$ denotes the divisor of the rational function of $f_x$
 seen as a cycle in $X$. Indeed, by assumption on $(f_x)$,
 this cycle has support in $Z$. 
 
Thus, we will represent the element $\ch_{1,n}([u]) \in \HB^{2n-1,n}(U)$ as
the class of a sum:
$
\sum_{x \in U^{(n-1)}} f_{u,x}
$
satisfying the conditions above.\footnote{It
is rather delicate to give a formula for the $f_{u,x}$.
 One can only say that they describe the part of the element
 $[u] \in K_1(U)$
 (topologically) supported in codimension $(n-1)$.} 

Assume $n=c$ and $Z$ is connected with generic point $\eta$.
Then the residual Riemann-Roch formula lands again in $\HB^{0,0}(Z)=\QQ$
 and reads:
$$
\boxed{
\sum_{x \in U^{(n-1)}, \eta < x} \mathrm{ord}^x_\eta(f_{u,x})
 =\mathrm{rk}_{A/I}\big([\coker(u)]\big)
}
$$
where $\mathrm{ord}_\eta^x$ denotes the order of the rational function at $x$.
Observe this is again an integral formula.

Let us finally consider the case where the codimension $c$ of $Z$ in $X$
 is arbitrary less than $n$. Then the residual formula can be stated as 
 the following equality of cycles in $CH^{n-c}(Z)_\QQ$:
$$
\boxed{
\sum_{x \in U^{(n-1)}} \mathrm{div}_X(f_{u,x})
 =\sum_{i+j=n-c} \td_i(-N_ZX).\ch_{0,j}\big([\coker(u)]\big).
}
$$
(Recall that by assumption, the cycle on the left has support in $Z$.)
\end{ex}

\section{The axiomatic of Panin revisited}
\label{sec:Panin}

\subsection{Axioms for (arithmetic) cohomologies}

Our theory is obviously linked with the more classical theory of
 \emph{oriented cohomology theory} developed by Panin
 and Smirnov (see \cite{Panin1,PaninRR,Panin2} and \cite{Smir1,Smir2}). 
 Our axioms are more restrictive
 as we ask for representability of a cohomology theory.
 Nevertheless, one can extract from Section \ref{sec:abs_coh&pur}
 the following generalization of the axioms used by Panin:
\begin{df}\label{df:coh_thy}
A \emph{ringed cohomology theory (with supports)} $\E$ on $\base$ is the datum 
 for each closed pair $(X,Z)$ in $\base$ of a bigraded
 abelian group $\E^{**}_Z(X)$ equipped with the following structures:
\begin{itemize}
\item contravariant functoriality as described in \ref{num:pullback_support},
\item covariant functoriality as described in \ref{num:pushforward_support},
\item refined products as described in \ref{num:products},
\item for each closed pair $(X,Z)$,
 a \emph{boundary morphism}
$$\E^{n,m}(X-Z) \xrightarrow{\delta_{X,Z}} \E^{n+1,m}_Z(X)$$
contravariantly natural and fitting in a Gysin long exact sequence
of the form \eqref{eq:loc_exact_seq},
\end{itemize}
which satisfies the axioms (E1)-(E7) described in
 Prop. \ref{prop:ppty_product} together with the following additional
 properties:
\begin{itemize}
\item \textit{Homotopy}: for any scheme $X$,
 $\E^{**}(X) \rightarrow \E^{**}(\AA^1_X)$ is an isomorphism,
\item \textit{Stability}: For any scheme $S$,
 let $\tilde \E^{2,1}(\PP^1_S):=\E^{2,1}(\PP^1_S)/\E^{2,1}(\{\infty\})$.
 There exists a family of classes $\eta_S \in \tilde \E^{2,1}(\PP^1_S)$
 indexed by schemes in $\base$ which is stable by pullbacks 
 and such that for any scheme $S$ and
 any integers $(n,m)$ the following map is an isomorphism:
$$
\E^{n,m}(S) \rightarrow \tilde \E^{n+2,m+1}(\PP^1_S), x \mapsto \eta_S.p^*(x).
$$
\item \textit{Excision}: for any excisive 
 morphism of closed pairs $f:(Y,T) \rightarrow (X,Z)$
 (see Def. \ref{df:excisive_morph}),
 the pullback $f^*:\E^{**}_Z(X) \rightarrow \E^{**}_T(Y)$ is an isomorphism.
\end{itemize}
A morphism of ringed cohomology theories with support is a natural
 transformation compatible with contravariant and covariant
 functorialities, with refined products and with the operator
 $\delta_{X,Z}$ for closed pairs $(X,Z)$ in $\base$.

We will say that $\E$ is \emph{oriented} if there exists a natural
 transformation of presheaves of sets on $\base$:
$$
c:\Pic \rightarrow \E^{2,1}
$$
such that for any scheme $S$, $c_{\PP^1_S}(\lambda)=\eta_S$
 where $\lambda=\mathcal O(-1)$ is the canonical line bundle on $\PP^1_S$.
 
We will say that a closed pair $(X,Z)$ is $\E$-pure
 if the morphisms
$$
\E^{**}(X,Z) \xleftarrow{\sigma_1^*} \E^{**}(D_ZX,\AA^1_Z)
 \xrightarrow{\sigma_0^*} \E^{**}(N_Z,Z)
$$
induced by the deformation diagram \eqref{eq:deformation}
 are isomorphisms. We say that $\E$ is \emph{absolutely pure}
 if any regular closed pair in $\base$ is $\E$-pure.
 
For short, we will say \emph{arithmetic cohomology}
 for a ringed cohomology with support which is oriented
 and absolutely pure.
 A morphism of arithmetic cohomology is defined likewise,
 but beware we do not require the compatibility with the given
 orientations.
\end{df}

\begin{num} \label{num:extend}
With this definition, we can extend all the results of sections
 \ref{sec:orientation}, \ref{sec:Gysin} and \ref{sec:RR} as follows:
\begin{enumerate}
\item One has to pay a special attention to the projective bundle formula
 (\ref{thm:PBF}) and realize that the lemma of Morel \ref{lm:Morel}
 can be stated and proven using cohomology with support (instead of
 working in the unstable homotopy category).
 Then one gets the working theory of Chern classes
 and formal group laws as established
 in Section \ref{sec:orientation&Chern}.
\item  For sections \ref{sec:fdl_class}, \ref{sec:inter_theory},
 \ref{sec:Gysin} and \ref{sec:RR}, the arguments just go through
 as we have been careful to rely only on the axiomatic described
 in section 1 and restated in the previous definition.
\end{enumerate}
The results obtained here cover the one proved earlier by Panin.
\end{num}

\begin{rem}
Note that the axiomatic described here differs especially 
 with that of Panin because of two points:
\begin{itemize}
\item we have devised another axiom, \emph{absolute purity},
 especially relevant in the arithmetic case;
\item we asked for the existence of a refined product.
\end{itemize}
Both properties are very strong and the natural examples 
 are given by representable cohomology theories
 -- but see also the next section.
 Note however that in the case of algebraic K-theory,
  they should be obtained without using representability:
  the case of absolute purity is of course
   the localization theorem of Quillen but the case
   of refined products is less obvious.
\end{rem}

\subsection{\'Etale cohomology}

\begin{ex}
\begin{enumerate}
\item Let $\Lambda$ be a torsion ring of characteristic exponent $N$.
 Let $\base$ be the category of regular schemes
 on which $N$ is invertible.
Then it follows from \cite{SGA4},
 using the method described in Section \ref{sec:abs_coh}
 together with the functoriality of \'etale sheaves
 established in \cite{SGA4},
 that for any closed immersion $i:Z \rightarrow X$,
 the bigraded cohomology groups:
$$
H^n_Z(X_\et,\Lambda(m))
$$
of the twisted sheaf $i^!\Lambda(m)$,
 computed in the small \'etale site of $X$,
 is a ringed cohomology with support over $\base$
 in the sense of the above definition.
 Recall that according to \cite[IX, Th. 3.3]{SGA4},
	this cohomology theory is oriented with an additive
	formal group law.
Moreover, it is absolutely pure over $\reg$
 according to the (absolute purity) theorem of Gabber.

Thus our constructions apply to this cohomology,
 which has an additive formal group law.
 In particular, we get maps
 for projective morphisms of regular $\ZZ[1/N]$-schemes
 on \'etale cohomology with coefficients in
 $\Lambda$.\footnote{This Gysin morphisms
  agree with the one constructed by Gabber-Riou:
	see Remark \ref{rem:compare_Riou}.}
\item Let $l$ be a prime number
 and $\base$ be the category of $\ZZ[1/l]$-regular schemes.

Then we can apply the construction of \cite{Jan}
 to get that \emph{continuous $l$-adic \'etale cohomology with support}:
$$
H^n_{\cont,Z}(X,\ZZ_l(m))
$$
defined in \emph{loc. cit.}, Section 3, (after Remark 3.5)
 is an arithmetic cohomology in the previous sense.
 In fact, homotopy, stability and excision follow from the known results
 of \cite{SGA4}. The refined product can be defined using the
 method of Section 6 given that we have a pairing
$$
\Gamma_T(Z,-) \otimes \Gamma_Z(X,-) \rightarrow \Gamma_T(X,-)
$$
of the functors of global sections with support for torsion sheaves
(as in the above example). Axioms (E1)-(E7) then follow.
 Note this theory is oriented: this is \emph{loc. cit.}
 (3.26). Finally, using the absolute purity theorem of Gabber
 in the form of the computation of $\Gamma_Z(X,-)$
 for $Z \subset X$ regular schemes, we get that
 this cohomology theory with support is absolutely $\reg$-pure.
The same construction works for $\QQ_l$-coefficients.

Thus we can also apply the constructions of this paper
 to $l$-adic \'etale cohomology (integral and rational).
\end{enumerate}
\end{ex}

\begin{num}
Using the more sophisticated theory of \cite{CD3},
 we can get many examples as follows.

Let $\mathscr T$ be a motivic triangulated
 category over $\base$ in the sense of \cite[Def. 2.4.45]{CD3}:
 in other words, this is a category fibered over the category $\base$
  which satisfies the axioms (A1)-(A4) of Par. \ref{num:functoriality_properties_SH}
  together with the homotopy and stability property (in fact, we will not
  use the adjoint properties of \emph{loc. cit.}).
Then we can associate with $\T$ a ringed cohomology theory with support:
 for any closed immersion $i:Z \rightarrow X$ in $\base$, we put:
$$
H^{n,m}_Z(X,\mathscr T):=\Hom_{\T(X)}(i_*(\un_Z),\un_X(m)[n])
$$
where $\un_?$ is the cartesian section of $\T$ made by the unit for the
 tensor product and $\un_X(m)$ denotes the $m$-th Tate twist
 (\cite[2.4.17]{CD3}).
Then exactly the same arguments as in the proof of
 Prop. \ref{prop:ppty_product} shows that
 this theory satisfies axioms (E1)-(E7).

Assume moreover that one has a premotivic adjunction
 (\cite[Def. 1.4.2]{CD3})
$$
\varphi^*:\mathscr T \rightarrow \mathscr \T'
$$
where $\T$ and $\T'$ are motivic categories.
Then, according to \cite[2.3.11 or 2.4.53]{CD3},
 for any closed immersion $i$, $\varphi^*i_* \simeq i_* \varphi^*$
 through a canonical isomorphism (called an \emph{exchange isomorphism}).
Thus, given any closed pair $(X,Z)$ in $\base$,
 one gets by applying $\varphi^*$ a morphism:
$$
H^{n,m}_Z(X,\mathscr T) \xrightarrow{\varphi} H^{n,m}_Z(X,\mathscr T')
$$
which is compatible with contravariant functoriality
 (resp. covariant functoriality, refined product, boundary)
 because $\varphi^*$ commutes with $f^*$ for any morphism $f$
 (resp. $i_*$ for any closed immersion $i$,
  tensor product, localization triangle).
\end{num}

\begin{ex} \label{ex:etale&continuous}
\begin{enumerate}
\item \textit{(Motivic) \'etale cohomology}.--
 Let $R$ be any ring. For any scheme $S$,
 Cisinski and the author have introduced in \cite[5.1.3]{CD4},
 following Voevodsky,
 the category $\DM_\h(S,R)$ of $\h$-motives.
 We proved in \emph{loc. cit.}, Th.~5.6.2, that it forms,
 for various $S$, a motivic triangulated category.
In particular, the cohomology theory
$$
H^{n,m}_{\et,Z}(X,R):=\Hom_{\DM_\h(X,R)}(i_*(\un_Z),\un_X(m)[n])
$$
is a ringed cohomology with support, defined over the category
 of all schemes. According to \emph{loc. cit.}, 5.6.2,
 it is even an arithmetic cohomology.
 According to Voevodsky, this cohomology theory is called
 the \emph{\'etale motivic cohomology} with coefficients
 in $R$.
 In view of the second computation below,
  we think that it should simply be called the \emph{\'etale cohomology}
  with coefficients in $R$.

According to the fundamental results of \emph{loc. cit.},
 one gets for any regular scheme $X$:
\begin{itemize}
\item if $R$ is a $\QQ$-algebra,
$$
H^{n,m}_{\et}(X,R)=H^{n,m}_{\mathcyr B}(X,R)=K_{2m-n}^{(m)}(X)
$$
is Beilinson motivic cohomology;
\item if $R$ is a torsion ring with characteristic exponent
 $N$, for any scheme $X$,
$$
H^{n,m}_{\et}(X,R)=H^{n}_{\et}(X[1/N],R(m)),
$$
where $X[1/N]$ is the open part of $X$ where $N$ is invertible,
 and the right hand side is the usual \'etale cohomology of $X[1/N]$
 with coefficients in $R$ twisted $m$-times.
\end{itemize}
Moreover, for any ring extension $R'/R$, there is
 a premotivic adjunction:
$$
\varphi:\DM_\h(S,R) \rightarrow \DM_\h(S,R')
$$
so that we get a morphism of arithmetic cohomologies:
\begin{equation} \label{eq:coh_morph1}
\varphi:H^{**}_{\et}(-,R) \rightarrow H^{**}_{\et}(-,R').
\end{equation}
\item \textit{Continuous \'etale cohomology}.--
 Let $R$ be any valuation ring with parameter $\ell$.
 Then, according to \emph{loc. cit.}, 7.2.11,
 the homotopy $\ell$-adic completion of $\DM_\h(-,R)$
 gives a motivic triangulated category $\DM_\h(-,\hat R_\ell)$
 and in particular a ringed cohomology theory with support,
 defined over the category of all schemes:
$$
H^{n,m}_{\cont,Z}(X,\hat R_\ell)
:=\Hom_{\DM_\h(X,\hat R_\ell)}(i_*(\un_Z),\un_X(m)[n]).
$$
According to \emph{loc. cit.}, this is an arithmetic cohomology,.

Note that when $R$ is a discrete valuation ring, and $X$ a scheme 
 such that the exponent characteristic of $R/\ell$ is invertible on $X$,
 according to \emph{loc. cit.}, 7.2.21, the triangulated category
 $\DM_\h(X,\hat R_\ell)$ agree with Ekedahl category of $\ell$-adic
 complexes. Thus the cohomology $H^{**}_{\cont,Z}(X,\hat R_\ell)$
 is Jannsen continuous \'etale $\ell$-adic cohomology
 and deserves the name of \emph{continuous
 \'etale cohomology} with coefficients in $\hat R_l$.

From the obvious premotivic adjunction 
 $\hat \rho^*_\ell:\DM_\h(S,R) \rightarrow \DM_\h(S,\hat R_\ell)$
 (see \cite[7.2.4]{CD4}),
 we get a morphism of arithmetic cohomologies:
\begin{equation} \label{eq:coh_morph2}
\rho_\ell:H^{**}_{\et}(-,R) \rightarrow H^{**}_{\cont}(-,\hat R_\ell).
\end{equation}
Let $Q$ be the fraction field of $R$. We now easily get the rational
 version of continuous \'etale cohomology by taking
 tensor product by $Q$ over $R$:
$$
H^{n,m}_{\cont,Z}(X,Q_\ell):=H^{n,m}_{\cont,Z}(X,\hat R_\ell) \otimes_R Q.
$$
which is again an arithmetic cohomology theory.
 And finally, a rational version of the previous morphism:
\begin{equation} \label{eq:coh_morph3}
\rho_\ell:H^{**}_{\et}(-,Q) \rightarrow H^{**}_{\cont}(-,Q_\ell).
\end{equation}
\item For a prime number $\ell$,
 combining \eqref{eq:coh_morph1} and \eqref{eq:coh_morph3},
 we get a morphism of ringed cohomologies with support
 on the category of regular $\ZZ[1/\ell]$-schemes:
\begin{equation} \label{eq:coh_morph4}
\rho_\ell:H^{**}_{\mathcyr B}(-) \rightarrow H^{**}_{\cont}(-,\QQ_\ell)
\end{equation}
from Beilinson motivic cohomology to continuous rational $\ell$-adic
 cohomology.
\end{enumerate}
\end{ex}

As a corollary of the preceding examples
 and the constructions of this paper, we thus obtain:
\begin{cor}\label{cor:etale_RR}
Assume $\base$ is one of the following categories of schemes:
\begin{enumerate}
\item[(a)] regular noetherian schemes of finite dimension;
\item[(b)] smooth schemes over a noetherian (singular) scheme
 of finite dimension.
\end{enumerate}
 Let $R$ be a ring (resp. discrete valuation ring with parameter
 $\ell$). In the respective case, we also denote by $Q_l$
 the fraction field of $\hat R_\ell$.
\begin{enumerate}
\item The ring cohomology theory $H^{**}_\et(-,R)$
 (resp. $H^{**}_\cont(-,\hat R_\ell)$, $H^{**}_\cont(-,Q_\ell)$)
 admits Chern classes,  Gysin morphisms
 for any projective morphism of schemes in $\base$,
 and residue morphisms associated with a closed immersion
 $i:Z \rightarrow S$ of schemes in $\base$ which fit into
 the usual localization long exact sequence.
 These residues and Gysin morphisms satisfy the following properties:
 compatibility with transversal pullback, excess of intersection,
 projection formula.
\item The five natural transformations
 of Example \ref{ex:etale&continuous} are functorial with
 respect to Gysin morphisms and localization long exact sequences.
\item For any integer $r\geq 0$,
 there exists a well defined higher Chern character:
$$
\ch_r:K_r(X) \rightarrow \oplus_{n \geq 0} H^{2n-r,n}_\et(Y,Q_\ell)
$$
from Quillen (resp. Thomason-Trobaugh in case (b))
 algebraic K-theory such that for any projective morphism
 $f:Y \rightarrow X$ in $\base$, the following diagram commutes:
$$
\xymatrix@R=16pt@C=22pt{
K_r(Y)_\QQ\ar_{\td(\tau_f).\ch_t}[d]\ar^{f_*}[r]
 & K_r(X)_\QQ\ar^{\ch_t}[d] \\
\oplus_{n \geq 0}H^{2n-r,n}_\et(Y,Q_\ell)\ar^{f_*}[r]
 & \oplus_{n\geq 0} H^{2n-r,n}_\et(X,Q_\ell) \\
}
$$
where on the top line, $f_*$ is the usual covariant
 functoriality of algebraic K-theory.
 Moreover for any closed immersion $i:Z \rightarrow X$
 in $\base$, one gets:
 $$
\xymatrix@R=20pt@C=34pt{
K_r(X-Z)_\QQ\ar_{\quad\td(-N_ZX).\ch_t}[d]\ar^{\partial_{X,Z}}[r]
 & K_r(Z)_\QQ\ar^{\ch_t}[d] \\
\oplus_{n \geq 0}H^{2n-r,n}_\et(Y,Q_\ell)\ar^-{\partial_{X,Z}}[r]
 & \oplus_{n\geq 0} H^{2(n-c)-r+1,n-c}_\et(X,Q_\ell), \\
}
$$
where $N_ZX$ is the normal bundle of $Z$ in $X$.

Under assumption (b), the functor $K_r$ can be replaced 
 by Weibel homotopy invariant K-theory $KH_r$ in the two
 previous diagram.
\end{enumerate}
\end{cor}
As explained in \ref{num:extend},
 Point (1) is a compact form of the results of sections
 \ref{sec:orientation}, \ref{sec:Gysin}
 (recall excess of intersection: \ref{prop:basic_Gysin2},
  projection formula: \ref{prop:basic_Gysin1}(b)).
Point (2) follows from Th. \ref{thm:RR1} and
 Th. \ref{thm:projective_RR}
 because all theories have additive formal group law
 and there is only one strict isomorphism of formal group law:
 this implies the Todd class involved in each formulas
 is necessarily equal to $1$
 (see Section \ref{sec:principle_comput}).
Point (3) finally follows from Prop. \ref{prop:classical_GRR}
 and Prop. \ref{prop:H_classical_GRR}.

\begin{rem} \label{rem:compare_Riou}
When $R$ is a torsion ring with characteristic exponent $N$,
 in \cite{GRiou},
 Riou following a construction of Gabber
 has defined Gysin morphisms on \'etale cohomology of
 $\ZZ[1/N]$-schemes with coefficients in $R$,
 with respect to all lci projective maps between
 any noetherian schemes.
If one restricts to regular schemes,
 we obtain using Th.~\ref{thm:compare_Gysin} that our Gysin
 maps coincide with the construction of Gabber-Riou.

Let us be more precise. First, let us compare our conventions
 with that of \emph{op. cit.}
 Let $X$ be a scheme and $\mathcal E$ be a locally free
 $\mathcal O_X$-module. Then the vector bundle associated with
 $\mathcal E$ is $E=V(\mathcal E^\vee)$, the spectrum over $X$
 of the symmetric algebra induced by the \emph{dual} of $\mathcal E$.
 Because of this convention
 one relates Chern classes used in \cite{GRiou} with ours by
 the formula:
$$
c_r(\mathcal E)=(-1)^r.c_r(E)
$$
(compare with relation \eqref{eq:each_Chern_dual}).

Once this convention is settled, one can apply
 Theorem \ref{thm:compare_Gysin} as Riou proved
 the excess intersection formula in \cite[Prop. 2.3.2]{GRiou}.
 Note also that,
  because of Cor.~\ref{cor:uniqueness_orientation},
	Chern classes are uniquely determined by the choice of
	a stability isomorphism:
$$
H^{2}_\et(\PP^1_S,R(1)) \simeq H^{0}_\et(S,R(0))=R
$$
(which necessarily appears as a particular case
 of the projective bundle formula).
\end{rem}

\bibliographystyle{amsalpha}
\bibliography{RR}

\section*{Acknowledgement}

I want to warmly thank Denis-Charles~Cisinski, Ofer~Gabber, Henri~Gillet, 
 Marc~Levine, Jo\"el~Riou and Alexander~Vishik, for discussions and ideas
 that motivated and made this work possible
 as well as Alberto~Navarro for a very useful proof-reading.
 Last but not least, I will never thank enough the anonymous referee of the
 current version for his incredibly careful reading and his numerous comments.
 In particular, he indicated to me a sign mistake in \cite[Cor. 5.31]{Deg8}
 which is corrected in formula \eqref{eq:Myschenko} of the present paper.
 He helped me to get a truly neater version of this paper.

\end{document}